\documentclass[12pt]{amsart}
\usepackage{amsmath,latexsym,amsfonts,amssymb,amsthm}
\usepackage{geometry}
\usepackage{hyperref}
\usepackage{mathrsfs}
\usepackage{graphicx,color}
\usepackage[alphabetic]{amsrefs}
\usepackage[all,cmtip]{xy}
\usepackage{bbm, stmaryrd}
\usepackage{tabu}

\newtheorem{prop}{Proposition}[section]
\newtheorem{thm}[prop]{Theorem}
\newtheorem{cor}[prop]{Corollary}

\newtheorem{lem}[prop]{Lemma}

\theoremstyle{definition}

\newtheorem{defn}[prop]{Definition}
\newtheorem{expl}[prop]{Example}
\newtheorem{rem}[prop]{\it Remark}

\newtheorem{claim}[prop]{Claim}

\newtheorem*{claim*}{Claim}

\newcommand{\bP}{\mathbb{P}}

\newcommand{\bR}{\mathbb{R}}
\newcommand{\bA}{\mathbb{A}}
\newcommand{\bQ}{\mathbb{Q}}
\newcommand{\bZ}{\mathbb{Z}}
\newcommand{\bN}{\mathbb{N}}

\newcommand{\bG}{\mathbb{G}}
\newcommand{\bT}{\mathbb{T}}
\newcommand{\bk}{\mathbbm{k}}

\newcommand{\tX}{\widetilde{X}}

\newcommand{\tD}{\widetilde{D}}

\newcommand{\cX}{\mathcal{X}}

\newcommand{\cO}{\mathcal{O}}

\newcommand{\cI}{\mathcal{I}}
\newcommand{\cM}{\mathcal{M}}
\newcommand{\cF}{\mathcal{F}}
\newcommand{\cG}{\mathcal{G}}
\newcommand{\cJ}{\mathcal{J}}

\newcommand{\cD}{\mathcal{D}}

\newcommand{\cA}{\mathcal{A}}

\newcommand{\DMR}{\mathcal{DMR}}

\newcommand{\fa}{\mathfrak{a}}
\newcommand{\fb}{\mathfrak{b}}
\newcommand{\fc}{\mathfrak{c}}
\newcommand{\fm}{\mathfrak{m}}

\newcommand{\fab}{\fa_{\bullet}}

\newcommand{\rd}{\mathrm{d}}

\newcommand{\Spec}{\mathrm{Spec}}
\newcommand{\Supp}{\mathrm{Supp}}
\newcommand{\Hom}{\mathrm{Hom}}
\newcommand{\mult}{\mathrm{mult}}
\newcommand{\lct}{\mathrm{lct}}

\newcommand{\Aut}{\mathrm{Aut}}
\newcommand{\vol}{\mathrm{vol}}
\newcommand{\ord}{\mathrm{ord}}
\newcommand{\Val}{\mathrm{Val}}

\newcommand{\Diff}{\mathrm{Diff}}

\newcommand{\gr}{\mathrm{gr}}
\newcommand{\Fut}{\mathrm{Fut}}
\newcommand{\wt}{\mathrm{wt}}
\newcommand{\QM}{\mathrm{QM}}
\newcommand{\Ex}{\mathrm{Ex}}

\newcommand{\triv}{\mathrm{triv}}

\newcommand{\colength}{\mathrm{colength}}

\newcommand{\Proj}{\mathrm{Proj}}

\newcommand{\tC}{\widetilde{C}}
\newcommand{\fX}{\mathfrak{X}}
\newcommand{\cP}{\mathcal{P}}
\newcommand{\bfP}{\mathbf{P}}

\numberwithin{equation}{section}

\title[Finite generation for valuations and K-stability]{Finite generation for valuations computing stability thresholds and applications to K-stability}

\date{}

\author{Yuchen Liu}
\address{Department of Mathematics, Northwestern University, Evanston, IL 60208, USA.}
\email{yuchenl@northwestern.edu}

\author{Chenyang Xu}
\address{Department of Mathematics, Princeton University, Princeton, NJ 08544, USA}
\email     {chenyang@princeton.edu}
\address   {Department of Mathematics, MIT, Cambridge, MA 02139, USA}
\email     {cyxu@math.mit.edu}
\address   {Beijing International Center for Mathematical Research,       Beijing 100871, China}
\email     {cyxu@math.pku.edu.cn}

\author{Ziquan Zhuang}
\address   {Department of Mathematics, MIT, Cambridge, MA 02139, USA}
\email     {ziquan@mit.edu}

\begin{document}

\begin{abstract}
We prove that on any log Fano pair of dimension $n$ whose stability threshold is less than $\frac{n+1}{n}$, any valuation computing the stability threshold has a finitely generated associated graded ring. Together with earlier works, this implies: (a) a log Fano pair is uniformly K-stable (resp. reduced uniformly K-stable) if and only if it is K-stable (resp. K-polystable); (b) the K-moduli spaces are proper and projective; and combining with the previously known equivalence between the existence of K\"ahler-Einstein metric and reduced uniform K-stability proved by the variational approach, (c) the Yau-Tian-Donaldson conjecture holds for general (possibly singular) log Fano pairs.
\end{abstract}

\maketitle

\section{Introduction}

In recent years, the algebro-geometric study of the K-stability of Fano varieties has made remarkable progress. See \cite{Xu-survey} for a comprehensive up-to-date survey.

The theory has naturally driven people's attention to valuations which are not necessarily divisorial. In fact, to further advance the theory, one main question is to show the finite generation property of the associated graded rings for quasi-monomial valuations of higher (rational) rank which minimize functions on the space of valuations arisen from K-stability theory.  

While the finite generation property for divisorial valuations follows from \cite{BCHM}, the higher rank case posts a completely new problem. In fact, there were very few studies on higher rank quasi-monomial valuations from the viewpoint of the minimal model program (MMP), which is our fundamental tool to study K-stability.

In this paper, we prove that quasi-monomial valuations that compute the stability thresholds (or $\delta$-invariants) of log Fano pairs satisfy the finite generation property (see \cite{Xu-HRFG}*{Conjecture 1.2}).  

\begin{thm}[=Theorem \ref{thm:f.g. for delta minimizer}, Higher Rank Finite Generation Conjecture]\label{t-HRFG} 
Let $(X,\Delta)$ be a log Fano pair of dimension $n$ and let $r>0$ be an integer such that $r(K_X+\Delta)$ is Cartier. Assume that $\delta(X,\Delta)<\frac{n+1}{n}$. Then for any valuation $v$ that computes $\delta(X,\Delta)$, the associated graded ring $\gr_v R$, where $R=\bigoplus_{m\in\bZ_{\geq 0}} H^0(X,-mr(K_X+\Delta))$, is finitely generated.
\end{thm}

The assumption $\delta(X,\Delta)< \frac{n+1}{n}$ might look a bit surprising since the original conjecture only assumes $\delta(X,\Delta)\le 1$. However, the improvement becomes quite natural using the  trick of compatible divisors invented in \cite{AZ-index2}.

\subsection{Corollaries of the main theorem}
Together with many earlier works in recent years, Theorem \ref{t-HRFG} solves some central questions in the field of K-stability theory. Combining with \cite{BLZ-opt-destabilization}, the first consequence we have is the following, which says that any log Fano pair that is not uniformly K-stable has an optimal destabilizing degeneration (in terms of preserving the stability threshold).

\begin{thm}[{=Theorem \ref{thm:ODC n+1/n}, Optimal Destabilization Conjecture}]\label{t-ODC}
Let $(X,\Delta)$ be a log Fano pair of dimension $n$ such that $\delta(X,\Delta)<\frac{n+1}{n}$, then $\delta(X,\Delta)\in \bQ$ and there exists a divisorial valuation $E$ over $X$ such that
$\delta(X,\Delta)=\frac{A_{X,\Delta}(E)}{S_{X,\Delta}(E)}$. 

In particular, if $\delta(X,\Delta)\le 1$, then there exists a non-trivial special test configuration $(\cX,\Delta_{\cX})$ with a central fiber $(X_0,\Delta_0)$ such that $\delta(X,\Delta)=\delta(X_0,\Delta_0)$, and $\delta(X_0,\Delta_0)$ is computed by the $\mathbb{G}_m$-action induced by the test configuration structure.
\end{thm}

The second main application of Theorem \ref{t-HRFG} is on the general construction of the K-moduli space. Indeed, it is proved in \cite{BHLLX-theta} that Theorem \ref{t-ODC} implies that there is a $\Theta$-stratification on the stack $\mathcal M^{\rm Fano}_{n,V,C}$ of $\mathbb{Q}$-Gorenstein families of $n$-dimensional log Fano pairs $(X,\Delta)\to B$ with a fixed volume $V$ and $C\cdot \Delta$ being integral. By the general theory of $\Theta$-stratification \cite{AHLH-S&theta}, this yields the properness of the K-moduli space. We also conclude the projectivity following \cite{XZ-CM-positive}.

\begin{thm}[{=Theorem \ref{thm:K-moduli conj}, Properness and Projectivity of K-moduli spaces}]\label{t-kmoduli}
The K-moduli space $M^{\rm Kps}_{n,V,C}$ is proper, and the CM line bundle on $M^{\rm Kps}_{n,V,C}$ is ample. 
\end{thm}

\begin{rem}
When the base field $\bk=\mathbb{C}$ and we restrict ourselves to the component whose fibers parametrize $\bQ$-Gorenstein smoothable K-polystable Fano varieties, then the properness follows from the analytic work of \cites{DS-GHlimit,CDS15, Tia15} (see e.g. \cite{LWXI}).
The ampleness of the CM line bundle restricting to this component is a consequence of \cite{XZ-CM-positive} combining with the analytic result in \cite{Li19}.
\end{rem}

\begin{rem}
Theorem \ref{t-kmoduli} is the last step in the general purely algebraic construction of the K-moduli space. We briefly review the previously known steps here: 

A notion of a family of log pairs over a general base was introduced in \cite{Kollar-pair}. \cite{Jiang-boundedness} (which heavily relies on  \cites{Birkar-boundedcomplements,Birkar-BAB}) proved the boundedness of K-semistable Fano varieties with a fix volume. A different proof, which only uses the solution to Batyrev's conjecture \cite{HMX-ACC}, is given in \cite{XZ-uniqueness}. Then using boundedness of complements \cite{Birkar-boundedcomplements}, \cite{BLX-openness} and \cite{Xu-quasimonomial} gave two different proofs for the openness of K-semistability, and as a consequence the open subfunctor $\mathcal M^{\rm Kss}_{n,V,C}\subseteq \mathcal M^{\rm Fano}_{n,V,C}$ parametrizing K-semistable log Fano pairs yields an Artin stack of finite type, which is called {\it the K-moduli stack}. 

By \cite{LWX-cone}, \cite{BX-uniqueness} and \cite{ABHLX-goodmoduli}, we know that $\mathcal M^{\rm Kss}_{n,V,C}$ admits a separated good moduli space $M^{\rm Kps}_{n,V,C}$ which parametrizes K-polystable log Fano pairs and is called {\it the K-moduli space}.

The remaining part is the properness and projectivity of the K-moduli space. In \cite{BHLLX-theta}, it is proved that Theorem \ref{t-ODC} implies the properness of $M^{\rm Kps}_{n,V,C}$; and in \cite{XZ-CM-positive}, it is shown that Theorem \ref{t-HRFG} implies the ampleness of the CM line bundle (introduced in \cite{Tian-97, PT-CM}) on $M^{\rm Kps}_{n,V,C}$.
\end{rem}

The next major consequence of Theorem \ref{t-HRFG} is the complete solution of the Yau-Tian-Donaldson conjecture for log Fano pairs, including singular ones.

\begin{thm}[{=Theorem \ref{thm:YTD}, Yau-Tian-Donaldson Conjecture}] A log Fano pair $(X,\Delta)$ is uniformly K-stable $($resp. reduced uniformly K-stable$)$ if and only if it is K-stable $($resp. K-polystable$)$. In particular, when the base field $\bk=\mathbb C$, $(X,\Delta)$ admits a weak K\"ahler-Einstein $($KE$)$ metric if and only if it is K-polystable. 
\end{thm}

\begin{rem}
In this generality, the direction which says that the existence of KE metrics implies K-polystability was settled in \cite{Berman-Kpoly}. 

For the converse direction, \cite{BBJ-variational} initiated a variational approach to the Yau-Tian-Donaldson conjecture, and the analytic side of this approach was completed in \cites{Li-Tian-Wang, Li19}, which shows that a log Fano pair $(X,\Delta)$ admits a weak KE metric if and only if it is {\it reduced uniformly K-stable}. Therefore, what remains to show is the purely algebro-geometric statement that K-polystability is equivalent to reduced uniform K-stability. 

When $|{\rm Aut}(X,\Delta)|<\infty$, this means for any log Fano pair $(X,\Delta)$ which is not {\it uniformly} K-stable, i.e. $\delta(X,\Delta)\le 1$, we need to show it is also not K-stable by producing a degeneration that destabilizes the log Fano pair. In \cite{BLX-openness} (see also \cite{BLZ-opt-destabilization}), a step toward constructing such a degeneration was made. More precisely, it was shown that the valuation which computes $\delta(X,\Delta) (\le 1)$ is quasi-monomial, and the degeneration should be the ${\rm Proj}$ of the graded ring associated to this valuation, provided it is finitely generated. Theorem \ref{t-HRFG} addresses the finite generation and then it follows that the sought degeneration as in Theorem \ref{t-ODC} exists. As a consequence, it establishes the equivalence between K-stability and uniform K-stability. In the more general case when the automorphism group is positive dimensional, it is shown in \cite{XZ-CM-positive} that Theorem \ref{t-HRFG} implies the equivalence between K-polystability and {\it reduced uniform} K-stability, by refining the argument from \cite{BLX-openness} and finding a quasi-monomial valuation that is not induced by a one parameter subgroup of $\Aut(X,\Delta)$, which yields a non-product type degeneration. 

\smallskip

We note that when $X$ is smooth and $\Delta=0$, the above theorem was first proved in \cites{CDS15, Tia15} using Cheeger-Colding-Tian theory, which seems difficult to generalize to the case of general (possibly singular) log Fano pairs. 
\end{rem}

We also prove the following statement, which is a (necessarily) modified version of a conjecture first raised by Donaldson in \cite{Don-cone}*{Conjecture 1} (see also \cite{Sze-conic} and \cite{BL18}*{Section 7} for some further discussions of the problem). 

\begin{thm}[{=Theorem \ref{t-DonConj}}]
Let $(X,\Delta)$ be a log Fano pair such that $\delta:=\delta(X,\Delta)<1$. For any sufficiently divisible integer $m>0$ and any general member $D_0$ of the linear system $|-m(K_X+\Delta)|$, if we take $D=\frac{1}{m}D_0$, then the pair $(X,\Delta+(1-\delta)D)$ is K-semistable. In particular, $(X,\Delta+(1-\delta')D)$ is uniformly K-stable for any $0\le \delta'<\delta$. 
\end{thm}
For a smooth Fano manifold $X$, the above theorem was essentially implied by a combination of \cite{Z-equivariant} with analytic results in \cite{CDS15, Tia15}. 
We also note that, as one can easily see from our proof of Theorem \ref{t-DonConj}, the integer $m$ can be chosen uniformly for any bounded family $(X,\Delta)\to S$ of log Fano pairs, e.g. the family of all smooth Fano manifolds with a fixed dimension. 



\subsection{Outline of the proof of Theorem \ref{t-HRFG}}

Recall that by \cite{BLX-openness}, any valuation computing $\delta(X,\Delta)\le 1$ is an lc place of a $\bQ$-complement $\Gamma$, and every divisorial lc place $w$ of the complement (parametrized by the rational points of the dual complex $\DMR(X,\Delta+\Gamma)$) induces a weakly special degeneration of the log Fano pair. In addition, if finite generation holds for a quasi-monomial valuation $v$, it is observed in \cite{LX-higher-rank} (see also Lemma \ref{l-finitegeneratedimpliesisom}) that as the valuation $w$ in the minimal rational affine space containing $v$ gets sufficiently close to $v$, the central fibers of the induced degenerations would be isomorphic to each other and in particular are bounded. Our first observation is that the converse of this implication is true, namely if the central fibers of the degenerations induced by nearby rational valuations belong to a bounded family, then indeed the associated graded ring is finitely generated. This is achieved in Theorem \ref{thm:f.g. criterion}. In fact, we consider the expected vanishing order function $w\to S(w)$ on a minimal rational affine space containing $v$. It is easy to show this function is concave, and the boundedness assumption implies that there exists a uniform $M$ such that $MS(w)$ is integer, if $w$ is an integral point. An elementary Diophantine approximation argument implies that $S$ is linear in a neighborhood of $v$. Given two arbitrary valuations $v_0$ and $v_1$, there is another natural construction (which is called the geodesic ray) of filtrations $\cF_t$ connecting $\cF_{v_0}$ and $\cF_{v_1}$. It always satisfies $S(\cF_t)=(1-t)S(v_0)+tS(v_1)$. If $v_0$ and $v_1$ are two valuations on a simplex, then for the natural valuation $v_t$ on the ray connecting $v_0$ and $v_1$ on the simplex, $\cF^{\lambda}_t R\subseteq \cF_{v_t}^{\lambda}R$ for any $\lambda$. This implies in general $S(\cF_t)\le S(\cF_{v_t})$. So our assumption of the boundedness of the degeneration implies that $S(\cF_t)=S(\cF_{v_t})$ as $S(v_t)$ is linear on $t$. When $t\in \mathbb Q$ and $v_0, v_1$ are lc places of a fixed $\mathbb Q$-complement, one can show that $\cF_t$ arises from a weakly special test configuration. Using the explicit description of filtrations that correspond to weakly special test configurations in \cite{BLX-openness}, it is then not hard to deduce from $S(\cF_t)=S(\cF_{v_t})$ that $\cF_t=\cF_{v_t}$ for all $t$, i.e. the filtration induced by $v_t$ is the same as $\cF_t$ which is finitely generated.

\smallskip

Therefore, in order to prove Theorem \ref{t-HRFG}, the remaining main technical goal is to construct a specific $\mathbb{Q}$-complement $\Gamma$ with the given minimizer $v$ as an lc place, and show that in a neighborhood of $v$ in the rational affine subspace in the dual complex $\DMR(X,\Delta+\Gamma)$, the degenerations corresponding to the rational points have bounded central fibers. The construction of the complement $\Gamma$ is indeed tricky, and the argument takes several steps.

In the first step, using compatible basis type divisors as first introduced in \cite{AZ-index2}, we prove an improvement of \cite{BLX-openness} that the complement $\Gamma$ can be chosen to contain some fixed multiple of any effective divisor $D\sim_{\mathbb Q} -(K_X+\Delta)$. This allows us to extend the bound in our assumption from $\delta(X,\Delta)\le 1$ to $\delta(X,\Delta)<\frac{n+1}{n}$. Another improvement of \cite{BLX-openness} is that we prove any quasi-monomial $v$ that computes $\delta(X,\Delta)<\frac{n+1}{n}$ is an lc place of a $\mathbb Q$-complement, using the global ACC Conjecture proved in \cite{HMX-ACC} (but not the boundedness of complements proved in \cite{Birkar-boundedcomplements}).


Next, starting with a fixed log resolution $\pi\colon (Y,E)\to (X,\Delta)$ such that $v\in \QM(Y,E)$, we run the construction in the previous step to get a $\mathbb{Q}$-complement $\Gamma\sim_{\mathbb Q}-(K_X+\Delta)$ containing a multiple of the pushforward of some generally positioned ample $\mathbb{Q}$-divisor $G$ on $Y$. The key idea is to look at the intersection $\QM(Y,E)\cap \DMR(X,\Delta+\Gamma)$. It consists of valuations that we call \emph{monomial lc places of special $\bQ$-complements} (with respect to the log resolution $(Y,E)$), see Definition \ref{d-specialcomplement} for the precise definition. The advantages of the valuations in the above intersections are twofold: firstly, rational points in this new dual complex corresponds simply to weighted blowups on the fixed log smooth model $Y$, which are  easier to analyze; secondly and more importantly, any rational point $w\in \QM(Y,E)\cap \DMR(X,\Delta+\Gamma)$ corresponds to special divisors $F_w$, i.e. they induce klt degenerations of the log Fano pair (this follows from a tie-breaking argument, see Lemma \ref{lem:special DMR}). For general complements the second property is hard to come by. As it turns out in the end (see Theorem \ref{thm:f.g. for special valuation}), a valuation satisfies the finite generation property and induces a klt degeneration of the log Fano pair if and only if it is a monomial lc place of some special $\bQ$-complement.

\smallskip

Recall that we have reduced the finite generation property of $v$ to the boundedness of the central fibers of the degenerations induced by nearby rational valuations in $\QM(Y,E)\cap \DMR(X,\Delta+\Gamma)$. Given that the degenerations have klt central fibers in this case, we can invoke results from \cites{Jiang-boundedness, XZ-uniqueness} to prove the boundedness once we establish a uniform positive lower bound on the $\alpha$-invariants of the central fibers.

The $\alpha$-invariants of the central fibers can be detected on $(X,\Delta)$ without writing down the explicit degeneration, see Lemma \ref{lem:(alpha,v)-complement}. For our purpose, it suffices to show that there exists a uniform $\alpha>0$ such that for any $w$ sufficiently close to $v$, and any effective $\mathbb Q$-divisor $D\sim_{\mathbb Q}-K_X-\Delta$, we can find another effective $\mathbb Q$-divisor $D'$ such that $(X,\Delta+\alpha D+(1-\alpha)D')$ is log canonical and $w$ is an lc place. This guarantees that the $\alpha$-invariants of central fibers are bounded from below by $\alpha$.

Let us consider the special case where $w(D)=A_{X,\Delta}(w)$ (it is not hard to see that this implies the general case, see the proof of Lemma \ref{lem:alpha bound degeneration}). In this case, we necessarily have $w(D')=A_{X,\Delta}(w)$ in order for $w$ to be an lc place. For simplicity, let us also assume that $w$ is centered at a closed point on $X$. By \cite{BCHM}, there exists a plt blowup $\rho\colon X_w\to X$ which extracts the Koll\'ar component $E_w$ that corresponds to the divisorial valuation $w$. If the complementary divisor $D'$ as above exists, then by adjunction the pair $(E_w,\Delta_w+\alpha D_w+(1-\alpha)D'_w)$ has to be log canonical (here the divisors $\Delta_w=\Diff_{E_w}(\rho^{-1}_*\Delta)$ and $D_w=\rho^{-1}_*D|_{E_w}$ etc. are what we naturally get by doing adjunction along $E_w$). In particular, $(E_w,\Delta_w+\alpha D_w)$ needs to be lc. Therefore, we should first guarantee that the $\alpha$-invariant of the Koll\'ar component $(E_w,\Delta_w)$ is at least bounded from below from $\alpha$. 

This is still not enough as the complementary divisor $D'$ adds another component $(1-\alpha)D'_w$ to the pair $(E_w,\Delta_w+\alpha D_w)$, potentially making its singularities worse. However, this is not a serious problem if a significant proportion of $D'_w$ comes from a basepoint-free linear system. The easiest way to make this happen is to impose the condition that $-\rho^*(K_X+\Delta)-\varepsilon A_{X,\Delta}(E_w) E_w$ is nef for some fixed constant $\varepsilon>0$. If $D''$ is the pushforward of a general divisor in the $\bQ$-linear system $|-\rho^*(K_X+\Delta)-\varepsilon A_{X,\Delta}(E_w) E_w|_\bQ$, then $D''_w$ has no contribution to the singularities. With some extra work, we can absorb $D''$ into $D'$ while keeping the singularities mild. For more details, see Lemma \ref{lem:alpha bound degeneration}.

\smallskip

With the above discussions, we are reduced to showing that the $\alpha$-invariant of the Koll\'ar component and the nef threshold as above have uniform lower bounds that are independent of $w$. This is not hard to see if $E_w$ is a weighted blowup at a smooth point of $X$, as the Koll\'ar component $E_w$ is simply a weighted projective space and both the $\alpha$-invariant and the nef threshold can be controlled explicitly. In general, we don't know how the map $\rho\colon X_w\to X$ or the Koll\'ar component $E_w$ look like, but at least we have a log smooth model $(Y,E)$ and $E_w$ can be extracted by a weighted blowup on $Y$. The last key idea of the proof is to transfer all the local lower bound statement to the log resolution $(Y,E)$, where things are easier to verify using the explicit geometry of the weighted blowup, and then descend to the original log Fano pair $(X,\Delta)$ using the Koll\'ar-Shokurov connectedness theorem. These are the final technical steps of the proof, see Lemma \ref{lem:alpha bound plt blowup} (for the lower bounds of $\alpha$-invariants of divisors over $Y$), Lemma \ref{lem:seshadri bound} (for the nef thresholds) and Lemma \ref{lem:alpha bound degeneration} (for the descent to $X$).

\medskip

\noindent {\bf Acknowledgement}: We would like to thank Harold Blum, Xiaowei Wang and Chuyu Zhou for helpful discussions. We also would like to thank the anonymous referees for many helpful comments. 
YL is partially supported by NSF Grant DMS-2148266 (formerly DMS-2001317).
CX is partially supported by NSF Grant DMS-1901849 and DMS-1952531. 
ZZ is partially supported by NSF Grant DMS-2055531.

\section{Preliminaries}\label{s-prelim}

Throughout this paper, we work over an algebraically closed field $\bk$ of characteristic $0$. We follow the standard terminology from \cites{KM98, Kol13}.

\begin{defn}\label{d-logsmooth}
A \emph{pair} $(X,\Delta)$ is a normal variety $X$ together with an effective $\bQ$-divisor $\Delta$ such that $K_X+\Delta$ is $\bQ$-Cartier.  A \emph{log Fano pair} $(X, \Delta)$ is a pair such that $X$ is proper, $-K_X-\Delta$ is ample, and $(X,\Delta)$ is klt. When $\Delta=0$, $X$ is also called a \emph{$\mathbb Q$-Fano variety}. A normal variety $X$ is said \emph{of Fano type}, if there exists an effective $\mathbb{Q}$-divisor $\Delta$ such that $(X,\Delta)$ is a log Fano pair. 

A \emph{log smooth model} $(Y,E)$ over a pair $(X,\Delta)$ consists of a log resolution $\pi: Y\to (X,\Delta)$  and a reduced divisor $E$ on $Y$, such that $E+\Ex(\pi)+\pi_*^{-1}\Delta$ has simple normal crossing (SNC) support. 
\end{defn}

\subsection{Valuations}

In this subsection, we assume that $X$ is a normal variety.

\begin{defn}
 A \emph{valuation} $v$ on $X$ is a $\bR$-valued valuation $v: K(X)^{\times}\to \bR$ such that $v$ has a center on $X$ and $v|_{\bk^\times}=0$.  By convention, we set $v(0)=+\infty$. We denote by $\Val_X$ the set of all valuations on $X$. Recall that the \emph{center} of $v$ on $X$, denoted by $c_X(v)$, is a scheme-theoretic point $\zeta\in X$ such that $v\geq 0$ on $\cO_{X,\zeta}$, and $v>0$ on the maximal ideal $\fm_{X,\zeta}$. Denote by $C_X(v):=\overline{c_X(v)}$. Since $X$ is separated, a center of $v$ on $X$ is unique if exists. If, in addition, $X$ is proper, then every valuation $v$ has a center on $X$. The \emph{trivial valuation} $v_{\triv}$ is defined as $v_{\triv}(f)=0$ for any $f\in K(X)^{\times}$.
 
 For a valuation $v\in \Val_X$, we define its \emph{valuation ideal sheaf} $\fa_p(v)$ for $p\in \bR_{\geq 0}$ as 
 \[
 \fa_p(v):=\{f\in \cO_X\mid v(f)\geq p\}. 
 \]
 We also define the \emph{valuation ideal sequence} of $v$ as $\fa_\bullet(v):=(\fa_m(v))_{m\in\bZ_{\geq 0}}$.
 
 For a section $s$ of a line bundle $L$ on $X$ and a valuation $v$ on $X$, we fix an isomorphism $\phi\colon L|_{U}\cong \mathcal{O}_U$ for some open set $U$ containing $c_X(v)$. We define $v(s):=v(f)$, where $f=\phi(s|_U)\in \Gamma(\mathcal{O}_U)$. It is clear that this does not depend on the choice of the trivialization $\phi$ or the open set $U$.
 
 For a valuation $v\in \Val_X$, we define its \emph{valuation semigroup} $\Phi^+_v:=\{v(f)\mid f\in \cO_{X,c_X(v)}\setminus \{0\}\}$ and its \emph{valuation group} $\Phi_v:=\{v(f)\mid f\in K(X)^\times\}$. The \emph{rational rank} of $v$ is defined as $\mathrm{rat.rk}(v):=\mathrm{rank}_{\bZ}\Phi_v$.
\end{defn}

\begin{defn}
 Let $\pi:Y\to X$ be a proper birational morphism where $Y$ is normal. A prime divisor $E$ on $Y$ is called a \emph{prime divisor over} $X$. It induces a valuation $\ord_E: K(X)^{\times}\to \bZ$ by taking the vanishing order along $E$. A valuation $v\in \Val_X$ is called \emph{divisorial} if $v=c\cdot \ord_E$ for some prime divisor $E$ over $X$ and some $c\in \bR_{\ge 0}$.
\end{defn}

\begin{defn} \label{d-quasimonomial}
 Let $\pi: Y\to X$ be a birational morphism where $Y$ is normal. Let $\eta\in Y$ be a scheme-theoretic point such that $Y$ is regular at $\eta$. For a regular system of parameters  $(y_1,\cdots, y_r)$ of $\cO_{Y,\eta}$ and $\alpha\in \bR_{\geq 0}^r $, we define a valuation $v_\alpha$ as follows. For $f\in \cO_{Y,\eta}\setminus \{0\}$, we may write $f$ in $\widehat{\cO_{Y,\eta}}\cong \kappa(\eta)\llbracket y_1, \cdots, y_r\rrbracket$ as $f=\sum_{\beta\in \bZ_{\geq 0}^r} c_\beta y^\beta$, where $c_\beta\in \kappa(\eta)$ and $y^\beta=y_1^{\beta_1}\cdots y_r^{\beta_r}$ with $\beta=(\beta_1,\cdots, \beta_r)$. We set
 \[
 v_\alpha (f):= \min \{\langle \alpha,\beta\rangle\mid c_{\beta}\neq 0\}.
 \]
 A valuation $v\in \Val_X$ is called \emph{quasi-monomial} if $v=v_\alpha$ as above for some $\pi:Y\to X$, $\eta$, $(y_1,\cdots, y_r)$ and $\alpha$. It is proven in \cite{ELS03} that a valuation $v$ is quasi-monomial if and only if it is an Abhyankar valuation, i.e. $v$ satisfies $\mathrm{tr.deg}(v)+\mathrm{rat.rk}(v)=\dim X$ where $\mathrm{tr.deg}(v)$ is the transcendental degree of $v$. From the above definition, we have that for any $f\in \cO_{Y,\eta}\setminus \{0\}$, the function $\alpha\mapsto v_\alpha (f)$ is piecewise rational linear and is concave, i.e.
 \[
 v_{t\alpha_1+(1-t)\alpha_2}(f)\ge t\cdot v_{\alpha_1}(f)+(1-t)\cdot v_{\alpha_2}(f)
 \]
 for any $0\le t\le 1$ and any $\alpha_1,\alpha_2\in \bR^r_{\ge 0}$. In particular, for any non-trivial effective $\bQ$-Cartier divisor $D$ (resp. graded sequence $\fab$ of ideals) on $X$, the function $\alpha\mapsto v_\alpha(D)$ (resp. $\alpha\mapsto v_\alpha(\fab)$) is piecewise rational linear and concave. This fact will be frequently used in this paper.
 
 If, in addition, $\pi:(Y,E=\sum_{i=1}^l E_i)\to X$ is a log smooth model where $(y_i=0)=E_i$ for $1\leq i\leq r$ as an irreducible component of $E$, then we denote the set $\{v_{\alpha}\mid \alpha \in \bR_{\geq 0}^r\}$ by $\QM_{\eta}(Y,E)$. We also set $\QM(Y,E):=\cup_\eta \QM_{\eta}(Y,E)$ where $\eta$ runs through all generic points of $\cap_{i\in J}E_i$ for some non-empty subset $J\subseteq \{1,\cdots, l\}$. For later reference, we note that if $v$ is a quasi-monomial valuation and $q$ is its rational rank, then the log resolution $\pi\colon Y\to X$ can be chosen (by passing to a further blowup) such that $v\in \QM_\eta (Y,E)$ for some codimension $q$ point $\eta$.
\end{defn}

\subsection{Stability thresholds}

We first define the log discrepancy function for valuations.
\begin{defn}[\cites{JM12, BdFFU15}] \label{d-log discrepancy}
For a pair $(X,\Delta)$, we define the \emph{log discrepancy} function $A_{X,\Delta}:\Val_X\to \bR\cup \{+\infty\}$ as follows.
\begin{enumerate}
    \item If $v=c\cdot \ord_E$ is divisorial, then 
    \[
    A_{X,\Delta}(v):=c\cdot A_{X,\Delta}(E)= c(1+ \mathrm{coeff}_E(K_Y-\pi^*(K_X+\Delta))).
    \]
    \item If $v=v_\alpha$ is quasi-monomial for a log smooth model $(Y,E)$ over $(X,\Delta)$, then 
    \[
    A_{X,\Delta}(v):=\sum_{i=1}^r \alpha_i\cdot A_{X,\Delta}(E_i).
    \]
    It is clear that $A_{X,\Delta}$ is linear on $\QM_{\eta}(Y,E)$.
    \item According to \cite{JM12}, there is a retraction map $r_{Y,E}: \Val_X\to \QM(Y,E)$ for any log smooth model $\pi:(Y,E)\to (X,\Delta)$ satisfying $\Supp (\Ex(\pi)+\pi_*^{-1}\Delta)\subseteq E$. For any $v\in \Val_X$, we define 
    \[
    A_{X,\Delta}(v):=\sup\{A_{X,\Delta}(r_{Y,E}(v))\mid (Y,E)\textrm{ is a log smooth model over }(X,\Delta)\}.
    \]
\end{enumerate}
\end{defn}

From the definition, we know that $A_{X,\Delta}(\lambda v)=\lambda\cdot A_{X,\Delta}(v)$ for any $\lambda\in \bR_{\geq 0}$. Note that a pair $(X,\Delta)$ is lc (resp. klt) if and only if $A_{X,\Delta}(v)\geq 0$ (resp. $>0$) for all valuations $v\in \Val_X\setminus \{v_{\triv}\}$. We also set
\[
\Val_X^\circ:= \{v\in \Val_X\mid v\neq v_{\triv} \textrm{ and } A_{X,\Delta}(v)<+\infty\}.
\]
Then it is clear that $\Val_X^\circ$ contains all non-trivial quasi-monomial valuations on $X$. If $(X,\Delta)$ is lc, then $v\in \Val_X$ is  an \emph{lc place} of $(X,\Delta)$ if $A_{X,\Delta}(v)=0$. If $(Y,E)$ is a log smooth model over an lc pair $(X,\Delta)$ satisfying $\Supp (\Ex(\pi)+\pi_*^{-1}\Delta)\subseteq E$, then by \cite{JM12}*{Corollary 5.4} we know that the set of all lc places of $(X,\Delta)$ coincides with $\QM(Y, E')$ where $E'$ is the sum of irreducible components $E_i$ of $E$ satisfying $A_{X,\Delta}(E_i)=0$. In particular, any lc place of $(X,\Delta)$ is a quasi-monomial valuation in $\QM(Y,E)$.

In the rest of this subsection, we assume that $(X,\Delta)$ is a log Fano pair. Let $r$ be a positive integer such that $L:=-r(K_X+\Delta)$ is  Cartier. Then the \emph{section ring} of $(X,L)$ is given by \[
R(X,L):=R=\bigoplus_{m\in\bZ_{\geq 0}} R_m =\bigoplus_{m\in\bZ_{\geq 0}} H^0(X, \cO_X(mL)).
\]

\begin{defn}[\cite{FO18}]
 Let $m$ be a positive integer such that $N_m:= h^0(X, \cO_X(mL))>0$. An \emph{$m$-basis type divisor} $D$ on $X$  is a divisor of the following form
 \[
 D=\frac{1}{mrN_m}\sum_{i=1}^{N_m} (s_i=0),
 \]
 where $(s_1,\cdots, s_{N_m})$ is a basis of the vector space $R_m=H^0(X, \cO_X(mL))$. It is clear that $D\sim_\bQ -(K_X+\Delta)$.
\end{defn}

\begin{defn}[\cite{BJ-delta}]
 Let $v\in \Val_X^\circ$ be a valuation. Let $m$ be a positive integer. We define the invariants 
 \begin{align*}
     T_m(v)&:= \max \{\tfrac{1}{mr}v(s)\mid s\in R_m\setminus\{0\}\},\\
     S_m(v)&:=\max \{v(D)\mid \textrm{$D$ is of $m$-basis type}\}.
 \end{align*}
 We define the \emph{$T$-invariant} and \emph{$S$-invariant} of $v$ as
 \[
 T_{X,\Delta}(v):= \sup_{m\in \bZ_{>0}} T_m(v)\quad \textrm{and}\quad S_{X,\Delta}(v):=\lim_{m\to\infty} S_m(v). 
 \]
 Note that the above limit exists as finite real numbers by \cite{BJ-delta}*{Corollary 3.6}.
 

\end{defn}

\begin{defn}[\cite{FO18, BJ-delta}]
 Let $m$ be a positive integer. Then we define 
 \[
 \delta_m(X,\Delta):= \inf \{\lct(X,\Delta;D)\mid \textrm{$D$ is of $m$-basis type}\}.
 \]
 The above infimum is indeed a minimum since $m$-basis type divisors are bounded, and lct takes finitely many values on a bounded $\bQ$-Gorenstein family \cite{Amb16}*{Corollary 2.10}. In particular, there exists some divisor $E$ over $X$ such that $\delta_m(X,\Delta)=\frac{A_{X,\Delta}(E)}{S_m(E)}$.

 The \emph{stability threshold} (also called \emph{$\delta$-invariant}) of a log Fano pair $(X,\Delta)$ is defined as
 \[
 \delta(X,\Delta):=\lim_{m\to\infty} \delta_m(X,\Delta).
 \]
   Equivalently, we have 
 \[
 \delta(X,\Delta)= \inf_{v\in\Val_X^\circ}\frac{A_{X,\Delta}(v)}{S_{X,\Delta}(v)}.
 \]
 We say that $\delta(X,\Delta)$  is \emph{computed} by a valuation $v\in \Val_X^\circ$  if $\delta(X,\Delta)=\frac{A_{X,\Delta}(v)}{S_{X,\Delta}(v)}$.
\end{defn}

\begin{thm}[{Fujita-Li valuative criterion, see \cite{Fujita-criterion, Li-criterion, BX-uniqueness}}]\label{t-fujitali}
Let $(X,\Delta)$ be a log Fano pair. Then $(X,\Delta)$ is K-semistable $($resp. uniformly K-stable$)$ if and only if $\delta(X,\Delta)\ge 1$ $($resp. $\delta(X,\Delta)>1)$; and $(X,\Delta)$ is K-stable if and only if $A_{X,\Delta}(E)>S_{X,\Delta}(E)$ for all divisorial valuations $E$.
\end{thm}

\begin{rem}
The original definitions of K-stability notions of a log Fano pair $(X,\Delta)$ use test configurations and Futaki invariants (see Definitions \ref{d-testconfiguration} and \ref{d-futaki}). In this paper, we will mostly use the equivalent characterization given by Theorem \ref{t-fujitali}. So one can take this as the definitions of corresponding concepts.  See Lemma \ref{l-Fut=beta} for the connection. 
\end{rem}

When ${\rm Aut}(X,\Delta)$ is positive dimensional and $T\subseteq \Aut(X,\Delta)$ is a torus, there is also a reduced version $\delta_T(X,\Delta)$. See \cite{XZ-CM-positive}*{Appendix A} for more discussions.

\begin{defn}[\cite{Tia87, CS08, BJ-delta}]
The \emph{$\alpha$-invariant} of a log Fano pair $(X,\Delta)$ is defined as 
\[
\alpha(X,\Delta):=\inf \{\lct(X,\Delta;D)\mid D\in |-K_X-\Delta|_{\bQ}\}.
\]
Equivalently, we have 
\[
 \alpha(X,\Delta)= \inf_{v\in\Val_X^\circ}\frac{A_{X,\Delta}(v)}{T_{X,\Delta}(v)}.
\]
\end{defn}

\begin{defn}
More generally, for any projective klt pair $(X,\Delta)$ and any effective $\bQ$-Cartier $\bQ$-divisor $G$ on $X$ we define 
\[
\lct(X,\Delta;|G|_\bQ):=\inf \{\lct(X,\Delta;D)\mid D\in |G|_{\bQ}\}.
\]
For any valuation $v\in\Val_X^\circ$, if we let 
\[
T(G;v)=\sup\{\tfrac{1}{m}v(s)\mid m\in\bN\mbox{ sufficiently divisible, } s\in H^0(X, \cO_X(mG))\},
\]
then as above we have 
\[
\lct(X,\Delta;|G|_\bQ)=\inf_{v\in\Val_X^\circ}\frac{A_{X,\Delta}(v)}{T(G;v)}.
\]
\end{defn}

In our argument later, we need the following theorem.
\begin{thm}\label{t-alphabounded}
Fix positive integers $n, C$ and three positive numbers $V,\alpha_0,\delta_0$. If we consider the set $\mathcal{P}$ of all $n$-dimensional log Fano pairs $\{(X,\Delta)\}$ such that $C\cdot \Delta$ is integral, $(-K_X-\Delta)^n=V$ and $\alpha(X,\Delta)\ge \alpha_0$ $($resp. $\delta(X,\Delta)\ge \delta_0)$. Then $\mathcal{P}$ is bounded.
\end{thm}
\begin{proof}
When $\Delta=0$, this is first proved in \cite{Jiang-boundedness} which heavily relies on \cite{Birkar-boundedcomplements,Birkar-BAB}. See also \cite{Chen-boundedness}. Later a proof which only uses the boundedness result from \cite{HMX-ACC} was given in \cite{XZ-uniqueness}. 
\end{proof}

\subsection{Filtrations and compatible basis type divisors}

In this subsection, we assume that $(X,\Delta)$ is a log Fano pair, and $L=-r(K_X+\Delta)$ is an ample Cartier divisor for some $r\in \bZ_{>0}$. Let $R=\oplus_{m\in\bZ_{\geq 0}} H^0(X,\cO_X(mL))$ be the section ring of $(X,L)$. 

\begin{defn}
A \emph{filtration} $\cF$ on $R$ is a collection of vector subspaces $\cF^\lambda R_m \subseteq R_m$ for any $m\in \bZ_{\geq 0}$ and $\lambda\in \bR_{\geq 0}$ satisfying the following properties.
\begin{enumerate}
    \item $\cF^\lambda R_m\subseteq \cF^{\lambda'} R_m$ if $\lambda\geq \lambda'$;
    \item $\cF^\lambda R_m=\cap_{\lambda'<\lambda} \cF^{\lambda'} R_m$ if $\lambda>0$;
    \item $\cF^0 R_m=R_m$ and $\cF^\lambda R_m=0$ for $\lambda\gg 0$;
    \item $\cF^\lambda R_m \cdot \cF^{\lambda'} R_{m'}\subseteq \cF^{\lambda+\lambda'}R_{m+m'}$.
\end{enumerate}

A filtration $\cF$ induces a function $\ord_{\cF}: R_m\to \bR_{\geq 0}$ as $\ord_{\cF}(s):=\max\{\lambda\mid s\in \cF^\lambda R_m\}$. By convention, we set $\ord_{\cF}(0)=+\infty$.

\end{defn}

In this paper, we are mainly interested in the following two types of filtrations coming from valuations or divisors.

\begin{expl}
Any valuation $v\in \Val_X$ induces a filtration $\cF_v$ on $R$  as 
\[
\cF_v^\lambda R_m := \{s\in R_m \mid v(s)\geq \lambda\}.
\]
Any non-zero effective $\bQ$-divisor $G$ on $X$ induces a filtration $\cF_G$ on $R$ as 
\[
\cF_G^\lambda R_m := \{s\in R_m \mid (s=0)\geq \lambda G\}. 
\]
\end{expl}

\begin{defn}
Let $\cF$ be a filtration on $R$. The \emph{associated graded ring} $\gr_{\cF} R$ of $\cF$ is defined as
\[
\gr_{\cF}R:=\bigoplus_{m\in\bZ_{\geq 0}} \bigoplus_{\lambda\in \bR_{\geq 0}} \gr^{\lambda}_{\cF} R_m, \quad \textrm{where } \gr^{\lambda}_{\cF} R_m:= \cF^\lambda R_m/\cup_{\lambda'>\lambda} \cF^{\lambda'} R_m.
\]
We say that $\cF$ is \emph{finitely generated} if $\gr_{\cF}R$ is a finitely generated $\bk$-algebra. For a valuation $v\in \Val_X$, we define the \emph{associated graded ring} of $v$ by $\gr_v R:= \gr_{\cF_v} R$. Note that the grading of $\gr_{v}R$ can be chosen as $(m,\lambda)\in \bZ_{\geq 0}\times \Phi_v^{+}$ where $\Phi_v^+$ is the valuation semigroup of $v$.
\end{defn}

\begin{defn}\label{d-compatible}
Let $\cF$ be a filtration on $R$. A basis $(s_1, \cdots, s_{N_m})$ of $R_m$ is said to be \emph{compatible with} $\cF$ if $\cF^\lambda R_m$ is spanned by some of the $s_i$'s for every $\lambda\in\bR_{\geq 0}$. An $m$-basis type divisor $D=\frac{1}{mrN_m}\sum_{i=1}^{N_m} (s_i=0)$ is said to be \emph{compatible with} $\cF$ if $(s_1, \cdots, s_{N_m})$ is compatible with $\cF$. By abuse of notation, we will say that an $m$-basis type divisor $D$ is compatible with a valuation $v$ (resp. an effective  $\bQ$-divisor $G$) if $D$ is compatible with the filtration induced by $v$ (resp. $G$).
\end{defn}

From the definition, it is easy to see that for any $v\in\Val_X^\circ$, we have $v(D)=S_m(v)$ for any $m$-basis type divisor $D$ that is compatible with $v$. Another useful fact about compatible divisors is the following.

\begin{lem}[\cite{AZ-index2}*{Lemma 3.1}]\label{lem:compatible}
Let $\cF$ and $\cG$ be two filtrations of $R$. Then for any $m\in \bZ_{>0}$ there exists an $m$-basis type divisor that is compatible with both $\cF$ and $\cG$.  
\end{lem}

\begin{defn}
Let $\cF$ be a filtration of $R$. We define the \emph{$T$-invariant} of $\cF$ as
\[
T_{X,\Delta}(\cF):=\sup_{m\in\bZ_{>0}}T_m(\cF)\in [0, +\infty],\quad \textrm{where }T_m(\cF):= \max\{\tfrac{\lambda}{mr}\mid \cF^\lambda R_m\neq 0 \}. 
\]
By Fekete's lemma, we know that $T_{X,\Delta}(\cF)=\lim_{m\to\infty} T_m(\cF)$.
We say that $\cF$ is \emph{linearly bounded} if $T_{X,\Delta}(\cF)<+\infty$. 
\end{defn}

\begin{defn}\label{defn-Sinvariant}
Let $\cF$ be a linearly bounded filtration. Let $(s_1, \cdots, s_{N_m})$ be a basis of $R_m$ that is compatible with $\cF$. We define the invariant 
\[
S_m(\cF):= \frac{1}{mrN_m} \sum_{i=1}^{N_m} \ord_{\cF}(s_i).
\]
It is clear that $S_m(\cF)$ does not depend on the choice of the compatible basis. 
The \emph{$S$-invariant} of $\cF$ is defined as
\[
S_{X,\Delta}(\cF):= \lim_{m\to\infty} S_m(\cF).
\]
Note that the above limit exists as finite real numbers by \cite{BJ-delta}*{Lemma 2.9}. In fact, by \emph{loc. cit.}, we have
\begin{equation} \label{eq:S-inv as integral}
    S_{X,\Delta}(\cF)=\frac{1}{(-K_X-\Delta)^n}\int_0^\infty \vol(\cF^{(t)}R) \rd t
\end{equation}
where $\vol(\cF^{(t)}R):=\lim_{m\to\infty}\frac{\dim \cF^{mrt}R_m}{(mr)^n/n!}$.
\end{defn}

By \cite{BJ-delta}*{Lemma 3.1}, any valuation $v\in \Val_X^\circ$ induces a linearly bounded filtration $\cF_v$. From our definitions, it is easy to see that each invariant from $T_m$, $S_m$, $T_{X,\Delta}$, and $S_{X,\Delta}$ has the same value for $v$ and $\cF_v$. For an effective nonzero $\bQ$-divisor $G$ on $X$, we define $S_m(G):=S_m(\cF_G)$ and $S_{X,\Delta}(G):=S_{X,\Delta}(\cF_G)$. As before, we note that $D\ge S_m(G)\cdot G$ for any $m$-basis type divisor $D$ that is compatible with $G$. The following calculation is also very useful for us.

\begin{lem} \label{lem:S=1/n+1}
Let $\lambda\in\bQ_{>0}$ and let $G\sim_\bQ -\lambda(K_X+\Delta)$ be an effective $\bQ$-divisor. Then $S_{X,\Delta}(G)=\frac{1}{\lambda(n+1)}$ where $n=\dim X$.
\end{lem}

\begin{proof}
We have $\vol(\cF_G^{(t)}R)=\vol(-K_X-\Delta-tG)=(1-\lambda t)^n\vol(-K_X-\Delta)$, thus the result follows from \eqref{eq:S-inv as integral}.
\end{proof}

\subsection{Special divisors and complements}

Let $(X,\Delta)$ be a log Fano pair. We first recall the concepts of (weakly) special test configurations. Note that we omit the polarization in the following definition, because (weakly) special test configurations are naturally anti-canonically polarized.

\begin{defn}\label{d-testconfiguration}
A \emph{weakly special test configuration} $(\cX, \Delta_{\cX})$ of $(X,\Delta)$ consists of the following data:
\begin{itemize}
    \item a normal variety $\cX$ together with a flat proper morphism $\pi: \cX\to \bA^1$;
    \item a $\bG_m$-action on $\cX$ such that $\pi$ is $\bG_m$-equivariant with respect to the standard $\bG_m$-action on $\bA^1$ by multiplication;
    \item $\cX\setminus\cX_0$ is $\bG_m$-equivariantly isomorphic to $X\times (\bA^1\setminus\{0\})$ where the $\bG_m$-action on $X$ is trivial;
    \item an effective $\bQ$-divisor $\Delta_{\cX}$ on $\cX$ such that $\Delta_{\cX}$ is the component-wise closure of $\Delta\times (\bA^1\setminus\{0\})$ under the identification between $\cX\setminus\cX_0$ and $X\times (\bA^1\setminus\{0\})$. 
    \item $-(K_{\cX}+\Delta_{\cX})$ is $\bQ$-Cartier and $\pi$-ample;
    \item $(\cX,\cX_0+\Delta_{\cX})$ is log canonical.
\end{itemize}
A weakly special test configuration $(\cX,\Delta_{\cX})$ of $(X,\Delta)$ is \emph{special} if $(\cX, \cX_0+\Delta_{\cX})$ is plt. The central fiber $(\cX_0, \Delta_0)$ of a special (resp. weakly special) test configuration $(\cX,\Delta_{\cX})$ is called a \emph{special} (resp. \emph{weakly special}) \emph{degeneration} of $(X,\Delta)$. A weakly special test configuration is \emph{trivial} if $\cX$ is $\bG_m$-equivariantly isomorphic to $X\times \bA^1$ where the $\bG_m$-action on $X$ is trivial. 
\end{defn}

\begin{defn}\label{d-futaki}
For a weakly special test configuration $(\cX,\Delta_{\cX})$ of an $n$-dimensional log Fano pair $(X,\Delta)$, we consider its gluing with $(X,\Delta)\times (\mathbb{P}^1\setminus \{0\})\to \mathbb{P}^1\setminus \{0\}$ to get a $\mathbb{G}_m$-equivariant family $\overline{\pi}\colon (\overline{\cX},\Delta_{\overline{\cX}})\to \mathbb{P}^1$. Then we define
$${\rm Fut}(\cX,\Delta_{\cX}):=-\frac{(-K_{\overline{\cX}/\mathbb P^1}-\Delta_{\overline{\cX}})^{n+1}}{(n+1)(-K_X-\Delta)^n}.$$
\end{defn}

It is well known from \cite{LX-special} that, to test K-stability of a log Fano pairs, it suffices to consider weakly special test configurations or even special test configurations, and in this paper we will not need more general test configurations. 

\begin{defn}
We define $(X,\Delta)$ to be K-polystable if and only for any weakly special test configuration $(\cX,\Delta_{\cX})$ of $(X,\Delta)$, we have ${\rm Fut}(\cX,\Delta_{\cX})\ge 0$ and the equality holds if and only if $(\cX,\Delta_{\cX})$ is a product test configuration, i.e. $(\cX,\Delta_{\cX})\cong (X,\Delta)\times\bA^1$.
\end{defn}

\begin{defn} \label{d-specialdivisor}
Let $E$ be a prime divisor over $X$. We say that $E$ is \emph{weakly special} (resp. \emph{special}) over $(X,\Delta)$ if there exists a weakly special (resp. special) test configuration $(\cX,\Delta_{\cX})$ with integral central fiber $\cX_0$, such that the restriction of the valuation $\ord_{\cX_0}$ to the subfield $K(X)\subseteq K(\cX)(=K(X\times \bA^1))$ is equal to $b\cdot \ord_E$ for some $b\in\bZ_{>0}$. By abuse of notation, we will say that $\ord_E$ or any valuation $v$ proportional to $\ord_E$ is {\it weakly special} (resp. {\it special}) if $E$ is weakly special (resp. special).
\end{defn}

\begin{lem}[\cite{Fujita-criterion}*{Theorem 6.13}] \label{l-Fut=beta}
Notation as in Definition \ref{d-specialdivisor}. Then we have 
$${\rm Fut}(\cX,\Delta_{\cX})=b(A_{X,\Delta}(E)-S_{X,\Delta}(E)).$$
\end{lem}

If $E$ is a weakly special divisor over $(X,\Delta)$ such that $\ord_{\cX_0}|_{K(X)}=b\cdot \ord_E$, then the central fiber $(\cX_0,\Delta_0)$ is uniquely determined by $E$ up to isomorphism, and $\cX_0\cong \Proj\,\gr_E R$ (see e.g. \cite[Section 3.6 and Lemma 3.7]{Xu-survey}). 

\begin{defn}\label{d-DMR}
A \emph{$\bQ$-complement} of $(X,\Delta)$ is an effective $\bQ$-Cartier $\bQ$-divisor $D\sim_{\bQ} -K_X-\Delta$ such that $(X, \Delta+D)$ is log canonical. A $\bQ$-complement $D$ is called an \emph{$N$-complement} for $N\in\bZ_{>0}$ if $N(K_X+\Delta+D)\sim 0$, and $N(\Delta+D)\geq N\lfloor\Delta\rfloor + \lfloor(N+1)\{\Delta\}\rfloor$ where $\{\Delta\}:=\Delta-\lfloor\Delta\rfloor$.

For any $\bQ$-complement $D$ of $(X,\Delta)$, we define the \emph{dual complex} of $(X,\Delta+D)$ as 
\[
\DMR(X, \Delta+D):=\{v\in \Val_X^\circ\mid A_{X,\Delta+D}(v)=0\textrm{ and }A_{X,\Delta}(v)=1\}.
\]
(The notation $\DMR$ comes from \cite{dFKX-dualcomplex}.) In particular, the space of all lc places of $(X,\Delta+D)$ is a cone over $\DMR(X,\Delta+D)$. By abuse of notation, we often write $v\in \DMR(X, \Delta+D)$ if $v$ is an lc place of $(X,\Delta+D)$.

A log resolution $(Y,E)$ over $(X,\Delta+D)$ will yield a rational piecewise linear (PL) structure of $\DMR(X, \Delta+D)$ (see \cite{dFKX-dualcomplex}). More generally, for any $v\in \DMR(X,\Delta+D)$ and any log smooth model $(Y,E)$ over $(X,\Delta)$ such that $v$ is contained in the interior of $\QM(Y,E)$, since $A_{X,\Delta}$ is linear on $\QM(Y,E)$, the function $v\mapsto v(D)$ is piecewise rational linear and concave (see Definitions \ref{d-quasimonomial} and \ref{d-log discrepancy}), and $A_{X,\Delta}(v)\ge v(D)$, the intersection 
\[
\QM(Y,E)\cap \DMR(X,\Delta+D)=\{v\in\QM(Y,E)\mid A_{X,\Delta}(v)=1=v(D)\}
\]
is convex and spans a rational linear subset of $\QM(Y,E)$. We call $\QM(Y,E)\cap \DMR(X,\Delta+D)$ {\it the minimal rational PL subspace of $\DMR(X, \Delta+D)$ in $\QM(Y,E)$ containing $v$}.
\end{defn}

As in the previous section, let $R=\bigoplus_{m\in\bN} H^0(X,-mr(K_X+\Delta))$ for some integer $r>0$ such that $r(K_X+\Delta)$ is Cartier. The following fact will be used throughout this paper.

\begin{lem} \label{l-fg for divisorial lc place}
Assume that $v$ is an divisorial lc place of some $\bQ$-complement. Then $\gr_v R$ is finitely generated.
\end{lem}

\begin{proof}
We may assume that $v=\ord_E$ for some prime divisor $E$ over $X$. Let $D$ be a $\bQ$-complement that realizes $E$ as an lc place. Since $(X,\Delta+(1-\varepsilon)D)$ is log Fano and $A_{X,\Delta+(1-\varepsilon)D}(E)<1$ for $0<\varepsilon\ll 1$, by \cite{BCHM}*{Corollary 1.4.3} there exists a projective birational morphism $\pi\colon Y\to X$ that extracts $E$ as the unique exceptional divisor. Moreover, we have $\pi^*(K_X+\Delta+(1-\varepsilon)D)=K_Y+\Gamma$ for some $\Gamma\ge 0$ such that $(Y,\Gamma)$ is klt and $-(K_Y+\Gamma)$ is nef and big, thus $Y$ is of Fano type. It then follows from \cite{BCHM}*{Corollary 1.3.1} that $\bigoplus_{m,k\in\bN} H^0(-mr\pi^*(K_X+\Delta)-kE)$ is finitely generated. As $\gr_v R$ is a quotient of this algebra, we conclude that $\gr_v R$ is also finitely generated.
\end{proof}

\begin{thm}[\cite{BLX-openness}*{Theorems 3.5 and A.2}]\label{t-complements}
There exists $N\in\bZ_{>0}$ depending only on $\dim(X)$ and $\mathrm{Coeff}(\Delta)$ such that the following statements are equivalent  for a prime divisor $E$ over $X$.
\begin{enumerate}
    \item $E$ is a weakly special divisor over $(X,\Delta)$;
    \item $E$ is an lc place of a $\bQ$-complement of $(X,\Delta)$;
    \item $E$ is an lc place of an $N$-complement of $(X,\Delta)$.
\end{enumerate}
\end{thm}
The equivalence (2) and (3) of the above characterization largely relies on the existence of bounded complements established in  \cite{Birkar-boundedcomplements}.

The following observation is made by the third named author (see \cite{Xu-survey}*{Theorem 4.12}).
\begin{thm}\label{thm:zhuang-special}
The following statements are equivalent for a prime divisor $E$ over $X$.
\begin{enumerate}
    \item $E$ is a special divisor over $(X,\Delta)$;
    \item $A_{X,\Delta}(E)<T_{X,\Delta}(E)$ and there exists a $\bQ$-complement $D'$ of $(X,\Delta)$ such that, up to rescaling, $E$ is the only lc place of $(X, \Delta+D')$;
    \item there exists an effective $\bQ$-divisor $D\sim -K_X-\Delta$ and $t\in (0,1)$ such that $(X, \Delta+tD)$ is lc with $E$ as the only lc place $($up to rescaling$)$.
\end{enumerate}
\end{thm}

\section{Log Fano pairs with \texorpdfstring{$\delta(X,\Delta)<\frac{n+1}{n}$}{delta<(n+1)/n}}

In this section, for a valuation $v$ which computes $\delta(X,\Delta)$, we carefully construct a $\mathbb{Q}$-complement $\Gamma$ such that $v$ is an lc place of $(X,\Delta+\Gamma)$. In fact, $\Gamma$ satisfies a number of other technical properties (we call it {\it a special complement}, see Definition \ref{d-specialcomplement} for its definition), which are indispensable for our proof of Theorem \ref{t-HRFG}.

As a by-product, in Section \ref{ss-n+1/n}, we show that the various results in \cite{BLX-openness} can be improved using the construction of compatible basis type divisors introduced in \cite{AZ-index2} (see Definition \ref{d-compatible}). 

\subsection{Complements for higher rank valuations}\label{ss-complement}
Recall that when $\delta(X,\Delta)\le 1$, any valuation computing $\delta(X,\Delta)$ is an lc place of a $\bQ$-complement \cite[Theorem A.7]{BLX-openness}. Using compatible divisors, we first generalize this result to log Fano pairs with $\delta(X,\Delta)<\frac{n+1}{n}$ and investigate the degree of freedom when choosing such complements.
\begin{lem} \label{lem:minimizer's complement}
Let $(X,\Delta)$ be a log Fano pair of dimension $n$ such that $\delta(X,\Delta)=\delta<\frac{n+1}{n}$, and let $v$ be a valuation that computes $\delta(X,\Delta)$. Let $\alpha\in (0,\min\{\frac{\delta}{n+1},1-\frac{n\delta}{n+1}\})\cap \bQ$. Then for any effective divisor $D\sim_\bQ -(K_X+\Delta)$, there exists some $\bQ$-complement $\Gamma$ of $(X,\Delta)$ such that $\Gamma\ge \alpha D$ and $v$ is an lc place of $(X,\Delta+\Gamma)$.
\end{lem}

\begin{proof}
Up to rescaling, we may assume that $A_{X,\Delta}(v)=1$. 
By \cite{BJ-delta}*{Proposition 4.8(ii)} and  \cite{Xu-quasimonomial}*{Theorem 1.1}, the valuation $v$ is quasi-monomial. Let $r$ be the rational rank of $v$. Let $\pi\colon Y\to X$ be a log resolution such that $v$ lies in the interior of $\QM_\eta(Y,E)$ for some simple normal crossing divisor $E=E_1+\cdots+E_r$ on $Y$ where $\eta=c_Y(v)$ is the generic point of a connected component of $\cap_{i=1}^r E_i$.  By \cite{LX-higher-rank}*{Lemma 2.7}, for any $\varepsilon>0$ there exists divisorial valuations $v_1,\cdots,v_r\in \QM_\eta(Y,E)$ and positive integers $q_1,\cdots,q_r$ such that 
\begin{itemize}
    \item $v$ is in the convex cone generated by $v_i$,
    \item for all $i=1,\cdots,r$, the valuation $q_i v_i$ is $\bZ$-valued and has the form $\ord_{F_i}$ for some divisor $F_i$ over $X$, and
    \item $|v_i-v|<\frac{\varepsilon}{q_i}$ for all $i=1,\cdots,r$. Here $|v_i-v|$ denotes the Euclidean distance of $v_i$ and $v$ in $\QM_\eta(Y,E)\cong \bR_{\geq 0}^r$.
\end{itemize}
We claim that when $\varepsilon$ is sufficiently small, there exists a $\bQ$-complement $\Gamma\ge \alpha D$ of $(X,\Delta)$ that has all $v_i$ as lc places. Taking this for granted, let us finish the proof of the lemma. Recall that $A_{X,\Delta}$ is linear on $\QM_\eta(Y,E)$, while $v\mapsto v(\Gamma)$ is concave (see the remarks in Definitions \ref{d-quasimonomial} and \ref{d-log discrepancy}), hence $v\mapsto A_{X,\Delta+\Gamma}(v)=A_{X,\Delta}(v)-v(\Gamma)$ is convex. Since the $v_i$'s are lc places, we have $A_{X,\Delta+\Gamma}(v_i)=0$ and hence as $v$ is contained in their convex hull we get $A_{X,\Delta+\Gamma}(v)\le 0$. As $(X,\Delta+\Gamma)$ is lc, this implies that $v$ is an lc place of $(X,\Delta+\Gamma)$ and therefore the statement of the lemma follows.

Returning to the proof of the claim, we first argue as in \cite{LX-higher-rank}*{Lemma 2.51}\footnote{There is a small error in the proof of \cite{LX-higher-rank}*{Lemma 2.51}, where the log discrepancy function $A_{X,\Delta+\fa_\bullet^c}(\cdot)$ was treated as a linear function on $\QM_\eta(Y,E)$; nonetheless it is Lipschitz (since it is convex) and that is enough for the argument there.}. Let $\fa_\bullet$ be the graded sequence of valuation ideals of $v$, i.e. $\fa_m=\fa_m(v)$. By the same argument as above, we know that the log discrepancy function $w\mapsto A_{X,\Delta+\fa_\bullet}(w)$ is convex on $\QM_\eta(Y,E)$. In particular, it is Lipschitz in a neighborhood of $v$, hence there exist some constants $C>0$ and $\varepsilon_1>0$ such that 
\[
|A_{X,\Delta+\fa_\bullet}(w)-A_{X,\Delta+\fa_\bullet}(v)|\le C|w-v|
\]
for any $w$ in the closed convex cone generated by $\{v_i\}_{i=1}^r$ and any $0<\varepsilon \leq \varepsilon_1$. Applying this to the divisorial valuations $v_i$ above, we find 
\[
A_{X,\Delta+\fa_\bullet}(F_i)=q_i A_{X,\Delta+\fa_\bullet}(v_i)\le C q_i |v_i-v|\le  C\varepsilon.
\]
Therefore, for some $0<\varepsilon_0\ll 1$ (depending on $\varepsilon$) we have 
\begin{equation} \label{eq:A<epsilon}
    A_{X,\Delta+\fa_\bullet^{1-\varepsilon_0}}(F_i)<2C\varepsilon
\end{equation}
for all $1\le i\le r$. Let $0\le D'\sim_\bQ -(K_X+\Delta)$ be general (so that it does not contain the center of $v$ in its support). Now for any $m\in\bN$ such that $-m(K_X+\Delta)$ is very ample, let $G=\beta D'+(1-\beta)D$ where $\beta=\max\{0,\frac{(n+1)(\delta-1)}{\delta}\}$, and let $D_m$ be an $m$-basis type $\bQ$-divisor that is compatible with both $G$ and $v$. Then we have $D_m\ge S_m(G)\cdot G$ and $v(D_m)=S_m(v)$. 

Denote by
\[
D'_m:=D_m-S_m(G)\cdot \beta D'\sim_\bQ -(1-\beta S_m(G))(K_X+\Delta).
\]
Note that $G\sim_\bQ -(K_X+\Delta)$, thus $\lim_mS_m(G)=S_{X,\Delta}(G)=\frac{1}{n+1}$ (see Lemma \ref{lem:S=1/n+1}) and $\lim_{m\to \infty} (1-\beta S_m(G))=\min\{1,\frac{1}{\delta}\}$ by a direct calculation. It follows that we can choose a sequence of rational numbers $\delta_m>0$ $(m\in\bN)$ such that $\delta_m<\delta_m(X,\Delta)$, $\lim_{m\to \infty} \delta_m = \delta$ and $\delta_m(1-\beta S_m(G))<1$ for all $m$. In particular, $(X,\Delta+\delta_m D'_m)$ is log Fano and by our assumption on $\alpha$, we have 
\[
\delta_m D'_m\ge (1-\beta)\delta_m S_m(G)\cdot D\ge \alpha D
\]
as $m\gg 0$. 

Since $v$ computes $\delta(X,\Delta)$ and $D'$ is general, we also see that as $m\gg 0$, 
\[
\delta_m v(D'_m)=\delta_m v(D_m)\ge (1-\varepsilon_0)\delta(X,\Delta)S_{X,\Delta}(v)=(1-\varepsilon_0)A_{X,\Delta}(v)=1-\varepsilon_0.
\]
Combined with \eqref{eq:A<epsilon} we obtain $A_{X,\Delta+\delta_m D'_m}(F_i)\le A_{X,\Delta+\fa_\bullet^{1-\varepsilon_0}}(F_i)<2C\varepsilon<1$ as long as $\varepsilon<\min\{\frac{1}{2C},\varepsilon_1\}$. By \cite{BCHM}*{Corollary 1.4.3} we know that there exists a $\bQ$-factorial birational model $p\colon \tX\to X$ that extracts exactly the divisors $F_i$. Let $\tD$ denote the strict transform of a divisor $D$ on $X$. Let
\[
K_{\tX}+\widetilde{\Delta}+\delta_m \tD'_m +\sum_{i=1}^r (1-a_i) F_i=p^*(K_X+\Delta+\delta_m D'_m)
\]
be the crepant pullback, then $a_i\in (0,2C\varepsilon)$ and as $(X,\Delta+\delta_m D'_m)$ is log Fano, we see that $(\tX,\widetilde{\Delta}+\alpha \tD+\sum_{i=1}^r (1-a_i) F_i)$ has a $\bQ$-complement. 

By the following Lemma \ref{lem:ACC}, if $\varepsilon$ is sufficiently small (depending only on $C$, $\varepsilon_1$, and the coefficients of $\Delta$ and $\alpha D$), then $(\tX,\widetilde{\Delta}+\alpha \tD+\sum_{i=1}^r F_i)$ also have a $\bQ$-complement. Pushing it forward to $X$, we obtain a $\bQ$-complement $\Gamma\ge \alpha D$ of $(X,\Delta)$ that realizes all $F_i$ as lc places. The proof is now complete.
\end{proof}

We have used the following well-known consequence of \cite{HMX-ACC} in the proof above.

\begin{lem} \label{lem:ACC}
Let $(X,\Delta)$ be a projective pair and let $G$ be an effective $\bQ$-Cartier $\bQ$-divisor on $X$. Assume that $X$ is of Fano type. Then there exists some $\varepsilon>0$ depending only on the dimension of $X$ and the coefficients of $\Delta$ and $G$ such that: if $(X,\Delta+(1-\varepsilon)G)$ has a $\bQ$-complement, then the same is true for $(X,\Delta+G)$. 
\end{lem}

\begin{proof}
This should be well-known to experts (see e.g. \cite{Birkar-BAB}*{Proof of Proposition 3.4}), but we provide a proof for readers' convenience. 

Replacing $X$ by a small $\bQ$-factorial modification, we may assume that $X$ itself is $\bQ$-factorial. Let $n=\dim X$ and let $I\subseteq \bQ_+$ be the coefficient set of $\Delta$ and $G$. By the ACC of log canonical thresholds and global ACC of log Calabi-Yau pairs \cite{HMX-ACC}*{Theorems 1.1 and 1.5}, we know that there exists some rational constant $\varepsilon>0$ depending only on $n,I$ which satisfies the following property: for any projective pair $(X,\Delta)$ of dimension at most $n$ and any $\bQ$-Cartier divisor $G$ on $X$ such that the coefficients of $\Delta$ and $G$ belong to $I$, we have $(X,\Delta+G)$ is lc as long as $(X,\Delta+(1-\varepsilon)G)$ is lc; if in addition there exists some divisor $D$ with $(1-\varepsilon)G\le D\le G$ such that $K_X+\Delta+D\sim_\bQ 0$, then $D=G$. 

Now let $(X,\Delta+(1-\varepsilon)G)$ be a pair with a $\bQ$-complement $\Gamma$. As $X$ is of Fano type, we may run the $-(K_X+\Delta+G)$-MMP $f\colon X\dashrightarrow X'$. Let $\Delta',G',\Gamma'$ be the strict transforms of $\Delta,G,\Gamma$. Note that by construction
\[
K_X+\Delta+G \le f^*(K_{X'}+\Delta'+G'),
\]
hence $(X,\Delta+G)$ has a $\bQ$-complement if and only if $(X',\Delta'+G')$ has one. Since 
\[
K_X+\Delta+(1-\varepsilon)G+\Gamma\sim_\bQ 0,
\]
the MMP is crepant for the lc pair $(X,\Delta+(1-\varepsilon)G+\Gamma)$, hence $(X',\Delta'+(1-\varepsilon)G'+\Gamma')$ is also lc. It follows that $(X',\Delta'+(1-\varepsilon)G')$ is lc, thus by our choice of $\varepsilon$, $(X',\Delta'+G')$ is lc as well. Suppose that $X'$ is Mori fiber space $g\colon X'\to S$. Then $K_{X'}+\Delta'+G'$ is $g$-ample. Since $K_{X'}+\Delta'+(1-\varepsilon)G'\sim_\bQ -\Gamma'\le 0$ and $\rho(X')=\rho(S)+1$, there exists some $\varepsilon'\in (0,\varepsilon]$ such that $K_{X'}+\Delta'+(1-\varepsilon')G'\sim_{g,\bQ} 0$. But if we restrict the pair to the general fiber of $X'\to S$ we would get a contradiction to our choice of $\varepsilon$. Thus $X'$ is a minimal model and $-(K_{X'}+\Delta'+G')$ is nef. As $X'$ is also of Fano type, we see that $-(K_{X'}+\Delta'+G')$ is semiample, hence  $(X',\Delta'+G')$ has a $\bQ$-complement. By the previous discussion, this implies that $(X,\Delta+G)$ has a $\bQ$-complement as well.
\end{proof}

To proceed, we make the following definition. Recall that a log smooth model $(Y,E)$ over $(X,\Delta)$ consists of a log resolution $\pi: Y\to (X,\Delta)$ and a reduced divisor $E$ on $Y$ such that $E+\Ex(\pi)+\pi_*^{-1}\Delta$ has SNC support (see Definition \ref{d-logsmooth}).

\begin{defn}\label{d-specialcomplement}
 A $\bQ$-complement $\Gamma$ of $(X,\Delta)$ will be called \emph{special with respect to a log smooth model $\pi\colon(Y,E)\to (X,\Delta)$} 
 if $\Gamma_Y=\pi_*^{-1}\Gamma\ge G$ for some effective ample $\bQ$-divisor $G$ on $Y$ whose support does not contain any stratum of $(Y,E)$. Any valuation $v\in \QM(Y,E)\cap \DMR(X, \Delta+\Gamma)$ is  called a \emph{monomial lc place} of the special $\bQ$-complement $\Gamma$ with respect to $(Y,E)$.
\end{defn}

The following immediate consequence of Lemma \ref{lem:minimizer's complement} says a valuation computing $\delta(X,\Delta)$ when $\delta(X,\Delta)<\frac{n+1}{n}$ is a monomial lc place of a special complement. Later we will show for a (possibly higher rank) valuation, if it is a monomial lc place of a special complement, its associated graded ring is finitely generated (see Theorem \ref{thm:f.g. for special valuation}).


\begin{cor} \label{cor:minimizer has special complement}
Let $(X,\Delta)$ be a log Fano pair of dimension $n$ such that $\delta(X,\Delta)=\delta<\frac{n+1}{n}$, and let $v$ be a valuation that computes $\delta(X,\Delta)$. Then there exists a log smooth model $\pi\colon (Y,E)\to (X,\Delta)$ and a special $\bQ$-complement $0\le \Gamma\sim_\bQ -(K_X+\Delta)$ with respect to $(Y,E)$, such that $v\in \QM(Y,E)\cap \DMR(X, \Delta+\Gamma)$.
\end{cor}

\begin{proof}
Since $v$ is quasi-monomial, we may find a log smooth model $\pi\colon (Y,E)\to (X,\Delta)$ whose exceptional locus supports a $\pi$-ample divisor $F$ such that $v\in \QM(Y,E)$. Choose some $0<\varepsilon\ll 1$ such that $L:=-\pi^*(K_X+\Delta)+\varepsilon F$ is ample and let $G$ be a general divisor in the $\bQ$-linear system $|L|_\bQ$ whose support does not contain any stratum of $(Y,E)$. Let $D=\pi_*G\sim -(K_X+\Delta)$ and let $\alpha<\min\{\frac{\delta}{n+1},1-\frac{n\delta}{n+1}\}$ be a fixed rational positive number. By Lemma \ref{lem:minimizer's complement}, there exists some complement $\Gamma$ of $(X,\Delta)$ such that $\Gamma\ge \alpha D$ and $v$ is an lc place of $(X,\Delta+\Gamma)$. Replace $G$ by $\alpha G$. By construction, the strict transform of $\Gamma$ is larger or equal to $G$, so $\Gamma$ is a special $\bQ$-complement with respect to $(Y,E)$.
\end{proof}

To help understand the importance of special complements, we prove the following statement. In Section \ref{s-finitegeneration}, we will need a stronger version of it (see Theorem \ref{t-alphabelow}).

\begin{lem} \label{lem:special DMR}
Let $\Gamma$ be a special $\bQ$-complement of a log Fano pair $(X,\Delta)$ with respect to a log smooth model $(Y,E)$. Denote by $\Pi:= \QM(Y,E)\cap \DMR(X,\Delta+\Gamma)$ the set of monomial lc places. Then every divisorial valuation $v\in \Pi$ is special $($see Definition \ref{d-specialdivisor}$)$.
\end{lem}

\begin{proof}
Fix $v\in \Pi(\bQ)$. By Theorem \ref{thm:zhuang-special}, it suffices to find an effective divisor $D\sim_\bQ -(K_X+\Delta)$ such that $\lambda=\lct(X,\Delta;D)\in(0,1)$ and that (up to rescaling) $v$ is the unique lc place of $(X,\Delta+\lambda D)$. To see this, let $W=C_Y(v)$. If $W$ has codimension at least two in $Y$, we let $\rho\colon Z\to Y$ be the weighted blowup corresponding to $v$, so that $v=c\cdot\ord_{F}$ for some $c>0$, where $F$ is the exceptional divisor of $\rho$. Otherwise (i.e. $W$ is an irreducible component of $E$) we let $Z=Y$, $\rho=\mathrm{id}$, and $F=W$.

By assumption, there exists an effective ample $\bQ$-divisor $G$ on $Y$ such that $\Gamma_Y\ge G$ and $C_Y(v)\not\subseteq \Supp(G)$. Since $(Y,E)$ is log smooth, it is clear that $F$ is a weighted projective space bundle over the smooth center $W$ and $-F$ is $\rho$-ample. In particular, there exists some $\varepsilon>0$ such that $\rho^*G-\varepsilon F$ is ample on $Z$. Let $G_1$ be a general divisor in the $\bQ$-linear system $|\rho^*G-\varepsilon F|_\bQ$ and consider the effective divisor $D\sim_\bQ -(K_X+\Delta)$ on $X$ satisfying $\rho^*\pi^*D=\rho^*(\pi^*\Gamma-G)+G_1+\varepsilon F$ (this is possible since the right hand side is $\sim_\bQ -\rho^*\pi^*(K_X+\Delta)$). We claim that this divisor $D$ satisfies the desired conditions.

Let $K_Y+\Delta_Y=\pi^*(K_X+\Delta)$ and $K_Z+\Delta_Z=\rho^*(K_Y+\Delta_Y)$ be the crepant pullbacks. We first show that the above claim is a consequence of the following two properties:
\begin{enumerate}
    \item $(Y,\Delta_Y+\pi^*\Gamma-G)$ is sub-lc and $v$ is an lc place of this sub-pair;
    \item $F$ is the only divisor that computes $\lct(Y,\Delta_Y;\rho_*(G_1+\varepsilon F))$.
\end{enumerate}
This is because, (1) implies that 
\[
A_{X,\Delta}(w)=A_{Y,\Delta_Y}(w)\ge w(\pi^*\Gamma-G)
\]
for all valuations $w\in \Val_X^\circ$, and equality holds when $w$ is proportional to $v$; on the other hand, if we let $\mu=\lct(Y,\Delta_Y;\rho_*(G_1+\varepsilon F))>0$, then (2) implies that 
\[
A_{X,\Delta}(w)=A_{Y,\Delta_Y}(w)\ge \mu\cdot w(\rho_*(G_1+\varepsilon F))
\]
for all valuations $w\in \Val_X^\circ$, and equality holds if and only if $w$ is proportional to $v$. Combining the two inequalities we have
\[
w(D)=w(\pi^*D)=w(\pi^*\Gamma-G+\rho_*(G_1+\varepsilon F))\le (1+\mu^{-1})A_{X,\Delta}(w)
\]
for all valuations $w\in \Val_X^\circ$, and equality holds if and only if $w$ is proportional to $v$. In particular, $\lct(X,\Delta;D)=\frac{1}{1+\mu^{-1}}\in(0,1)$ and up to rescaling, $v$ is the unique lc place that computes this lct, which is exactly what we want.

It remains to prove the two properties above. Point (1) is quite straightforward since by assumption $v$ is an lc place of the sub-lc sub-pair $(Y,\Delta_Y+\pi^*\Gamma)$ and $G$ does not contain $C_Y(v)$. To see point (2), we note that by assumption, $\Delta_Y$ has simple normal crossing support, $\lfloor \Delta_Y \rfloor \le 0$ (as $(X,\Delta)$ is klt), and therefore one can easily check that the sub-pair $(Z,\Delta_Z\vee F)$ is plt. Here we denote by $D_1\vee D_2$ the smallest $\bQ$-divisor $D$ such that $D\ge D_i$ for $i=1,2$. Let $t=\frac{A_{Y,\Delta_Y}(F)}{\varepsilon}$. Then
\[
\rho^*(K_Y+\Delta_Y+t\rho_*(G_1+\varepsilon F))=K_Z+\Delta_Z\vee F+t G_1
\]
by construction. Since $G_1$ is general, the pair $(Z,\Delta_Z\vee F+t G_1)$ is also plt. This proves (2). The proof is now complete.
\end{proof}

\subsection{Minimizers and constructibility}\label{ss-n+1/n}
The existence of a valuation computing $\delta(X,\Delta)$ is proved to exist in \cite{BJ-delta}*{} if the ground field $\bk$ is uncountable, and in \cite{BLX-openness} when $\delta(X,\Delta)\le 1$ for a general ground field, where it is also shown that in this case any minimizer is an lc place of a $\mathbb{Q}$-complement. Here we extend these results to the case when $\delta(X,\Delta)$ is bounded by $\frac{n+1}{n}$.

\begin{thm}\label{t-n+1/n}
Let $(X,\Delta)$ be a log Fano pair of dimension $n$ such that $\delta(X,\Delta)<\frac{n+1}{n}$. Then 
\begin{enumerate}
    \item there exists a valuation computing $\delta(X,\Delta)$; and
    \item there exists a positive integer $N$ depending only on $\dim(X)$ and the coefficients of $\Delta$ such that for any valuation $v$ computing $\delta(X,\Delta)$, there exists an $N$-complement $D$ of $(X,\Delta)$ which satisfies that $v$ is an lc place of $(X,\Delta+D)$.
\end{enumerate}
\end{thm}
\begin{proof}
First we prove (1). For any sufficiently divisible $m\in\bN$, let $\delta_m:=\delta_m (X,\Delta)$, and let $E_m$ be a divisor over $X$ such that $\frac{A_{X,\Delta}(E_m)}{S_m(E_m)}=\delta_m$.  

Fix a positive integer $m_0$ such that $-m_0(K_X+\Delta)$ is a very ample Cartier divisor. For each sufficiently divisible $m\in\bN$, let $H_m$ be a smooth divisor in the linear system $|-m_0(K_X+\Delta)|$ that does not contain the center of $E_m$. For any such $m$, we can find an $m$-basis type divisor $D_m$ which is compatible with both $E_m$ and $H_m$ by Lemma \ref{lem:compatible}. We write $D_m=\Gamma_m+a_m H_m$ where ${\rm Supp}(\Gamma_m)$ does not contain $H_m$. Note that by Definition \ref{defn-Sinvariant},
\begin{align*}
a_m &=\frac{1}{m \dim H^0(-m(K_X+\Delta))}\sum_{j=1}^\infty \dim H^0(X,-m(K_X+\Delta)-jH_m)\\
&=\frac{1}{m \dim H^0(-m(K_X+\Delta))}\sum_{j=1}^\infty \dim H^0(X,-(m-jm_0)(K_X+\Delta))
\end{align*}
does not depend on the choice of $H_m$ and in particular $\lim_{m\to \infty} a_m=\frac{1}{m_0(n+1)}$ by Lemma \ref{lem:S=1/n+1}. We know that
$$\lct(X,\Delta;D_m)\le \frac{A_{X,\Delta}(E_m)}{\ord_{E_m}(D_m)} =\frac{A_{X,\Delta}(E_m)}{S_m(E_m)}=\delta_m,$$
where the equality $\ord_{E_m}(D_m)=S_m(E_m)$ follows from the fact that $D_m$ is chosen to be compatible with $E_m$. However, we have $\lct(X,\Delta;D_m)\ge \delta_m$ by the definition of $\delta_m$. Thus $\lct(X,\Delta;D_m) = \delta_m$ and the lct is computed by $E_m$. Since $H_m$ does not contain the center of $E_m$, it follows that $(X,\Delta+\delta_m\Gamma_m)$ is lc and $E_m$ is an lc place of this pair.

Note that $\lim_{m\to\infty} \delta_m=\delta(X,\Delta)<\frac{n+1}{n}$. So for sufficiently large $m$, we get
$$\delta_m\Gamma_m=\delta_m(D_m-a_m H_m)\sim_{\mathbb Q}-\lambda_m(K_X+\Delta)$$
where $\lambda_m=\delta_m(1-m_0 a_m)\in (0,1)$. Thus $E_m$ is an lc place of a $\mathbb{Q}$-complement. The rest of the proof is the same as in \cite{BLX-openness}*{Theorems 4.6 and A.7}: we know that $E_m$ is indeed an lc place of an $N$-complement for some $N$ that only depends on $\dim(X)$ and ${\rm Coeff}(\Delta)$. 
Therefore, after passing to a subsequence, we can find an $N$-complement $D$, together with lc places $F_m$ of $(X,\Delta+D)$, 
such that $\frac{A_{X,\Delta}(E_m)}{S_{X,\Delta}(E_m)}=\frac{A_{X,\Delta}(F_m)}{S_{X,\Delta}(F_m)}$ for all sufficiently divisible $m\in\bN_+$. If we take $v$ to be the limit of $(A_{(X,\Delta)}(F_m))^{-1} \ord_{F_m}$ in $\DMR(X,\Delta+D)$, then $v$ computes $\delta(X,\Delta)$, as 
$$\frac{A_{X,\Delta}(v)}{S_{X,\Delta}(v)}=\lim_{m\to\infty} \frac{A_{X,\Delta}(F_m)}{S_{X,\Delta}(F_m)}=\lim_{m\to\infty} \frac{A_{X,\Delta}(E_m)}{S_{X,\Delta}(E_m)}=\lim_{m\to\infty} \delta_m=\delta(X,\Delta).$$ 

For (2), it follows immediately from Lemma \ref{lem:minimizer's complement} that $v$ is an lc place of a $\mathbb Q$-complement $\Gamma$. There exists a log smooth model $(Y,E)\to (X,\Delta+\Gamma)$ such that every component $E_i$ $(i=1,\cdots,q)$ of $E$ is an lc place of $(X,\Delta+\Gamma)$ and every prime divisor on $Y$ with log discrepancy 0 with respect to $(X,\Delta+\Gamma)$ is contained in $E$. In particular,  $v\in \QM(Y,E)$. By \cite{BCHM}*{Corollary 1.4.3}, there exists a $\bQ$-factorial birational model $\mu\colon \tX\to X$ that extracts exactly the divisors $E_i$ and $Y\dasharrow \tX$ is isomorphic at the generic point of any component of all non-empty intersections of $\bigcap_{i\in J} E_j$ for $J\subset \{1,\cdots,q\}$. Let $a_i={\rm Coeff}_{E_i}(\Delta)$ if $E_i$ is a prime divisor on $X$, otherwise let $a_i=0$. Then we can argue as in the proof of \cite{BLX-openness}*{Lemma 3.2 and Theorem 3.5}: $\tX$ is of Fano type, $(\tX,\mu_*^{-1}\Delta+\sum_{i=1}^q (1-a_i)E_i)$ has a $\bQ$-complement, therefore also has an $N$-complement, whose pushforward on $X$ gives an $N$-complement $D$ of $(X,\Delta)$ that has all $E_i$ ($i=1,\cdots,q$) as lc places. In particular, it also has $v$ as an lc place (as in the proof of Lemma \ref{lem:minimizer's complement}).
 \end{proof}

Combining Theorem \ref{t-n+1/n} with the argument in \cite{BLX-openness}*{Section 4}, we get the following generalization of \cite{BLX-openness}*{Theorem 1.1}.

\begin{cor}
For a $\mathbb{Q}$-Gorenstein family of log Fano pairs $(X,\Delta)\to S$ over a normal base, the function
$$t\in S\mapsto \min \left\{ \frac{n+1}{n}, \delta(X_{\bar{t}},\Delta_{\bar{t}}) \right\}$$
is lower semi-continuous and constructible, where $(X_{\bar{t}},\Delta_{\bar{t}})$ is the base change to the algebraic closure of $k(t)$.
\end{cor}

\begin{rem}
It is proved in \cite{Z-equivariant} that
$$\min\{\delta(X_{\bar{t}},\Delta_{\bar{t}}),1\}=\min\{ \delta(X_t,\Delta_t),1\}.$$ 
However, in general we may have
$\delta(X_{\bar{t}},\Delta_{\bar{t}})<\delta(X_{{t}},\Delta_{{t}})$ (see \cite{CP-Proj}*{Remark 4.16}).
\end{rem}

\section{Finite generation}\label{s-finitegeneration}

In this section, we prove the following finite generation result.

\begin{thm} \label{thm:f.g. for delta minimizer}
Let $(X,\Delta)$ be a log Fano pair of dimension $n$ and let $r>0$ be an integer such that $r(K_X+\Delta)$ is Cartier. Let $R=\bigoplus_{m\in\bZ_{\geq 0}} H^0(X,-mr(K_X+\Delta))$. Assume that $\delta(X,\Delta)<\frac{n+1}{n}$. Then for any valuation $v$ that computes $\delta(X,\Delta)$, the associated graded ring $\gr_v R$ is finitely generated.
\end{thm}

To tackle Theorem \ref{thm:f.g. for delta minimizer}, we need some finite generation criterion for lc places of complements. As shown by the examples in \cite{AZ-index2}*{Theorem 1.4} and Section \ref{s-examples}, one needs extra assumptions on the valuation and the complement. It turns out that the special complements as defined in Definition \ref{d-specialcomplement} will be the correct one for proving Theorem \ref{thm:f.g. for delta minimizer}. In other words, monomial lc places of special complements have finitely generated associated graded rings (see Theorem \ref{thm:f.g. for special valuation}).

To prove this, we will show that in a neighborhood of a monomial lc place of a special complement, the divisorial ones induce degenerations to log Fano pairs whose alpha invariants are bounded from below by a positive constant (see Theorem \ref{t-alphabelow}), and this is sufficient for the finite generation (see Theorem \ref{thm:f.g. criterion}).

\subsection{Finite generation criterion}

The next statement gives us the finite generation criterion. Clearly, it implies Theorem \ref{thm:f.g. for delta minimizer} by Corollary \ref{cor:minimizer has special complement}.

\begin{thm} \label{thm:f.g. for special valuation}
Let $(X,\Delta)$ be a log Fano pair, and let $v$ be a quasi-monomial valuation on $X$. Let $R=\bigoplus_{m\in\bN} H^0(X,-mr(K_X+\Delta))$ for some integer $r>0$ such that $r(K_X+\Delta)$ is Cartier. Then the following are equivalent.
\begin{enumerate}
    \item The associated graded ring $\gr_v R$ is finitely generated and the central fiber $(X_v,\Delta_v)$ of the induced degeneration is klt.
    \item The valuation $v$ is a monomial lc place of a special $\bQ$-complement $\Gamma$ with respect to some log smooth model $(Y,E)$ $($see Definition \ref{d-specialcomplement}$)$. 
\end{enumerate}
\end{thm}

We explain some notation in the above theorem. Assuming the finite generation of  $\gr_v R$, we define $X_v:=\Proj~\gr_v R$, and $\Delta_v$ is the induced degeneration of $\Delta$ to $X_v$. More precisely, suppose $\Delta=\sum_{i=1}^{l} a_i \Delta_i$ where $\Delta_i$ is a prime divisor on $X$ and $a_i\in \bQ_{\geq 0}$. Let $I_{\Delta_i}\subseteq R$ be the graded ideal of $\Delta_i$. Let $\mathrm{in}(I_{\Delta_i})\subseteq \gr_v R$ be the initial ideal of $I_{\Delta_i}$. Then $\Delta_v:=\sum_{i=1}^l a_i \Delta_{v,i}$, where $\Delta_{v,i}$ is the divisorial part of the closed subscheme $V(\mathrm{in}(I_{\Delta_i}))\subseteq X_v$, i.e. $\Delta_{v,i}$ and $V(\mathrm{in}(I_{\Delta_i}))$ coincide away from a codimension $2$ subset of $X_v$. See  \cite[p. 648]{BX-uniqueness} for a similar definition.

The remaining part of this section will be devoted to the proof of this theorem. In this subsection, we reduce the proof to showing the boundedness of the degenerations induced by divisorial valuations that are sufficiently closed to $v$.


We first prove the easier direction in Theorem \ref{thm:f.g. for special valuation}. 

\begin{lem} \label{lem:complement for f.g. special valuation}
Assume that $\gr_v R$ is finitely generated and $(X_v,\Delta_v)$ is klt. Then $v$ is a monomial lc place of a special $\bQ$-complement.
\end{lem}

\begin{proof}
Let $q$ be the rational rank of $v$. Let $\pi\colon (Y,E)\to (X,\Delta)$ be a log smooth model such that $v\in \QM_\eta(Y,E)$ for some codimension $q$ point $\eta\in Y$. Let $D\sim_\bQ -(K_X+\Delta)$ be the divisor constructed in the proof of Corollary \ref{cor:minimizer has special complement} so that the strict transform $G=\pi_*^{-1}D$ is ample. We have $D=\frac{1}{m_0 r}\{f=0\}$ for some $m_0\in\bN$ and some $f\in H^0(X,-m_0 r(K_X+\Delta))$. By assumption, there exists some $f_0:=f,f_1,\cdots,f_\ell\in R$ whose restrictions form a (finite) set of generators $\bar{f}_0,\cdots, \bar{f}_\ell$ of $\gr_v R$ (in particular, $f_0,\cdots, f_\ell$ generates $R$). By enlarging the set of generators, we may also assume that all $I_{\Delta_i}\subseteq R$ are generated by the restrictions of some elements from $f_0,...,f_{\ell}$.

 By assumption, $(X_v,\Delta_v+\varepsilon D_v)$ is klt for some rational constant $0<\varepsilon\ll 1$, thus by Lemma \ref{l-finitegeneratedimpliesisom}, $(X_w,\Delta_w+\varepsilon D_w)\cong (X_v,\Delta_v+\varepsilon D_v)$ is also klt for divisorial valuations $w$ in a sufficiently small neighbourhood $U\subseteq \Sigma:=\QM_\eta(Y,E)$ of $v$. In particular, since $v$ lies in the interior of $\Sigma$ (by construction), we may assume that the closure $\bar{U}$ is a compact subset of ${\rm int}(\Sigma)$. By Theorem \ref{t-complements}, there exists an integer $N$ that only depends on $(X,\Delta+\varepsilon D)$ such that any divisorial valuation $w_0\in U$ is an lc place of an $N$-complement $0\le \Gamma_0\sim_\bQ -(K_X+\Delta+\varepsilon D)$. Recall that $v(f)$ is computed as the smallest weight of monomials in the power series expansion of $f$ at the point $\eta$. As $\Gamma$ varies among the $N$-complements and $w$ varies in a small neighbourhood of $v$, we have $a\cdot \mult_\eta \pi^*\Gamma\le w(\Gamma)<C$ for some constant $a,C>0$ that only depends on $v$. Since there are only finitely many monomials with bounded multiplicity, we conclude that the value of $w(\Gamma)$ is determined by only finitely many such monomials. Hence by shrinking the neighbourhood $U$, we may assume that whenever $\Gamma$ is an $N$-complement of $(X,\Delta+\varepsilon D)$ and $v(\Gamma)\neq A_{X,\Delta+\varepsilon D}(v)$ then $w(\Gamma)\neq A_{X,\Delta+\varepsilon D}(w)$ for any $w\in U$. In particular since $w_0(\Gamma_0)= A_{X,\Delta+\varepsilon D}(w_0)$ for the $N$-complement $\Gamma_0$ we constructed above, we have $v(\Gamma_0)= A_{X,\Delta+\varepsilon D}(v)$ and therefore $v$ is also an lc place of $(X,\Delta+\Gamma')$, where $\Gamma'=\varepsilon D+\Gamma_0$. Since $\pi_*^{-1}\Gamma'\ge \varepsilon G$ and $G$ is ample, it is a special $\bQ$-complement with respect to $(Y,E)$ by construction. In other words, $v$ is a monomial lc place of a special $\bQ$-complement as desired.\end{proof}

The following lemma follows from the proof of \cite{LX-higher-rank}*{Lemma 2.10}, which we reproduce here for reader's convenience.

\begin{lem}\label{l-finitegeneratedimpliesisom}
Let $q$ be the rational rank of $v$. Let $\pi\colon (Y,E)\to (X,\Delta)$ be a log smooth model such that $v\in \QM_\eta(Y,E)$ for some codimension $q$ point $\eta\in Y$. Assume that $\gr_v R$ is finitely generated by the restrictions of homogeneous elements $f_0,\cdots,f_\ell \in R$. Then for all valuations $w\in \Sigma:=\QM_\eta(Y,E)$ that are sufficiently close to $v$, we have an isomorphism $\gr_w R\cong \gr_v R$ sending the restrictions of $f_0,\cdots,f_\ell$ in $\gr_v R$ to their respective restrictions in $\gr_w R$.
\end{lem} 

\begin{proof}
Since $\bar{f}_0,\cdots, \bar{f}_\ell$ generate $\gr_v R$, we have a surjection $\pi_v\colon \bk[x_0,\cdots,x_\ell]\to \gr_v R$ sending $x_i$ to $\bar{f}_i$. Similarly we have a homomorphism $\pi_w\colon \bk[x_0,\cdots,x_\ell]\to \gr_w R$ sending $x_i$ to the restriction of $f_i$. Let us first show that $\pi_w$ factors through $\gr_v R$ when $w$ is sufficiently close to $v$.

The map $\pi_v$ is easily seen to be a doubly graded homomorphism if we set $\deg(x_i)=(m_i,v(f_i))$ where $f_i\in R_{m_i}$. Let $\Phi_1,\cdots,\Phi_p$ be homogeneous generators of its kernel. Let $(y_1,\cdots,y_q)$ be a regular system of parameters of $\cO_{Y,\eta}$ and let $\alpha\in\bR^q_+$ be such that $v=v_\alpha$. By construction, we have $v(\Phi_i(f_0,\cdots,f_\ell))>\wt_{\alpha}(\Phi_i)$ where we set $\wt_{\alpha}(x_i)=v_\alpha (f_i)$ and it induces a natural weight on every polynomial in $\bk[x_0,...,x_\ell]$. Each $f_i$ has a local expansion $f_i=\sum_{\beta\in \bZ_{\geq 0}^q} c^{(i)}_\beta y^\beta$ at $\eta$, where we use the same notation from Definition \ref{d-quasimonomial}. Since $v$ has rational rank $q$, for any homogeneous element $f\in R$ we have $v(f)=\langle \alpha,\beta_f\rangle$ for some uniquely determined $\beta_f\in \bZ_{\geq 0}^q$. In particular, we have $v(f_i)=\langle \alpha,\beta_i\rangle$ for some $\beta_i \in \bZ_{\geq 0}^q$; moreover, for any other $\beta\in \bZ_{\geq 0}^q$ with $c^{(i)}_\beta \neq 0$, we have $\langle \alpha,\beta\rangle > v(f_i)$. It follows that if $\alpha'\in \Sigma$ is sufficiently close to $\alpha$ then $w=v_{\alpha'}$ satisfies $w(f_i)=\langle \alpha',\beta_i\rangle$. Using a similar argument, we also see that if $\alpha'$ is sufficiently close to $\alpha$ then $w(\Phi_i(f_0,\cdots,f_\ell))>\wt_{\alpha'}(\Phi_i)$. This implies that the $\Phi_i$'s are contained in the kernel of $\pi_w$; in particular, the map $\pi_w$ factors through $\gr_v R$. 

Denote by $\varphi\colon \gr_v R\to \gr_w R$ the induced map. We proceed to show that $\varphi$ is an isomorphism. First we show that it is injective. For this it suffices to show that for any $\Phi\in \bk[x_0,\cdots,x_\ell]$ such that the monomials in $\Phi$ have the same weight with respect to $\wt_{\alpha'}$ and that $\pi_v(\Phi)\neq 0$, we have $\pi_w(\Phi)\neq 0$. After disregarding monomials in $\Phi$ with zero image under $\pi_v$, we may further write $\Phi=\Phi'+\Phi''$ where $\pi_v(\Phi')\neq 0$, all monomials in $\Phi'$ have $\wt_\alpha = \wt_\alpha (\Phi)$ while the ones in $\Phi''$ have $\wt_\alpha > \wt_\alpha (\Phi)$. Let $g=\Phi(f_0,\cdots,f_\ell)$ and let $u_i=c^{(i)}_{\beta_i}y^{\beta_i}$ be the monomial in the local expansion of $f_i$ that computes $v(f_i)$. We aim to show $w(g)=\wt_{\alpha'}(\Phi)$ which is equivalent to saying $\pi_w(\Phi)\neq 0$.

As $\pi_v(\Phi')\neq 0$ (i.e. $v(g)=\wt_\alpha(\Phi')$), we must have $\Phi'(u_0,\cdots,u_\ell)\neq 0$: this is the only monomial in the local expansion (at $\eta\in Y$) of $g$ that can have weight $\wt_{\alpha}(\Phi)$ (here we used the fact that $v=v_\alpha$ has rational rank $q$). 
Since the monomial $\Phi'(u_0,\cdots,u_\ell)$ appears in the expansion of $\Phi(u_0,\cdots,u_\ell)$, therefore, $w(\Phi(u_0,\cdots,u_\ell))\le w(\Phi'(u_0,\cdots,u_\ell))=\wt_{\alpha'}(\Phi')=\wt_{\alpha'}(\Phi)$ (the last equality holds since all monomials in $\Phi$ have the same $\wt_{\alpha'}$).

Since $w(f_i-c^{(i)}_{\beta_i}y^{\beta_i})>w(f_i)=w(u_i)$ by our choice of $w$, then
\[
w(\Phi(f_0,\cdots,f_\ell)-\Phi(u_0,\cdots,u_\ell))>\wt_{\alpha'}(\Phi)
\]
(again we use the assumption that all monomials in $\Phi$ have the same $\wt_{\alpha'}$).
So we have
\[
w(g)=w(\Phi(u_0,\cdots,u_\ell))\le \wt_{\alpha'}(\Phi).
\]
On the other hand, we necessarily have $w(g)\ge \wt_{\alpha'}(\Phi)$. So $w(g)= \wt_{\alpha'}(\Phi)$ and therefore $\pi_w(\Phi)\neq 0$. This proves that $\varphi\colon \gr_v R\to \gr_w R$ is injective. But $\varphi$ is a graded homomorphism, and both $\gr_v R$ and $\gr_w R$ have the same dimensions ($=\dim R_m$) in degree $m$, so $\varphi$ is also surjective. Clearly $\varphi$ sends the restrictions of $f_0,\cdots,f_\ell$ in $\gr_v R$ to their respective restrictions in $\gr_w R$. This finishes the proof of the lemma.
\end{proof}

The reverse direction of Theorem \ref{thm:f.g. for special valuation} is much harder. To this end, we drop the assumptions on the complements and prove a weaker finite generation criterion. 

\begin{thm} \label{thm:f.g. criterion}
Let $(X,\Delta)$ be a log Fano pair and let $0\le \Gamma\sim_\bQ -(K_X+\Delta)$ be a $\bQ$-complement. Let $v_0$ be an lc place of $(X,\Delta+\Gamma)$ and let $\Sigma\subseteq\DMR(X,\Delta+\Gamma)$ be the minimal rational PL subspace containing $v_0$ induced by a fixed log smooth model of $(X,\Delta)$ $($see Definition \ref{d-DMR}$)$. Then the following are equivalent:
\begin{enumerate}
    \item The associated graded ring $\gr_{v_0} R$ is finitely generated.
    \item There exists an open neighbourhood $v_0\in U\subseteq \Sigma$ such that the set
    \[
    \{(X_v,\Delta_v)\,|\,v\in U(\bQ):=U\cap \Sigma(\bQ)\}
    \]
    is bounded.
    \item The $S$-invariant function
    \[
    v\mapsto S_{X,\Delta}(v)
    \]
    is linear on a neighbourhood of $v_0$ in $\Sigma$.
\end{enumerate}
\end{thm}

We first handle the implication $(3)\Rightarrow (1)$. This is done in the next two lemmas by studying the concavity of the $S$-invariant function.

\begin{lem} \label{lem:S linear imply f.g. interval}
Let $0\le \Gamma\sim_\bQ -(K_X+\Delta)$ be a $\bQ$-complement of the log Fano pair $(X,\Delta)$. Let $v_0,v_1\in \DMR(X,\Delta+\Gamma)$ be divisorial valuations in the same simplex determined by a fixed log smooth model $(Y,E)$. There exists a natural linear map $[0,1]\to \QM(Y,E)$ sending $0\mapsto v_0$ and $1\mapsto v_1$. For $t\in (0,1)$ we then denote $v_t$ to be the valuation corresponding to $t$. Then
\begin{enumerate}
    \item $S_{X,\Delta}(v_t)\ge (1-t)S_{X,\Delta}(v_0)+t S_{X,\Delta}(v_1)$.
    \item When equality holds, we have
    \[
    \cF_{v_t}^\lambda R_m = {\rm Span}\{s\in R_m\,|\,(1-t)v_0(s)+t v_1(s)\ge \lambda\}
    \]
    for all $\lambda\in\bR$ and all $m\in\bN$. In particular, the filtration $\cF_{v_t}$ is finitely generated.
\end{enumerate}
\end{lem}

\begin{proof}
Let $\cF_0=\cF_{v_0}$, $\cF_1=\cF_{v_1}$ and let $\cF_t$ be the filtration given by 
\[
\cF_t^\lambda R_m = {\rm Span}\{s\in R_m\,|\,(1-t)v_0(s)+t v_1(s)\ge \lambda\}.
\]
We claim that $S(\cF_t)=(1-t)S(v_0)+t S(v_1)$ and that $\gr_{v_0}(\gr_{v_1} R)\cong \gr_{\cF_t}R$ for any $t\in (0,1)$, where by abuse of notation we denote by $\gr_{v_0}(\gr_{v_1} R)$ the graded ring associated to the filtration on $\gr_{v_1} R$ that is induced by $\cF_{v_0}$ on $R$ (in fact, both claims hold even without assuming that $v_0,v_1$ are divisorial). To see this, we note that 
\[
\gr^\alpha_{v_0}\gr^\beta_{v_1} R \cong \frac{\cF_0^\alpha R \cap \cF_1^\beta R}{(\cF_0^{>\alpha} R \cap \cF_1^\beta R) + (\cF_0^\alpha R \cap \cF_1^{>\beta} R)}
\]
naturally maps to $\gr_{\cF_t}^{(1-t)\alpha+t\beta}R$. This induces the map $\varphi\colon \gr_{v_0}(\gr_{v_1} R) \to \gr_{\cF_t}R$. Let us check that this is an isomorphism. Since both sides are $\bN$-graded and have the same dimensions ($=\dim R_m$) in degree $m$, it suffices to check that $\varphi$ is surjective; but this is clear from the definition of $\cF_t$. Hence $\varphi$ is an isomorphism. We next pick a basis $s_1,\cdots,s_{N_m}$ of $R_m$ that is compatible with both $\cF_0$ and $\cF_1$ by Lemma \ref{lem:compatible}. It is straightforward to check (using $\gr_{v_0}\gr_{v_1} R\cong \gr_{\cF_t}R$) that this basis is also compatible with $\cF_t$, thus 
\[
S_m(\cF_t)=\frac{1}{mrN_m}\sum_{i=1}^{N_m} ((1-t)v_0(s_i)+tv_1(s_i)) = (1-t)S_m(v_0)+t S_m(v_1).
\]
The claim then follows by letting $m\to \infty$. Clearly $\cF_t^\lambda R_m\subseteq \cF_{v_t}^\lambda R_m$, hence $S(v_t)\ge S(\cF_t)=(1-t)S(v_0)+t S(v_1)$, which proves (1).

Next assume that equality holds for some $t\in (0,1)$. Then by concavity it holds for all $t\in (0,1)$. By the above proof we also have $\vol(\cF_{v_t}^{(\lambda)} R)=\vol(\cF_{t}^{(\lambda)} R)$ for all $\lambda\in \bR$. Assume for the moment that $t\in \bQ$. We claim that 
\begin{claim*}
The filtration $\cF_t$ is finitely generated and induces a weakly special degeneration of $(X,\Delta)$.
\end{claim*}
To see this, let $o\in (X',\Delta'+\Gamma')$ be the affine cone over $(X,\Delta+\Gamma)$ using the polarization $-r(K_X+\Delta)$. Note that $(X',\Delta')$ is klt and $(X',\Delta'+\Gamma')$ is lc by \cite{Kol13}*{Lemma 3.1}. Let $w_i$ $(i=0,1)$ be the divisorial valuation over $(X',\Delta'+\Gamma')$ given by $w_i(s)=m+v_i(s)$ where $s\in R_m$, and let $E_i$ $(i=0,1)$ be the corresponding prime divisor over $X'$. Then $w_i$ are both lc place of $(X',\Delta'+\Gamma')$ that are centered at the vertex $o$. Let 
\[
\fa_m=\fa_{mA_{X',\Delta'}(w_0)}(w_0)\cap \fa_{mA_{X',\Delta'}(w_1)}(w_1)
\]
for some sufficiently divisible integer $m>0$. Then $\fa_m$ is $\fm_o$-primary. As $\cO_{X'}(-m\Gamma')\subseteq \fa_m$ and $(X',\Delta'+\Gamma')$ is log canonical,  it follows that $(X',\Delta'+\fa^{1/m}_m)$ is lc and has both $w_0,w_1$ as lc places. By \cite{Xu-HRFG}*{Theorem 3.5 and Proposition 3.6}, we then have a $\bG_m^2$-equivariant locally stable family $(\fX,\Delta_{\fX})$ over $\bA^2=\bA^2_{t_0,t_1}$ with general fiber $(X',\Delta')$  given by
\[
\fX=\Spec \bigoplus_{(i_0,i_1)\in\bZ^2} (\fa_{i_0}(\ord_{E_0})\cap \fa_{i_1}(\ord_{E_1}))t_0^{-i_0}t_1^{-i_1}.
\]
In particular, the above algebra (denoted by $\cA$) is finitely generated. On the other hand, we have $\cA/(t_0,t_1)\cong \gr_{v_0}(\gr_{v_1}R) \cong \gr_{\cF_t} R$, hence $\gr_{\cF_t} R$ is finitely generated. For suitable $a,b\in\bN$, the base change $(\fX,\Delta_{\fX})\times_{\bA^2} \bA^1$ via the map $\bA^1\to\bA^2: s\mapsto (s^a,s^b)$ is the test configuration induced by $\cF_t$. Indeed, by the definition of $\cF_t$ it is not hard to see that one can choose $a,b$ such that
\[
\fX\times_{\bA^2} \bA^1 \cong \Spec \bigoplus_{j\in \bZ} s^{-j}\cF_t^{\lambda j} R
\]
for some $\lambda\in \bQ_+$, and this determines the test configuration. Since $(\fX,\Delta_{\fX})\times_{\bA^2} \bA^1$ is a locally stable family with normal general fiber, we see that $\fX\times_{\bA^2} \bA^1$ is normal, thus the algebra on the above right hand side is integrally closed. In particular, for any $p\in \bN_+$, we have
\begin{equation} \label{eq:F_t saturated}
    f\in \cF_t^\lambda R_m \Longleftrightarrow f^p\in \cF_t^{\lambda p} R_{mp}.
\end{equation}
As $(\fX,\Delta_{\fX})$ is a locally stable family, its central fiber 
\[
(\Spec~ \cA/(t_0,t_1), \Delta_{\fX,(0,0)}) = (\Spec (\gr_{\cF_t} R), \Delta_{\fX,(0,0)})
\]
is slc. Since this is an orbifold cone over the central fiber of the test configuration induced by $\cF_t$, we deduce from \cite{Kol13}*{Lemma 3.1} (see also \cite{Kol04}) that the test configuration induced by $\cF_t$ is weakly special. This proves the claim.


Next, let $q>0$ be the smallest integer such that the valuation $qv_t$ is $\bZ$-valued; it is also the smallest integer such that the set of jumping numbers $\cF_t$ lies in $\frac{1}{q}\bZ$. By the above claim and (the proof of) \cite{BLX-openness}*{Theorem A.2}, there exists some $\bZ$-valued divisorial valuations $w_i$ ($i=1,\cdots,\ell$) and some $a_i\in \bQ$ such that 
\[
\cF_t^\lambda R_m=\{f\in R_m\,|\,w_i(f)\ge \lambda q+m a_i \mbox{ for all } 1\le i\le \ell\}.
\]
Indeed, by \cite{BHJ-DH-measure}*{Proposition 2.15}, the above equality holds at least for sufficiently divisible $m$; but then it holds for every $m\in\bN$ by \eqref{eq:F_t saturated}. Suppose that $\cF_t\neq \cF_{v_t}$. Then there exists some $f\in R_m$ such that $v_t(f)=\mu>0$ and $f\in \cF_t^\lambda R_m\setminus \cF_t^{>\lambda} R_m$ for some $\lambda<\mu$. In particular, $\lambda\in \frac{1}{q}\bZ$, and for at least one of the valuations $w_i$, say $w_1$, we have $w_1(f)=\lambda q+ma_1$. Let $\varepsilon\in \bQ_+$ be sufficiently small so that $\eta:=(\mu-\lambda)q+\varepsilon ma_1>0$. Then for sufficiently divisible integer $k\in\bN$, the kernel of the map 
\[
R_{\varepsilon mk} \stackrel{\cdot f^k}{\longrightarrow} \cF_{v_t}^{\mu k} R_{(1+\varepsilon) mk}/\cF_t^{\mu k} R_{(1+\varepsilon)mk}
\]
is contained in $\cF^{\eta k}_{w_1} R_{\varepsilon mk}$. It follows that 
\[
\dim (\cF_{v_t}^{\mu k} R_{(1+\varepsilon) mk}/\cF_t^{\mu k} R_{(1+\varepsilon)mk}) \ge \dim (R_{\varepsilon mk}/\cF^{\eta k}_{w_1} R_{\varepsilon mk})
\]
and thus dividing out by $k^n/n!$ and letting $k\to \infty$ we deduce from \cite{BHJ-DH-measure}*{Lemma 5.13} that 
\[
\vol(\cF_{v_t}^{(\lambda')} R)- \vol(\cF_{t}^{(\lambda')} R)>0\quad\mbox{where}\quad\lambda'=\frac{\mu}{(1+\varepsilon)m},
\]
a contradiction. Hence the two filtrations $\cF_t$ and $\cF_{v_t}$ coincides when $t\in \bQ$.

In general, for any fixed $m\in\bN$, $a,b\in\bQ$ and any irrational $t\in (0,1)$, we have 
\begin{equation} \label{eq:F_t perturb}
    \cF_t^{(1-t)a+tb} R_m = \cF_{t'}^{(1-t')a+t'b} R_m \quad\mbox{and}\quad
    \cF_{v_t}^{(1-t)a+tb} R_m = \cF_{v_{t'}}^{(1-t')a+t'b} R_m
\end{equation}
for some $t'\in \bQ$ that is sufficiently close to $t$. Indeed, recall that $s_1,\cdots,s_{N_m}$ is a basis of $R_m$ that's compatible with both $\cF_0$ and $\cF_1$, and $\cF_t^\lambda R_m$ is spanned by those $s_i$ with $(1-t)v_0(s_i)+tv_1(s_i)\ge \lambda$. As $t\not\in\bQ$, we have $(1-t)v_0(s_i)+tv_1(s_i) = (1-t)a+tb$ if and only if $v_0(s_i)=a$ and $v_1(s_i)=b$, thus for these $s_i$'s we necessarily have $(1-t')v_0(s_i)+t'v_1(s_i) = (1-t')a+t'b$ for all $t'\in (0,1)$. For the remaining $s_i$'s, as $t'$ is sufficiently close to $t$, we may assume 
$$(1-t')v_0(s_i)+t'v_1(s_i) <\mbox{ (resp. $>$) }(1-t')a+t'b$$ if $(1-t)v_0(s_i)+tv_1(s_i) <$ (resp. $>$) $(1-t)a+tb$. Hence $\cF_t^{(1-t)a+tb} R_m$ and $\cF_{t'}^{(1-t')a+t'b} R_m$ are spanned by the same set of $s_i$'s and we get the first equality in \eqref{eq:F_t perturb}. The second equality in \eqref{eq:F_t perturb} can be proved in a similar fashion as in Lemma \ref{lem:complement for f.g. special valuation}, noting that as $t'$ stays close to $t$, the monomial that computes $v_t(s)$ (where $s\in R_m$) also computes $v_{t'}(s)$ (and there are only finitely many such monomials to consider). Thus since $\cF_t=\cF_{v_t}$ when $t\in \bQ$, it follows that $\cF_t=\cF_{v_t}$ for irrational $t$ as well. Since $\gr_{\cF_t} R\cong \gr_{v_0}(\gr_{v_1} R)$ and the right hand side is finitely generated by the above discussion, the same holds for $\cF_{v_t}$. The proof is now complete.
\end{proof}

\begin{lem} \label{lem:S linear imply f.g. general}
Let $0\le \Gamma\sim_\bQ -(K_X+\Delta)$ be a $\bQ$-complement of the log Fano pair $(X,\Delta)$, let $(Y,E)$ be a log smooth model of $(X,\Delta)$, and let $P$ be a convex subset in $\QM(Y,E)\cap\DMR(X,\Delta+\Gamma)$. Assume that the linear span of $P$ is rational, and $S_{X,\Delta}(\cdot)$ is linear on $P$. Then for any $v\in P$, the associated filtration $\cF_v$ is finitely generated.
\end{lem}

\begin{proof}
We first claim that for any valuations $v_0,v_1\in P$ and any $t\in (0,1)$, we have
\begin{equation} \label{eq:F_v when S linear}
\cF_{v_t}^\lambda R_m = {\rm Span}\{s\in R_m\,|\,(1-t)v_0(s)+t v_1(s)\ge \lambda\}
\end{equation}
where $v_t=(1-t)v_0+tv_1$. As in the proof of Lemma \ref{lem:S linear imply f.g. interval}, \eqref{eq:F_v when S linear} is equivalent to saying that for any basis $s_1,\cdots,s_{N_m}$ of $R_m$ that's compatible with both $v_0$ and $v_1$, we have
\begin{equation} \label{eq:v_t linear}
    v_t\left(\sum_{i=1}^{N_m} a_i s_i\right)=\min\{(1-t)v_0(s_i)+tv_1(s_i)\mid a_i\neq 0\}
\end{equation}
for any $a_1,\cdots,a_{N_m}\in \bk$. As in the proof of \eqref{eq:F_t perturb}, we know that the basis $s_1,\cdots,s_{N_m}$ is also compatible with divisorial valuations in the minimal rational PL subspace containing $v_0$ (resp. $v_1$) that are sufficiently close to $v_0$ (resp. $v_1$). By assumption, these divisorial valuations are necessarily contained in $P$. By Lemma \ref{lem:S linear imply f.g. interval}, the equality \eqref{eq:F_v when S linear} (and hence \eqref{eq:v_t linear}) holds when $v_0,v_1$ are divisorial. Thus as the function $v\mapsto v(s)$ is continuous for any fixed $s\in R_m$, we conclude that \eqref{eq:v_t linear} (and hence \eqref{eq:F_v when S linear}) holds for all $v_i\in P$ and $t\in (0,1)$. 

Now for any $v\in P$, we may find some $t\in (0,1)$ and some valuation $v_1'\in P$ such that $v_t'=(1-t)v+tv_1'$ is divisorial. By \eqref{eq:F_v when S linear} and the discussions in Lemma \ref{lem:S linear imply f.g. interval}, we see that $\gr_{v_1'}(\gr_v R)\cong \gr_{v_t'} R$. Since $v_t'$ is divisorial, the associated graded ring $\gr_{v_t'} R$ is finitely generated by Lemma \ref{l-fg for divisorial lc place}. It follows that $\gr_{v_1'}(\gr_v R)$ is finitely generated, hence the same holds for $\gr_v R$ by lifting the generators of $\gr_{v_1'}(\gr_v R)$. Thus $\cF_v$ is finitely generated as desired.
\end{proof}

We now present the proof of Theorem \ref{thm:f.g. criterion}.

\begin{proof}[Proof of Theorem \ref{thm:f.g. criterion}]
(1) $\Rightarrow$ (2) in fact in a suitable neighbourhood $U$ all $(X_v,\Delta_v)$ are isomorphic, see \cite{LX-higher-rank}*{Lemma 2.10} or Lemma \ref{l-finitegeneratedimpliesisom}. (3) $\Rightarrow$ (1) by Lemma \ref{lem:S linear imply f.g. general}. So it remains to prove (2) $\Rightarrow$ (3). The key point is that if a concave function takes rational values with linearly bounded denominators on rational points, then it is linear.

Denote by $(Y,E)$ the log smooth model of $(X,\Delta)$ that appears in the definition of $\Sigma$. 
By Lemma \ref{lem:S linear imply f.g. interval}, the $S$-invariant function is concave on the simplex of $\QM(Y,E)$ containing $v_0$. In particular, it is Lipschitz in a neighborhood of $v_0$ in $\Sigma$, and we may find some constant $C>0$ such that 
\begin{equation} \label{eq:lipschitz}
    |S_{X,\Delta}(v)-S_{X,\Delta}(v_0)|\le C|v-v_0|
\end{equation}
for any $v$ in a neighbourhood of $v_0$ in $\Sigma$. Let $S_0=S_{X,\Delta}(v_0)$. By \cite{LX-higher-rank}*{Lemma 2.7}, for any $\varepsilon>0$ there exists divisorial valuations $v_1,\cdots,v_\ell\in \Sigma$, rational numbers $S_1,\cdots,S_\ell$ and positive integers $q_1,\cdots,q_\ell$ such that 
\begin{itemize}
    \item $(v_0,S_0)$ is in the convex cone generated by $(v_i,S_i)$, i.e., there exists some $\lambda_i>0$ such that $v_0=\sum_{i=1}^\ell \lambda_i v_i$ and $S_0=\sum_{i=1}^\ell \lambda_i S_i$.
    \item $(q_i v_i,q_i S_i)$ is an integer vector for all $i=1,\cdots,\ell$.
    \item $|v_i-v_0|+|S_i-S_0|<\frac{\varepsilon}{q_i}$ for all $i=1,\cdots,\ell$.
\end{itemize}
In particular, by the last condition, we may assume that $v_i\in U$, where $U\subseteq \Sigma$ is the open neighbourhood of $v_0$ in the condition (2). 

Since the set $\{(X_v,\Delta_v)\,|\,v\in U(\bQ)\}$ is bounded, there exists some integer $M>0$ such that $M(K_X+\Delta)$ is Cartier and that $M\cdot \Fut_{X_v,\Delta_v}(\xi)\in \bZ$ for any $v\in U(\bQ)$ and any one parameter subgroup $\xi\colon \bG_m\to \Aut(X_v,\Delta_v)$. If we let $q$ be an integer such that $qv$ is integral, then $Mq\cdot A_{X,\Delta}(v)\in\bZ$. Note that $qv$ induces a one parameter subgroup $\xi_v\colon \bG_m\to \Aut(X_v,\Delta_v)$ with $\Fut_{X_v,\Delta_v}(\xi_v)=q\cdot \beta_{X,\Delta}(v)$. It follows that $Mq\cdot S_{X,\Delta}(v)\in \bZ$ for any $v\in U(\bQ)$ and any integer $q$ such that $qv$ is integral. In particular, we have $Mq_i\cdot  S_{X,\Delta}(v_i)\in \bZ$. 

On the other hand, by \eqref{eq:lipschitz} we have $|S_{X,\Delta}(v_i)-S_0|\le \frac{C\varepsilon}{q_i}$ and hence 
\[
|Mq_i\cdot S_{X,\Delta}(v_i) - Mq_i S_i|\le Mq_i\cdot |S_{X,\Delta}(v_i)-S_0| + Mq_i \cdot |S_i-S_0| \le (C+1)M\varepsilon.
\]
Note that the constants $C$ and $M$ are independent of the choice of $\varepsilon$. Thus if we take $\varepsilon=\frac{1}{2(C+1)M}$, then as $Mq\cdot S_{X,\Delta}(v)$ and $q_i S_i$ are both integers we deduce that $S_{X,\Delta}(v_i)=S_i$. But as $S_{X,\Delta}(\cdot)$ is concave on $U$ and $v_0=\sum_{i=1}^\ell \lambda_i v_i$, we also have
\[
\sum_{i=1}^\ell \lambda_i S_{X,\Delta}(v_i)\le  S_{X,\Delta}(v_0)=S_0=\sum_{i=1}^\ell \lambda_i S_i = \sum_{i=1}^\ell \lambda_i S_{X,\Delta}(v_i).
\]
Hence the first inequality is an equality, which forces $S_{X,\Delta}(\cdot)$ to be linear on the cone generated by $v_1,\cdots,v_\ell$. In particular, it is linear in a neighbourhood of $v_0$ in $\Sigma$.
\end{proof}

\subsection{Estimate of alpha invariants}
We next proceed to check the condition of Theorem \ref{thm:f.g. criterion}(2) when the complement is special. In order to control the boundedness of $(X_v,\Delta_v)$ for $v\in \Sigma(\mathbb Q)$, we wish to apply the boundedness result as in Theorem \ref{t-alphabounded}. In light of Lemma \ref{lem:special DMR},  we already know $(X_v,\Delta_v)$ is klt.
Then we need to further analyze the $\alpha$-invariant of $(X_v,\Delta_v)$.

The following theorem is our main result in this section.
\begin{thm}\label{t-alphabelow}Let $(X,\Delta)$ be a log Fano pair and $\Gamma$ a special complement with respect to a log smooth model $(Y,E)$ (see Definition \ref{d-specialcomplement}). Let $K\subset \DMR(X,\Delta+\Gamma)$ be a compact subset that is contained in the interior of a simplicial cone in $\QM(Y,E)$. Then there exists some constant $\alpha>0$ such that for all rational points $v\in K$, the alpha invariants $\alpha(X_v,\Delta_v)$ of the induced degenerations $(X_{v},\Delta_v)$ are bounded from below by $\alpha$.  
\end{thm}
Our main tool is the following characterization of $\alpha$-invariants.

\begin{lem} \label{lem:(alpha,v)-complement}
Let $v$ be a divisorial valuation such that $\gr_v R$ is finitely generated and let $\alpha\in(0,1)$ be a rational number. Then $\alpha(X_v,\Delta_v)\ge \alpha$ if and only if for all $0\le D\sim_\bQ -(K_X+\Delta)$, there exists some $0\le D'\sim_\bQ -(K_X+\Delta)$ such that $(X,\Delta+\alpha D+(1-\alpha)D')$ is lc and has $v$ as an lc place.
\end{lem}

For ease of notation, we call such $D'$ {\it an $(\alpha,v)$-complement of $D$}.

\begin{proof}
Note that $(X_v,\Delta_v)$ has an induced $\bG_m$-action. By taking the limit under the $\bG_m$-action, we see that any effective $\bQ$-divisor $G\sim_\bQ -(K_{X_v}+\Delta_v)$ degenerates to some $\bG_m$-invariant divisor $G_0$. By the semi-continuity of log canonical thresholds, we have $\lct(X_v,\Delta_v;G)\ge \lct(X_v,\Delta_v;G_0)$, hence $\alpha(X_v,\Delta_v)\ge \alpha$ if and only if $\lct(X_v,\Delta_v;G_0)\ge \alpha$ for all $\bG_m$-invariant divisors $G_0\sim_\bQ -(K_{X_v}+\Delta_v)$. Any such $G_0$ is also the specialization of some divisor $0\le D\sim_\bQ -(K_X+\Delta)$ on $X$, and $\lct(X_v,\Delta_v;G_0)\ge \alpha$ means that $v$ induces a weakly special degeneration of $(X,\Delta+\alpha D)$. By Theorem \ref{t-complements}, this is the case if and only if for all $\varepsilon\in\bQ$, $0<\varepsilon\ll 1$, the valuation $v$ is an lc place of a $\bQ$-complement of the klt pair $(X,\Delta+(\alpha-\varepsilon) D)$, i.e. $D$ has an $(\alpha-\varepsilon,v)$-complement (note that $(X,\Delta+\alpha D)$ is only lc so we cannot directly apply Theorem \ref{t-complements}). We claim that this is equivalent to saying that $D$ has an $(\alpha,v)$-complement. Clearly only the forward direction needs a proof. We may write $v=c\cdot \ord_E$. Let $a=A_{X,\Delta}(E)$. Since $v$ is an lc place of $\bQ$-complement, there exists a birational model $\pi\colon Y\to X$ that extracts $E$ as the only possible exceptional divisor (by \cite{BCHM}*{Corollary 1.4.3}), and $Y$ is of Fano type. Moreover, it follows from the existence of $(\alpha-\varepsilon,v)$-complement that the pair $(Y,\pi^{-1}_*(\Delta+(\alpha-\varepsilon)D)\vee E)$ has a $\bQ$-complement for all $0<\varepsilon\ll 1$. By Lemma \ref{lem:ACC}, this implies that $(Y,\pi^{-1}_*(\Delta+\alpha D)\vee E)$ also has a $\bQ$-complement, whose pushforward to $X$ is an $(\alpha,v)$-complement of $D$. This proves the claim and also the statement of the lemma.
\end{proof}

\begin{cor} \label{cor:alpha<A/T}
Let $v$ be a divisorial valuation which is an lc place of a $\mathbb Q$-complement of $(X,\Delta)$. Then
\[
\alpha(X_v,\Delta_v)\le 1- \frac{A_{X,\Delta}(v)}{T_{X,\Delta}(v)}.
\]
\end{cor}

\begin{proof}
We may assume that $\alpha(X_v,\Delta_v)>0$, otherwise since $v$ is an lc place of some $\bQ$-complement we have $T_{X,\Delta}(v)\ge A_{X,\Delta}(v)$ and the result is clear. Let $\alpha=\alpha(X_v,\Delta_v)$. Note that $\alpha(X_v,\Delta_v)<1$ since otherwise it can not have any non-trivial $\bG_m$-action by \cite{LZ-alpha-sharpness}*{Corollary 3.6} (the same proof in \emph{loc. cit.} works for pairs). Choose some effective divisor $D\sim_\bQ -(K_X+\Delta)$ whose support does not contain $C_X(v)$. By Lemma \ref{lem:(alpha,v)-complement}, there exists some $0\le D'\sim_\bQ -(K_X+\Delta)$ such that $v$ is an lc place of $(X,\Delta+\alpha D+(1-\alpha)D')$. In particular, $(1-\alpha)v(D')=v(\alpha D+(1-\alpha)D')=A_{X,\Delta}(v)$, which implies $A_{X,\Delta}(v)\le (1-\alpha)T_{X,\Delta}(v)$. In other words, $\alpha(X_v,\Delta_v)\le 1- \frac{A_{X,\Delta}(v)}{T_{X,\Delta}(v)}$.
\end{proof}

In order to construct $(\alpha,v)$-complements for some uniform constant $\alpha$ (and therefore produce a uniform lower bound for $\alpha(X_v,\Delta_v)$ by Lemma \ref{lem:(alpha,v)-complement}), our strategy is to refine the proof of Lemma \ref{lem:special DMR} using the alpha invariants and nef thresholds of the corresponding exceptional divisors $F$. This is done in the next three lemmas. 

To this end, we introduce some notation. Under the notation and assumptions of Theorem \ref{t-alphabelow}, we fix an effective ample $\bQ$-divisor $G$ on $Y$ that does not contain any stratum of $E$ such that $\Gamma_Y\ge G$. For any divisorial valuation $v\in \DMR(X,\Delta+\Gamma)\cap \QM(Y,E)$, let $\rho\colon Z\to Y$ be the corresponding weighted blowup, $F$ the exceptional divisor (i.e. $v=c\cdot \ord_F$), and $(Z,\Delta_Z)$, $(Y,\Delta_Y)$ the crepant pullbacks as in the proof of Lemma \ref{lem:special DMR}. Let $\Delta^+:=\Delta_Z\vee 0\vee F$. Note that $(Z,\Delta^+)$ is plt. By adjunction we may write $K_F+\Phi=(K_Z+\Delta^+)|_F$. Let 
$$L:=-\rho^*\pi^*(K_X+\Delta)-A_{X,\Delta}(F)\cdot F.$$ 
Since $v=c\cdot \ord_F$ is an lc place of $(X,\Delta+\Gamma)$, $F$ is not contained in the support of $\rho^*\pi^*\Gamma-A_{X,\Delta}(F)\cdot F\sim_\bQ L$, thus the $\bQ$-linear system $|L_{|F}|_\bQ$ is non-empty and we may define
\[
\alpha_v:=\lct (F,\Phi; |L_{|F}|_{\mathbb Q}).
\]
We also let 
\[
\varepsilon_v:=\sup\{t\ge 0\,|\,\rho^*G-tA_{X,\Delta}(F)\cdot F \mbox{ is nef}\}.
\]
Note that as $-F$ is $\rho$-ample, we have $\varepsilon_v>0$ and for any $t\in (0,\varepsilon_v)$, the divisor $\rho^*G-tA_{X,\Delta}(F)\cdot F$ is ample.

\begin{lem} \label{lem:alpha bound degeneration}
Let $a,b>0$ be constants. Then there exists some constant $\alpha>0$ depending only on $a,b,(X,\Delta)$ and $\Gamma$ such that $\alpha(X_v,\Delta_v)\ge \alpha$ as long as $\alpha_v>a$ and $\varepsilon_v>b$.
\end{lem}

\begin{proof}
We may assume that $a<1$. By Lemma \ref{lem:(alpha,v)-complement}, it is enough to find some constant $\alpha>0$ such that $(\alpha,v)$-complement exists for any $0\le D\sim_\bQ -(K_X+\Delta)$. 

As a reduction step, we first show that it suffices to check the existence of $(\alpha,v)$-complement for those $D$ such that $v(D)=A_{X,\Delta}(v)$. Indeed, as $G$ is ample, we may find some constant $0<\lambda\ll 1$ such that $G+\lambda\pi^*(K_X+\Delta)$ remains ample. It follows that $T(G;v)\ge \lambda\cdot T_{X,\Delta}(v)$, thus
\begin{align*}
    T_{X,\Delta}(v)=T(\pi^*\Gamma;v) & \ge v(\pi^*\Gamma-G)+T(G;v) \\
    & \ge v(\Gamma)+\lambda\cdot T_{X,\Delta}(v)=A_{X,\Delta}(v)+\lambda\cdot T_{X,\Delta}(v),
\end{align*}
or $(1-\lambda)T_{X,\Delta}(v)\ge A_{X,\Delta}(v)$. On the other hand, $\alpha(X,\Delta)T_{X,\Delta}(v)\le A_{X,\Delta}(v)$ by the definition of alpha invariants. 

We claim there exists some constant $\mu\in(0,1)$ depending only on $\lambda$ and $\alpha(X,\Delta)$ such that for any $0\le D\sim_\bQ -(K_X+\Delta)$, we can always find some $0\le D_1\sim_\bQ -(K_X+\Delta)$ and $\nu\ge \mu$ such that $\nu v(D)+(1-\nu) v(D_1)=A_{X,\Delta}(v)$. In fact, since the possible values of $v(D_1)$ are dense between $0$ and $T_{X,\Delta}(v)$, this is a consequence of the following elementary fact: if $0<\alpha\le \frac{A}{T}\le 1-\lambda$, then there exists $\mu>0$ depending only on $\alpha$ and $\lambda$, such that for any $p\in [0,T]$, we can find some $q\in [0,T)$ and some $\nu\ge \mu$ such that $\nu p+(1-\nu)q=A$. 

Clearly if $(\alpha,v)$-complement exists for $\nu D+(1-\nu)D_1$, then $(\alpha\mu,v)$-complement exists for $D$. This proves the reduction.

Next we fix a sufficiently small $t\in(0,1)$ such that $s:=\frac{(1-a)t}{1-t}<b$. By assumption, $\rho^*G-sA_{X,\Delta}(F)\cdot F$ is ample. Fix any $0\le D\sim_\bQ -(K_X+\Delta)$ with $v(D)=A_{X,\Delta}(v)$, let $H'$ be a general member of the $\bQ$-linear system $|\rho^*G-sA_{X,\Delta}(F)\cdot F|_\bQ$, and let $H=\rho_*H'$. 

We claim that along $\rho(F)$ the pair 
\[
(Y,\Delta_Y+a\cdot\pi^*D+\frac{1-t}{t}H)
\]
is lc and has $F$ as its unique lc place. To see this, first we have
$$A_{Y,\Delta_Y}(F)-\ord_F(a\cdot\pi^*D+\frac{1-t}{t}H)=A_{X,\Delta}(F)-aA_{X,\Delta}(F)-(1-a)A_{X,\Delta}(F)=0.$$
Then let 
\[
D'=\rho^*\pi^*D-\ord_F(D)\cdot F=\rho^*\pi^*D-A_{X,\Delta}(F)\cdot F\sim_\bQ L.
\]
By assumption, $(F,\Phi+aD'|_F)$ is klt, hence since $H'$ is general, we see that $(F,\Phi+aD'|_F+\frac{1-t}{t}H'|_F)$ is also klt. By inversion of adjunction, $(Z,\Delta^+ +aD'+\frac{1-t}{t}H')$ is plt along $F$. Since $\Delta^+\ge \Delta_Z\vee F$, we deduce that $(Z,\Delta_Z\vee F+aD'+\frac{1-t}{t}H')$ is also plt along $F$. By construction and the above calculation, we can check that
\[
K_Z+\Delta_Z\vee F+aD'+\frac{1-t}{t}H'=\rho^*(K_Y+\Delta_Y+a\cdot\pi^*D+\frac{1-t}{t}H),
\]
thus $(Y,\Delta_Y+a\cdot\pi^*D+\frac{1-t}{t}H)$ is lc along $\rho(F)$ and $F$ is the only lc place there, proving the previous claim. 

We also know that $(Y,\Delta_Y+\pi^*\Gamma-G)$ is lc and $F$ is an lc place of the pair. Taking convex linear combination as in the proof of Lemma \ref{lem:special DMR}, it follows that 
\begin{align*}
    & \left(Y,\Delta_Y+t\left(a\cdot\pi^*D+\frac{1-t}{t}H\right)+(1-t)(\pi^*\Gamma-G)\right) \\
    & = \big(Y,\Delta_Y+at\cdot \pi^*D+(1-t)(\pi^*\Gamma-G+H)\big)
\end{align*}
is lc along $\rho(F)$ and $F$ is the only lc place of the pair in a neighbourhood of $\rho(F)$. In particular, $\rho(F)$ is a connected component of the non-klt locus of the pair. Note that
\[
K_Y+\Delta_Y+at\cdot\pi^*D+(1-t)(\pi^*\Gamma-G+H)=\pi^*(K_X+\Delta+atD+(1-t)(\Gamma-\pi_*G+\pi_*H)),
\]
thus $(X,\Delta+atD+(1-t)(\Gamma-\pi_*G+\pi_*H))$ is lc along $C_X(v)=\pi(\rho(F))$ and $C_X(v)$ is a connected component of its non-klt locus, since otherwise in some neighbourhood of $\pi^{-1}C_X(v)$ there would be another non-klt center of $(Y,\Delta_Y+at\cdot\pi^*D+(1-t)(\pi^*\Gamma-G+H))$ that is disjoint from $\rho(F)$, contradicting the Koll\'ar-Shokurov connectedness theorem (see e.g. \cite{K+92}*{17.4 Theorem}). 

Similarly, as
\[
-(K_X+\Delta+atD+(1-t)(\Gamma-\pi_*G+\pi_*H))\sim_\bQ -(1-a)t(K_X+\Delta)
\]
is ample, we deduce that $(X,\Delta+atD+(1-t)(\Gamma-\pi_*G+\pi_*H))$ is indeed lc everywhere, as otherwise there would be some non-klt center of the pair that is disjoint from $C_X(v)$, contradicting Koll\'ar-Shokurov connectedness. Since by construction $v=c\cdot\ord_F$ is an lc place of $(X,\Delta+atD+(1-t)(\Gamma-\pi_*G+\pi_*H))$, we may add some general divisor $0\le D_1\sim_\bQ -(1-a)t(K_X+\Delta)$ to the pair and conclude that $D$ has an $(at,v)$-complement. Since $D$ is arbitrary and $t$ only depends on $a,b$, this completes the proof.
\end{proof}

The argument for the following lemma is similar to the one in \cite{Z-product}.

\begin{lem} \label{lem:alpha bound plt blowup}
Notation as above. Let $K\subseteq\DMR(X,\Delta+\Gamma)$  be a compact subset that is contained in the interior of some simplicial cone in $\QM(Y,E)$. Then there exist some constants $a>0$ such that $\alpha_v\ge a$ for all divisorial valuations $v\in K$.
\end{lem}

\begin{proof}
Let $E_i$ ($i=1,\cdots,r$) be the irreducible components of $E$ so that $W=E_1\cap \cdots \cap E_r$ is the common center of valuations in $K$ on $Y$. Any divisorial valuation $v\in K$ corresponds to a weighted blowup at $W$ with weights $\wt(E_i)=a_i$ for some integers $a_i>0$ such that $\gcd(a_1,\cdots,a_r)=1$. Since $K$ is compact, there exists some constant $C>0$ such that $\frac{a_i}{a_j}<C$ for all $1\le i,j\le r$.

To describe the weighted blow up in more details, if locally around a point $x\in W$, $E_i$ is given by the equation $e_i=0$, then the weighted blow up is given by 
\begin{equation}\label{e-weightedblowup}
{\rm Proj}_{\mathcal{O}_Y}(\mathcal{O}_Y\oplus \cI_1\oplus \cI_2\oplus\cdots),
\end{equation}
where $\cI_d$ around $x$ is generated by monomials $ e_1^{d_1}\cdots e_r^{d_r}$ such that $\sum a_id_i\ge d$. The exceptional divisor $F$ is a weighted projective space bundle over $W$ with fiber $F_0$ isomorphic to $\mathbb{A}^r\setminus \{0\}/\mathbb {G}_m$ with the action $\lambda\cdot (y_1,...,y_r)=(\lambda^{a_1}y_1,..., \lambda^{a_r}y_r)$.
The corresponding well-formed weighted projective space is given in the following way: let $q_i=\gcd(a_1,\cdots,\hat{a}_i,\cdots,a_r)$, $q=q_1\cdots q_r$, and $a'_i=\frac{a_i q_i}{q}$, then $F_0\cong \bP(a'_1,\cdots,a'_r)$ which is well-formed (see \cite[Lemma 5.7]{Fletcher-weighted}).

Let $c_i=A_{X,\Delta}(E_i)>0$, $b_i=\max\{0,\ord_{E_i}(\Delta_Y)\}<1$. Then we have the following facts (with notation as in the paragraph before Lemma \ref{lem:alpha bound degeneration}):
\begin{enumerate}
    \item $A_{X,\Delta}(F)=\sum_{i=1}^r a_i A_{X,\Delta}(E_i) =a_1 c_1+\cdots+a_r c_r$,
    \item $L_{F_0}:=L|_{F_0}\sim_\bQ \frac{A_{X,\Delta}(F)}{q} L_0$ where $L_0$ is the class of $\cO(1)$ on $\bP(a'_1,\cdots,a'_r)$,
    \item $\Phi_{F_0}:=\Phi|_{F_0}=\sum_{i=1}^r \frac{q_i-1+b_i}{q_i}\{x_i=0\}$ where $x_1,\cdots,x_r$ are the weighted homogeneous coordinates on $\bP(a'_1,\cdots,a'_r)$,
    \item $\fb_m:=\rho_*\cO_Z(-mF)/\rho_*\cO_Z(-(m+1)F)\cong \bigoplus \cO_W(-(m_1 E_1+\cdots+m_r E_r))$, where the direct sum runs over all $(m_1,\cdots,m_r)\in\bN^r$ such that $a_1 m_1+\cdots a_r m_r=m$,
\end{enumerate}

\begin{rem}
All these facts can easily be seen if we view the weighted blow up as a (Deligne-Mumford) stack (see e.g. \cite[Section 3]{ATW19}). Then the stacky exceptional divisor $\cF$ is a weighted projective stack bundle over $W$ with the fiber $\cP(a_1,...,a_r):=[(\bA^{r}\setminus \{0\})/\bG_m]$ as a Deligne-Mumford stack where $\bG_m$ acts diagonally on $\bA^r$ by weights $(a_1,\cdots, a_r)$ (see e.g. \cite[Section 2.3]{RT-weighted}), and $\mathcal{O}(-\cF)|_{\cP(a_1,...,a_r)}\cong \cO(1)$. The natural morphism 
$$\cP(a_1,...,a_r) \to \bP(a'_1,\cdots,a'_r)$$
maps the Deligne-Mumford stack to its coarse space, with an orbifold divisor of the form  $(1-\frac{1}{q_i})(x_i=0)$ along the $i$-th coordinate hyperplane $(x_i=0)$ of $\bP(a'_1,\cdots,a'_r)$ (see \cite[Section 2.5]{RT-weighted}). 

Now the claim (2) follows from the fact that $\mathcal{O}(-\cF)$ is $\mathcal{O}(1)$ on $\cP(a_1,...,a_r)$, and the pull back of $\cO(1)$ on $\bP(a'_1,\cdots,a'_r)$  to $\cP(a_1,...,a_r)$ is $\mathcal{O}(q)$. In fact, for the map $$\cP(\lambda_1,...,\lambda_i,...,\lambda_n)\to \cP(\frac{\lambda_1}{\lambda},...,\lambda_i,...,\frac{\lambda_n}{\lambda})\qquad \mbox{where $\lambda=\gcd(\lambda_1,...,\hat{\lambda}_i,...,\lambda_r)$},$$
the pull back of the divisor given by the section $x_i$ of $\mathcal{O}_{\cP(\frac{\lambda_1}{\lambda},...,\lambda_i,...,\frac{\lambda_n}{\lambda})}(\lambda_i)$ is $\lambda$-multiple of the orbifold divisor given by the section $y_i$ of $\mathcal{O}_{\cP(\lambda_1,...,\lambda_i,...,\lambda_n)}(\lambda_i)$. 
For (3), the birational transform of $E_i$ restricting on $\cP(a_1,...,a_r)$ is $\{y_i=0\} $ on $\cP(a_1,...,a_r)$ with coordinates $y_i$, then 
$K_{\cP(a_1,...,a_r)}+\sum b_i\{y_i=0\}$ is the pull back of $K_{\bP(a'_1,\cdots,a'_r)}+\sum_{i=1}^r \frac{q_i-1+b_i}{q_i}\{x_i=0\}$ as the morphism $\cP(a_1,...,a_r)\to \bP(a'_1,\cdots,a'_r)$ precisely has codimension one orbifold components along $(y_i=0)$ with degree $q_i$ (alternatively, (3) can also be derived through a direct calculation of the different, noting that $Z$ has a cyclic quotient singularity of order $q_i$ along the codimension $2$ component $\{x_i=0\}\subseteq F$). Finally, let $\fc_m=\bigoplus \cO_W(-(m_1 E_1+\cdots+m_r E_r))$ be the direct sum that appears in (4). 
Since $\rho^*(\cO_Y(-(m_1 E_1+\cdots+m_r E_r))\subset \cO_Z(-mF)$ if $\sum a_im_i= m$, 
then $$\cO_Y(-(m_1 E_1+\cdots+m_r E_r))\subset \rho_*\cO_Z(-mF)=\cI_m.$$
Therefore, we have a natural map $\fc_m\to \fb_m$. A local computation shows that both two sides are generated as free module by the image of $e_1^{m_1}\cdots e_r^{m_r}$ with $a_1m_1+\cdots+a_rm_r=m$. Thus (4) follows.
\end{rem}

By (1) and (4),  we have for any $m\in\bN$ such that $mL$ is Cartier, 
\begin{align*}
\rho_*\cO_{F}(mL)&\cong\cO_Y(-m\pi^*(K_X+\Delta))\otimes \rho_*\cO_Z(-mA_{X,\Delta}(F)\cdot F)/\rho_*\cO_Z(-(mA_{X,\Delta}(F)+1)F)\\
&\cong \bigoplus_{(m_1,...,m_r)} \cO_W\big(-m\pi^*(K_X+\Delta)-(m_1 E_1+\cdots+m_r E_r)\big),    
\end{align*}
the direct sum running over all $\sum^r_{i=1} a_im_i=m\sum^r_{i=1} a_ic_i$. Recall that $\frac{a_i}{a_j}<C$ for all $1\le i,j\le r$, thus $m_1+\cdots+m_r\le C_0 m$ where $C_0=\lceil C\sum_{i=1}^r c_i\rceil$. If we choose a very ample line bundle $H_0$ such that $H_0+E_i$ ($1\le i\le r$) are all very ample and $H_0+\pi^*(K_X+\Delta)$ is ample, then for sufficiently divisible $m$, each direct summand in $\rho_*\cO_F(mL)$ admits an inclusion 
\begin{align*}
    \cO_W\big(-m\pi^*(K_X+\Delta)-(m_1 E_1+\cdots+m_r E_r)\big) & \hookrightarrow \cO_W\big((m+m_1+\cdots+m_r)H_0\big) \\
    & \hookrightarrow \cO_W\big((C_0+1)mH_0\big).
\end{align*}
Therefore for $H=(C_0+1)H_0$ and sufficiently divisible $m$, we have
$${\rho}_*\cO_{F}(mL) \hookrightarrow \cO_W(mH)^{\oplus N_m}$$ for some integer $N_m$ ($={\rm rank}(\fb_{mA_{X,\Delta}(F)})$). 

Since $F_0$ is toric, by \cite{BJ-delta}*{Theorem F} we know that $\lct(F_0,\Phi_{F_0};|L_{F_0}|_\bQ)$ is computed by one of torus invariant divisors $\{x_i=0\}$, thus by (1) (2) and (3), we get
\[
\lct(F_0,\Phi_{F_0};|L_{F_0}|_\bQ)=\frac{\min_{1\le i\le r} a_i(1-b_i)}{a_1 c_1+\cdots+a_r c_r} \ge a
\]
for some constants $a>0$ that only depend on $b_i$, $c_i$ and $C$. Let $\Delta_W:=(\Delta_Y\vee 0 - \sum_{i=1}^r b_i E_i)|_W$. Note that $(W,\Delta_W)$ is klt and $\rho^*\Delta_W$ is the vertical part of $\Phi$. We have $\lct(W,\Delta_W;|H|_\bQ)>0$ by Izumi's inequality (see e.g. \cite{Li-volume}*{Theorem 1.2} or \cite{BL-vol-lsc}*{Theorem 20}), since there exists constant $M>0$ such that $\mult_w D\le M$ for any effective divisor $D\sim_\bQ H$ and any $w\in W$. Thus after replacing $a$ by a smaller positive number, we may further assume that $\lct(W,\Delta_W;|H|_\bQ) \ge a$. We claim that $\alpha_v\ge a$.

To see this, let $t\in(0,a)$ and let $\Phi'\sim_\bQ L|_F$ be an effective divisor. Suppose that $(F,\Phi+t\Phi')$ is not lc. Then since $(F,\Phi+t\Phi')$ is lc along the general fiber of $F\to W$ by our choice of $a$, we know that there exists a divisorial valuation $v_0$ over $F$ such that $A_{F,\Phi+t\Phi'}(v_0)<0$ and the center of $v_0$ does not dominate $W$. By \cite{Z-product}*{Lemma 2.1}, $v_0$ restricts to a divisorial valuation $w$ on $W$. 

Let $g\colon W_1\to W$ be a birational morphism such that the center of $w$ is a divisor $Q$ on $W_1$, let $F_1=F\times_W W_1$, $\Phi_1=g^*(\Phi-\rho^*\Delta_W)$ (i.e. the pullback of the horizontal part of $\Phi$; here we also denote the projection $F_1\to F$ by $g$), and let $P$ be the preimage of $Q$ in $F_1$. 

Since $F\to W$ is locally a trivial product $F_0\times W$, it is not hard to see that the formation of $\rho_*\cO_F(mL)$ commutes with the base change $W_1\to W$. Thus by projection formula, we see that
\[
H^0(F_1,\cO_{F_1}(g^*mL-kP))=H^0(W_1,g^*\rho_*\cO_{F}(mL)\otimes \cO_{W_1}(-kQ)).
\]
Since ${\rho}_*\cO_{F}(mL)\hookrightarrow \cO_W(mH)^{\oplus N_m}$ for any sufficiently divisible $m$, we have
$$H^0(F_1,\cO_{F_1}(g^*mL-kP))\neq 0 \Rightarrow H^0(W_1,\cO_{W_1}(mg^*H-kQ))\neq 0$$
for any $k\in \bN$, which implies 
\begin{equation*}\label{eq-comparemult}
\ord_P(\Phi')\le \sup_{ H'\in |H|_{\mathbb Q} }\ord_Q(H').
\end{equation*}
Since  $\lct(W,\Delta_W;|H|_\bQ) \ge a$,
 \begin{equation*}\label{eq-TinvariantforW}
a\cdot \sup_{ H'\in |H|_{\mathbb Q} }\ord_Q(H')\le A_{W,\Delta_W}(Q)=A_{F,\Phi}(P).
 \end{equation*} 
Since $t<a$, we combine the above inequalities to conclude that $t\cdot \ord_P(\Phi')<A_{F,\Phi}(P)$. It follows that if we write 
\[
g^*(K_F+\Phi+t\Phi')=K_{F_1}+\Phi_1+\lambda P+D
\]
where $P\not\subseteq \Supp(D)$, then $\lambda\le 1$. 

On the other hand, since the divisor $P$ is vertical, over a general fiber of $P\to Q$ (which we still denote by $F_0$), we have $D|_{F_0}\sim_\bQ tg^*\Phi'|_{F_0}\sim_\bQ tL|_{F_0}$, thus by our choice of $a$, $(P,(\Phi_1+D)|_P)$ is lc along the general fibers of $P\to Q$, hence by inversion of adjunction we see that $(F_1,\Phi_1+\lambda P+D)$ is also lc along the general fibers of $P\to Q$. In particular, it is lc at the center of $v_0$, a contradiction. Thus $(F,\Phi+t\Phi')$ is always lc and $\alpha_v\ge a$ as desired.
\end{proof}

\begin{lem} \label{lem:seshadri bound}
Notation as above. Let $K\subseteq\DMR(X,\Delta+\Gamma)$ be a compact subset that is contained in the interior of some simplicial cone in $\QM(Y,E)$. Then there exist some constants $b>0$ such that $\varepsilon_v\ge b$ for all divisorial valuations $v\in K$.
\end{lem}

\begin{proof}
We continue to use the notation from the proof of Lemma \ref{lem:alpha bound plt blowup}. Let 
\[
\fa_m:=\rho_*\cO_Z(-mA_{X,\Delta}(F)\cdot F).
\]
Since $\frac{a_i}{a_j}<C$ for all $1\le i,j\le r$, there exists some constant $M\in\bN$ such that 
\[
\frac{1}{A_{X,\Delta}(F)}\ord_{F}(f)\ge \frac{1}{M}\mult_W(f)
\]
for all regular function $f$ around the generic point of $W$. In particular, $\cI^{Mm}_W\subseteq \fa_m$ for all $m\in\bN$. 

\begin{claim}
We can find a sequence of ideals 
\begin{equation} \label{eq:ideal filtration}
    \cO_Y\supseteq \cI_W\supseteq\cdots\supseteq \fa_m\supseteq\cdots\supseteq \cI^{Mm}_W
\end{equation}
on $Y$ such that the quotients of consecutive terms are all isomorphic to $\cO_W(-n_1 E_1-\cdots-n_r E_r)$ for some $(n_1,\cdots,n_r)\in\bN^r$ with $\sum^r_{i=1} n_i < Mm$. 
\end{claim}

The key point in this claim is that $\fa_m$ appears in this sequence (the remaining terms in the middle do not matter too much to us).

\begin{proof}
We know 
$$\cI_W^{p}/\cI_W^{p+1}={\rm Sym}^p(\cI_W/\cI_W^{2})\cong \bigoplus_{n_1+\cdots +n_r=p, (n_1,...,n_r)\in \bN^r} \mathcal{O}_W(-n_1E_1-\cdots -n_rE_r).$$
Therefore we can find a sequence of ideals
 $$\cJ_0(=\cO_Y)\supseteq \cJ_1(=\cI_W)\supseteq\cdots \cJ_k\supseteq\cdots\supseteq \cI^{Mm}_W,$$
such that $\cJ_k/\cJ_{k+1}\cong \mathcal{O}_W(-n_1E_1-\cdots -n_rE_r)$ for some $n_1+\cdots +n_r<Mm, (n_1,...,n_r)\in \bN^r$ and all such $(n_1,...,n_r)$ appears precisely once as a subquotient. Indeed, if we order $\bN^r$ such that $(n_1,\cdots,n_r)\prec (n'_1,\cdots,n'_r)$ if and only if $n_1+\cdots+n_r<n'_1+\cdots+n'_r$ or $n_1+\cdots+n_r=n'_1+\cdots+n'_r$ and $(n_1,\cdots,n_r)<(n'_1,\cdots,n'_r)$ in the lexicographic order, then we can choose $\cJ_k$ so that around any $x\in W$ it is locally generated by the monomials $e_1^{d_1}\cdots e_r^{d_r}$ where $e_i=0$ is the local equation of $E_i$ and $(d_1,\cdots,d_r)\in\bN^r$ is at least the $(k+1)$-th smallest under the above ordering. We claim,
\begin{equation*}
  (\fa_m\cap \cJ_k)/ (\fa_m\cap \cJ_{k+1})=
    \begin{cases}
      \mathcal{O}_W(-n_1E_1-\cdots -n_rE_r) & \text{if $a_1n_1+\cdots +a_rn_r\ge mA_{X,\Delta}(F)$ }\\
       0               & \text{if $a_1n_1+\cdots +a_rn_r < mA_{X,\Delta}(F)$}\ \ .
          \end{cases}       
\end{equation*}
In fact, $(\fa_m\cap \cJ_k)/ (\fa_m\cap \cJ_{k+1})$ is isomorphic to the image of $\fa_m\cap \cJ_k$ in $\cJ_k\to \cJ_k/\cJ_{k+1}$, and a local calculation then gives that the image is 0 if $e_1^{n_1}\cdots e_r^{n_r}\not\in \fa_m$ and  $\mathcal{O}_W(-n_1E_1-\cdots -n_rE_r)$ if $e_1^{n_1}\cdots e_r^{n_r}\in \fa_m$.

Similarly, 
 \begin{equation*}
  (\fa_m+ \cJ_k)/ (\fa_m+ \cJ_{k+1})=
    \begin{cases}
      0 & \text{if $a_1n_1+\cdots +a_rn_r\ge mA_{X,\Delta}(F)$ }\\
       \mathcal{O}_W(-n_1E_1-\cdots -n_rE_r)             & \text{if $a_1n_1+\cdots +a_rn_r< mA_{X,\Delta}(F)$}\ \ .
          \end{cases}       
\end{equation*}
Then we can choose the ideals in the sequence \eqref{eq:ideal filtration} which contains $\fa_m$ to be all $\cJ_k+\fa_m$ such that  $(\fa_m+ \cJ_k)/ (\fa_m+ \cJ_{k+1})\neq 0$; and the ideals contained in $\fa_m$ to be all $\cJ_k\cap \fa_m$ such that  $(\fa_m\cap \cJ_k)/ (\fa_m\cap \cJ_{k+1})\neq 0$. 
\end{proof}

We next choose some sufficiently large and divisible integer $m_0,p>0$ such that 
\begin{enumerate}
    \item the line bundles $\frac{p}{M}G$ and $(\frac{p}{M}G-E_i)|_W$ are globally generated for all $i=1,\cdots,r$;
    \item $H^i(W,\cO_W(mpG-n_1 E_1-\cdots-n_r E_r))=0$ for all $i,m\in\bN_+$ and all $(n_1,\cdots,n_r)\in\bN$ with $\sum^r_{i=1} n_i\le Mm$ (this is possible by Fujita vanishing); and
    \item $\cO_Y(mpG)\otimes \cI^{Mm}_W$ is globally generated and $H^j(Y, \cO_Y(mpG)\otimes \cI^{Mm}_W)=0$ for $m\ge m_0, j\in \mathbb N_+$ (this holds as long as $m_0\gg 0$ and $p\cdot h^*G-M\cdot E$ is ample on the blowup $h\colon Y'\to Y$ along $W$ with exceptional divisor $E$).
\end{enumerate}
Consider the filtration \eqref{eq:ideal filtration}. Let $\cI_1\supseteq \cI_2$ be two consecutive terms, then we have the exact sequence
$$0\to \cO_Y(mpG)\otimes \cI_2 \to \cO_Y(mpG)\otimes \cI_1 \to \cO_W(mpG|_W)\otimes (\cI_1/\cI_2) \to 0. $$
Since $(\cI_1/\cI_2)\cong \cO_W(-n_1E_1-\cdots -n_rE_r)$ for some $(n_1,...,n_r)\in \bN^r$ and $\sum n_i\le Mm$,  then by (2),
$$H^i(W, \cO_W(mpG|_W)\otimes (\cI_1/\cI_2))=0 \qquad \mbox{any \ } i>0,$$
thus if $H^i(Y, \cO_Y(mpG)\otimes \cI_2)=0$ for $i>0$, then 
$$H^i(Y, \cO_Y(mpG)\otimes \cI_1)=0 \qquad \mbox{any \ } i>0.$$ 
If moreover $\cO_Y(mpG)\otimes \cI_2$ is globally generated, then as $\cO_W(mpG|_W)\otimes (\cI_1/\cI_2)$ is globally generated by (1), we know that
$\cO_Y(mpG)\otimes \cI_2$ is globally generated by diagram chasing and the vanishing $H^1(Y, \cO_Y(mpG)\otimes \cI_2)=0$.

Therefore, working inductively and starting from $ \cI^{Mm}_W$ by (3), we conclude that for any ideal sheaf $\cI\subseteq \cO_Y$ that appears in the sequence \eqref{eq:ideal filtration}, the sheaf $\cO_Y(mpG)\otimes \cI$ is globally generated and $H^i(Y,\cO_Y(mpG)\otimes \cI)=0$ for any $i\in\bN_+$ and any $m\ge m_0$. In particular, $\cO_Y(mpG)\otimes \fa_m$ is globally generated for all $m\ge m_0$, which implies that $p\rho^*G-A_{X,\Delta}(F)\cdot F$ is nef. In other words, $\varepsilon_v\ge \frac{1}{p}$ and we are done since the integer $p$ does not depend on the valuation $v$.
\end{proof}

\begin{proof}[Proof of Theorem \ref{t-alphabelow}]
The result now follows from a combination of Lemmas \ref{lem:alpha bound degeneration}, \ref{lem:alpha bound plt blowup} and \ref{lem:seshadri bound}.
\end{proof}
We now have all the ingredients to prove Theorem \ref{thm:f.g. for special valuation} and Theorem \ref{thm:f.g. for delta minimizer}.

\begin{proof}[Proof of Theorem \ref{thm:f.g. for special valuation}]
By Lemma \ref{lem:complement for f.g. special valuation} we already have $(1) \Rightarrow (2)$, so it remains to prove $(2)\Rightarrow (1)$. Let $v$ be a valuation that satisfies (2). We first show that $\gr_v R$ is finitely generated. To this end, let $\Sigma\subseteq \QM(Y,E)\cap \DMR(X,\Delta+\Gamma)$ be the smallest rational PL subspace containing $v$, and let $U\subseteq \Sigma$ be a small open neighbourhood of $v$ such that the closure of $U$ is contained in the interior of the simplicial cone in $\QM(Y,E)$ that contains $v$. By Theorem \ref{thm:f.g. criterion}, it is enough to show that the set $\{(X_w,\Delta_w)\,|\,w\in U(\bQ)\}$ is bounded. By Theorem \ref{t-alphabounded}, this is true if $\alpha(X_w,\Delta_w)\ge \alpha$ for some constant $\alpha>0$ that does not depend on $w\in U(\bQ)$, which then follows from Theorem \ref{t-alphabelow}. Next, by Lemma \ref{l-finitegeneratedimpliesisom}, we have $(X_v,\Delta_v)\cong (X_w,\Delta_w)$. Since $\alpha(X_w,\Delta_w)\ge \alpha>0$ and in particular $(X_w,\Delta_w)$ is klt, we see that $(X_v,\Delta_v)$ is also klt. This finished the proof.
\end{proof}

\begin{proof}[Proof of Theorem \ref{thm:f.g. for delta minimizer}]
This follows immediately from Theorem \ref{thm:f.g. for special valuation} and Corollary \ref{cor:minimizer has special complement}.
\end{proof}

\section{Applications}

In this section we present some applications of the finite generation results from the previous section. As we mentioned, combining with earlier works, Theorem \ref{t-HRFG} solves a number of major questions on the study of K-stability of Fano varieties.

\subsection{Optimal degeneration}
\begin{thm}[Optimal Destabilization Conjecture]\label{thm:ODC n+1/n}
Let $(X,\Delta)$ be a log Fano pair of dimension $n$ such that $\delta(X,\Delta)< \frac{n+1}{n}$, then $\delta(X,\Delta)\in \bQ$ and there exists a divisorial valuation $E$ over $X$ such that
$\delta(X,\Delta)=\frac{A_{X,\Delta}(E)}{S_{X,\Delta}(E)}$. 

In particular, if $\delta(X,\Delta)\le 1$, then there exists a non-trivial special test configuration $(\cX,\Delta_{\cX})$ with a central fiber $(X_0,\Delta_0)$ such that $\delta(X,\Delta)=\delta(X_0,\Delta_0)$, and $\delta(X_0,\Delta_0)$ is computed by the $\mathbb{G}_m$-action induced by the test configuration structure.
\end{thm}

\begin{proof}
Let $v$ be a valuation on $X$ that computes $\delta(X,\Delta)$. By Corollary \ref{cor:minimizer has special complement}, there exists some complement $\Gamma$ of $(X,\Delta)$ such that $v\in \DMR(X,\Delta+\Gamma)$. Let $\Sigma\subseteq \DMR(X,\Delta+\Gamma)$ be the smallest rational PL subspace containing $v$. By Theorems \ref{thm:f.g. for delta minimizer} and  \ref{thm:f.g. criterion}, the $S$-invariant function $w\mapsto S_{X,\Delta}(w)$ on $\Sigma$ is linear in a neighbourhood of $v$. As $v$ computes $\delta(X,\Delta)$, we have
\[
A_{X,\Delta}(v) = \delta(X,\Delta)S_{X,\Delta}(v).
\]
Since the log discrepancy function $w\mapsto A_{X,\Delta}(w)$ is linear in a neighbourhood of $v\in \Sigma$ and by the definition of stability thresholds we have
\[
A_{X,\Delta}(w)\ge \delta(X,\Delta)S_{X,\Delta}(w)
\]
for all $w\in \Sigma$, we see that
\[
A_{X,\Delta}(w) = \delta(X,\Delta)S_{X,\Delta}(w)
\]
in a neighbourhood $U\subseteq \Sigma$ of $v$. In particular, any divisorial valuation $w\in U(\bQ)$ also computes $\delta(X,\Delta)$. Since $w$ is a divisorial lc place of a complement, it induces a weakly special test configuration of $(X,\Delta)$ by \cite{BLX-openness}*{Theorem A.2}. By \cite{Li-criterion}*{Proof of Theorem 3.7} or \cite{Fujita-criterion}*{Theorem 5.2} we know that $\beta_{X,\Delta}(w)=A_{X,\Delta}(w)-S_{X,\Delta}(w)$ is rational. Since $A_{X,\Delta}(w)$ is clearly rational, we see that $\delta(X,\Delta)$ is also rational.

For the last part, it follows from \cite{BLZ-opt-destabilization}*{Theorem 1.1} as the conjectural assumption there is verified by the first part.
\end{proof}

\begin{thm}[Yau-Tian-Donaldson conjecture]\label{thm:YTD}
A log Fano pair $(X,\Delta)$ is uniformly K-stable if and only if it is K-stable; and it is reduced uniformly K-stable if and only if it K-polystable. In particular, $(X,\Delta)$ admits a weak KE metric if and only if it is K-polystable. 
\end{thm}

\begin{proof}
Suppose first that $\delta(X,\Delta)\le 1$. Then by Theorem \ref{thm:ODC n+1/n}, the stability threshold is computed by some divisor $E$ over $X$. By \cite{BX-uniqueness}*{Theorem 4.1}, this implies that $(X,\Delta)$ is not K-stable. In other words, if $(X,\Delta)$ is K-stable, then $\delta(X,\Delta)>1$, i.e. $(X,\Delta)$ is uniformly K-stable.

Suppose next that $(X,\Delta)$ is K-polystable. Let ${\mathbb{T}}\subseteq \Aut(X,\Delta)$ be a maximal torus. We show that $\delta_{\mathbb{T}}(X,\Delta)>1$. Suppose not, then by \cite{XZ-CM-positive}*{Appendix A}, we know that $\delta_{\mathbb{T}}(X,\Delta)=1$ and $\delta(X,\Delta)$ is computed by some ${\mathbb{T}}$-invariant quasi-monomial valuation $v$ that is 
not of the form $\wt_\xi$ for any $\xi \in \Hom(\bG_m,{\mathbb{T}})\otimes_\bZ \bR$. Moreover, $v$ is an lc place of a complement. Let $m\in\bN$ be sufficiently divisible and consider the ${\mathbb{T}}$-invariant linear system
\[
\cM:=\{s\in H^0(-m(K_X+\Delta))\,|\,v(s)\ge m\cdot A_{X,\Delta}(v)\}.
\]
Let $D_0\in |\cM|$ be a general member and let $D=\frac{1}{m}D_0$. Then $(X,\Delta+\frac{1}{m}\cM)$ has the same set of lc places as $(X,\Delta+D)$ and thus by construction $v\in\DMR(X,\Delta+D)$. Since ${\mathbb{T}}$ is a connected algebraic group, every lc place of the ${\mathbb{T}}$-invariant pair $(X,\Delta+\frac{1}{m}\cM)$ is automatically ${\mathbb{T}}$-invariant. In particular, $\DMR(X,\Delta+D)$ consists only of ${\mathbb{T}}$-invariant valuations. 

By the same argument as in the proof of Theorem \ref{thm:ODC n+1/n}, we see that $\delta(X,\Delta)$ is also computed by some divisorial valuations $w\in \DMR(X,\Delta+D)$ that are sufficiently close to $v$ (in particular, $w\neq \wt_\xi$ as well). Since $w$ is ${\mathbb{T}}$-invariant, by \cite{BX-uniqueness}*{Theorem 4.1}, $w$ induces a ${\mathbb{T}}$-equivariant special test configuration $(\cX,\cD)$ of $(X,\Delta)$ with $\Fut(\cX,\cD)=0$. Since ${\mathbb{T}}\subseteq\Aut(X,\Delta)$ is a maximal torus and $w\neq \wt_\xi$ for any $\xi\in \Hom(\bG_m,{\mathbb{T}})\otimes_\bZ \bR$, we deduce that $(\cX,\cD)$ is not a product test configuration. But this contradicts the K-polystability assumption of $(X,\Delta)$. Therefore, we must have $\delta_{\mathbb{T}}(X,\Delta)>1$ and $(X,\Delta)$ is reduced uniformly K-stable. 

When $\bk=\mathbb C$,the existence of KE metric now follows from this equivalence and \cite{Li19}*{Theorem 1.2} (see also \cites{BBJ-variational, Li-Tian-Wang}).
\end{proof}

\begin{thm}[K-moduli conjecture]\label{thm:K-moduli conj}
The K-moduli space $M^{\rm Kps}_{n,V,C}$ is proper, and the CM line bundle on $M^{\rm Kps}_{n,V,C}$ is ample. 
\end{thm}

\begin{proof}
The properness follows from Theorem \ref{thm:ODC n+1/n} and \cite{BHLLX-theta}*{Corollary 1.4}. The ampleness of the CM line bundle follows from Theorem \ref{thm:YTD} and \cite{XZ-CM-positive}*{Theorem 1.1}.
\end{proof}

\subsection{Twisted stability by adding a general boundary}
Our final application is the proof of a modified version of a conjecture of Donaldson (see \cite{Don-cone}*{Conjecture 1}, \cite{Sze-conic} and \cite{BL18}*{Section 7}). Using the Optimal Degeneration Theorem \ref{t-ODC}, we can reduce the calculations to a maximal degeneration of the log Fano pair. Nevertheless, we still need a subtle analysis of valuations in a neighborhood of the $\delta$-minimizers. 

\begin{thm}\label{t-DonConj}
Let $(X,\Delta)$ be a log Fano pair such that $\delta:=\delta(X,\Delta)<1$. Then $(X,\Delta+(1-\delta)D)$ is K-semistable for any sufficiently divisible integer $m\in\bN$ and any general $D\in \frac{1}{m}|-m(K_X+\Delta)|$. In particular, $(X,\Delta+(1-\delta')D)$ is uniformly K-stable for any $0\le \delta'<\delta$.
\end{thm}

For the proof we first need a few auxiliary lemmas. 



\begin{lem} \label{lem:interpolation}
Let $(X,\Delta)$ be a log Fano pair and let $D\sim_\bQ -(K_X+\Delta)$ be an effective $\bQ$-divisor such that $(X,\Delta+D)$ is klt. Assume that $(X,\Delta+tD)$ is K-semistable for some $t\in(0,1)$. Then $(X,\Delta+sD)$ is uniformly K-stable for all $s\in (t,1)$.
\end{lem}

\begin{proof}This simple interpolation result is well known. We include it here for the sake of
completeness.

By assumption, for any valuation $v$ on $X$ with $A_{X,\Delta}(v)<\infty$ we have $$A_{X,\Delta+tD}(v)=A_{X,\Delta}(v)-t\cdot v(D) > (1-t) v(D)$$ and $A_{X,\Delta+tD}(v)\ge S_{X,\Delta+tD}(v)=(1-t)S_{X,\Delta}(v)$. Thus for any $s\in (t,1)$, we have
\begin{align*}
    A_{X,\Delta+sD}(v)=A_{X,\Delta+tD}(v)-(s-t)v(D) & >\frac{1-s}{1-t}A_{X,\Delta+tD}(v) \\
    & \ge (1-s)S_{X,\Delta}(v)=S_{X,\Delta+sD}(v)
\end{align*}
for any valuation $v\in \Val_{X}^\circ$. Hence $(X,\Delta+sD)$ is uniformly K-stable.
\end{proof}

\begin{lem} \label{lem:Bertini in family}
Let $X$ be a normal projective variety
with an ample line bundle $L$. Let $Z\subseteq X\times U$ be a bounded flat family of positive dimensional normal subvarieties of $X$ over a normal variety $U$. Let $\Gamma$ be an effective $\bQ$-divisor on $Z$ that does not contain any $Z_u$ $($so the restriction $\Gamma_u:=\Gamma|_{Z_u}$ is well-defined as a $\bQ$-divisor on $Z_u )$. Assume that $K_{Z/U}+\Gamma$ is $\bQ$-Cartier and $(Z_u,\Gamma_u)$ is klt for all $u\in U$. 

Then there exists some constant $a>0$ such that for all sufficiently large $m\in\bN$, a general member $D\in |mL|$ does not contain any $Z_u$ and $(Z_u,\Gamma_u+aD_u)$ is lc for all $u\in U$.
\end{lem}

\begin{proof}
For any closed point $y\in Y$ of a variety $Y$ and any effective Cartier divisor $G=(g=0)$ on $Y$, we define the order of vanishing of $G$ at $y$ as
\[
\ord_y(G)=\max\{j\in\bN\,\vert\,g\in \fm_y^j\}.
\]
By the family version of the Izumi type inequality (see e.g. \cite{BL-vol-lsc}*{Theorem 20}), there exists some constant $K_0>0$ depending only on the family $(Z,\Gamma)\to U$ such that 
\[
v(D_u)\le K_0\cdot A_{Z_u,\Gamma_u}(v)\cdot \ord_x(D_u)
\]
for any $u\in U$, any effective Cartier divisor $D_u$ on $Z_u$, and any $v\in \Val_{Z_u}^\circ$ such that $x=C_{Z_u}(v)$ is a closed point. In particular this implies 
\begin{equation}\label{eq:lct>=1/ord}
    \lct_x(Z_u,\Gamma_u;D_u)\ge \frac{1}{K_0\cdot \ord_x(D_u)}
\end{equation}
for all $x\in Z_u$ and all effective Cartier divisors $D_u$ on $Z_u$.

Let $m\in\bN$ be large enough so that the restrictions 
\[
\varphi_{u,x}\colon H^0(X,\cO_X(mL))\to H^0(Z_u,\cO_{Z_u}(mL))\to H^0(Z_u,\cO_{Z_u}(mL)\otimes (\cO_{Z_u}/\fm_x^{\dim Z+1}))
\]
are surjective for any closed point $u\in U$ and $x\in Z_u$ (this is possible since $L$ is ample).  Since $\dim Z_u\ge 1$, we have $h^0(\cO_{Z_u}(mL)\otimes (\cO_{Z_u}/\fm_x^{\dim Z+1}))=h^0(\cO_{Z_u}/\fm_x^{\dim Z+1})>\dim Z$. A simple dimension count using the incidence variety
\[
W=\{(x,f)\in  Z\times H^0(X,mL)\,|\,x\in Z_u,\,\varphi_{u,x}(f)=0\}\subseteq Z\times H^0(X,mL)
\]
shows that the second projection $W\to H^0(X,mL)$ is not surjective. 

Hence if $D\in |mL|$ is a general member, then $\ord_x(D_u)\le \dim Z$ for all $u\in U$ and $x\in Z_u$. By \eqref{eq:lct>=1/ord}, this implies $\lct(Z_u,\Gamma_u;D_u)\ge \frac{1}{K_0\cdot \dim Z}$. Thus if we take $a=\frac{1}{K_0\cdot \dim Z}$, then $(Z_u,\Gamma_u+aD_u)$ is lc for all $u\in U$.
\end{proof}

Consider next the following setup. Let $(X,\Delta)$ be a log Fano pair.  Let $\bT<\Aut(X,\Delta)$ be a maximal torus, let $N=\Hom(\bG_m,\bT)$ be the co-weight lattice and let $M=\Hom(\bT,\bG_m)$ be the weight lattice. For sufficiently divisible $r$, we have a weight decomposition 
$$R:=\bigoplus_{m\in \bN}H^0(X,-mr(K_X+\Delta))=\bigoplus_{(m,\alpha)\in \bN\times M} R_{m,\alpha}.$$ 
For each $\xi\in N_\bR$, we set $\lambda_\xi:=\inf_{(m,\alpha)}\{\frac{\langle \xi,\alpha\rangle}{mr}\,\vert\,R_{m,\alpha}\neq 0\}$. Since $R$ is finitely generated, the function $\xi\mapsto \lambda_\xi$ is piecewise linear with rational coefficients. Moreover, for sufficiently divisible $m$, $\lambda_\xi=\inf_{\alpha}\{\frac{\langle \xi,\alpha\rangle}{mr}\,\vert\,R_{m,\alpha}\neq 0\}$.
Then the valuation $\wt_\xi$ is given by $\wt_\xi (s) = \langle \xi,\alpha\rangle - \lambda_\xi mr$ for all $0\neq s\in R_{m,\alpha}$. In fact, let $s^*\in R_{m,\alpha}$ such that $\wt_\xi(s^*)=\lambda_\xi mr$. Then the trivialization of $-rm(K_X+\Delta)$ around $C_{X}(\wt_{\xi})$ is given by $s\to \frac{s}{s^*}$, thus  
$$\wt_{\xi}(s)=\wt_{\xi}(\frac{s}{s^*}\cdot s^*)= \wt_{\xi}(\frac{s}{s^*})=\langle \xi,\alpha\rangle - \lambda_\xi mr.$$

\begin{lem} \label{lem:A and S linear on N}
Let $V\subseteq N_\bR$ be a convex subset where $\xi\mapsto \lambda_\xi$ is linear. Then the functions
\[
\xi\mapsto A_{X,\Delta}(\wt_\xi) \quad\text{and}\quad \xi\mapsto S_{X,\Delta}(\wt_\xi)
\]
are both linear on $V$.
\end{lem}

\begin{proof}
We choose a basis $s_1,\cdots, s_{N_m}$ of $R_m$  such that each $s_i\in R_{m,\alpha_i}$ for some $\alpha_i\in M$. In other words, $\{s_1,\cdots, s_{N_m}\}$ is a disjoint union of bases of $R_{m,\alpha}$ over all $\alpha\in M$. From the definition of $\wt_\xi$ we know that $\cF_{\wt_\xi}^\lambda R_m$ is a direct sum of some of the $R_{m,\alpha}$'s for every $\lambda\in \bR_{\geq 0}$. Thus the basis $s_1,\cdots, s_{N_m}$ is compatible with $\wt_\xi$ for every $\xi\in N_{\bR}$. Hence we have 
\[
S_m(\wt_\xi) = \frac{1}{mrN_m}\sum_{i=1}^{N_m} \wt_\xi(s_i) = \frac{1}{mrN_m}\sum_{i=1}^{N_m} (\langle \xi,\alpha_i \rangle - \lambda_\xi mr).
\]
Since $\xi\mapsto \lambda_\xi$ is linear on $V$, the above equation implies that $\xi\mapsto S_m(\wt_\xi)$ is also linear on $V$. Therefore, $\xi\mapsto S_{(X,\Delta)}(\wt_\xi)$ is linear on $V$ as $S_{(X,\Delta)}(\wt_\xi)= \lim_{m\to\infty} S_m(\wt_\xi)$. Since $A_{X,\Delta}(\wt_\xi)-S_{X,\Delta}(\wt_\xi)=\Fut(\xi)$ is always linear on $N_\bR$, the lemma follows.
\end{proof}

Since $R=\oplus_{(m,\alpha)\in \bN\times M}R_{m,\alpha}$ is finitely generated, the sub-semigroup $\Lambda:=\{(m,\alpha)\in \bN\times M\mid R_{m,\alpha}\neq 0\}$ of $\bN\times M$ is also finitely generated. As a result,  $\lambda_\xi$ is equal to the minimum of finitely many rational linear functions in $\xi$ of the form $\frac{\langle \xi,\alpha\rangle}{mr}$. Thus $\lambda_\xi$ is a rational piecewise linear function in $\xi$ on $N_{\bR}$.
We may decompose $N_\bR$ into a fan consisting of finitely many rational simplicial cones such that $\xi\mapsto \lambda_\xi$ is linear on each cone.

For each cone $\sigma$ of the fan (other than the origin), we choose some $\xi_\sigma\in N_{\bR}$ in its interior and let $Z_\sigma:=C_X(\wt_{\xi_\sigma})$. 

\textbf{Claim}. $Z_\sigma$ does not depend on the choice of $\xi_\sigma$, and $Z_\sigma\subseteq Z_\tau$ if $\sigma\supseteq \tau$. 

Next, we prove the claim. Let $\Lambda_\xi:=\{(m,\alpha)\in \Lambda\mid \langle \xi, \alpha\rangle>\lambda_\xi mr\}$. We first show that 
\begin{equation}\label{eq:center-wt}
C_X(\wt_\xi) = \bigcap_{(m,\alpha)\in \Lambda_\xi} \mathrm{Bs}(R_{m,\alpha}),
\end{equation}
where $\mathrm{Bs}(\cdot)$ denotes the base locus of a linear system. By the definition of $\wt_\xi$, we know that $\wt_\xi(s)>0$ for any $s\in R_{m,\alpha}$ with $(m,\alpha)\in \Lambda_\xi$. Thus the  ``$\subseteq$'' direction of \eqref{eq:center-wt} is clear. For the ``$\supseteq$'' direction of \eqref{eq:center-wt}, let $\cI_{\xi}\subset\cO_X$ be the ideal sheaf of $C_X(\wt_\xi)$ with reduced scheme structure. Since $-r(K_X+\Delta)$ is an ample line bundle, we fix a sufficiently large $m_1\in \bN$ such that $\cI_\xi\otimes\cO_X(-m_1 r(K_X+\Delta))$ is globally generated.  Clearly,
\[
H^0(X,\cI_\xi\otimes\cO_X(-m_1 r(K_X+\Delta))) = \{s\in R_{m_1}\mid \wt_{\xi}(s)>0\} = \bigoplus_{(m_1,\alpha)\in \Lambda_{\xi}} R_{m_1,\alpha}.
\]
Thus the global generation of $\cI_\xi\otimes\cO_X(-m_1 r(K_X+\Delta))$ implies 
\[
C_X(\wt_\xi) = \mathrm{Bs}(\bigoplus_{(m_1,\alpha)\in \Lambda_{\xi}} R_{m_1,\alpha}) = \bigcap_{(m_1,\alpha)\in \Lambda_\xi} \mathrm{Bs}(R_{m_1,\alpha})\supseteq \bigcap_{(m,\alpha)\in \Lambda_\xi} \mathrm{Bs}(R_{m,\alpha}).
\]
This finishes the proof of \eqref{eq:center-wt}. 

Back to the proof of claim. Recall that $\xi_\sigma$ belongs to the interior of $\sigma$. To prove the claim, it suffices to show that $C_X(\wt_{\xi_\sigma})\subseteq C_X(\wt_{\xi'})$ whenever $\xi'\in \sigma$. By \eqref{eq:center-wt}, this reduces to showing  $\Lambda_{\xi_\sigma}\supseteq \Lambda_{\xi'}$. Let $(m,\alpha)\in \Lambda\setminus \Lambda_{\xi_\sigma}$, i.e. $\langle \xi_\sigma,\alpha\rangle = \lambda_{\xi_\sigma} mr$. 
Since the function $\xi\mapsto \langle \xi,\alpha\rangle - \lambda_\xi mr$ is linear and nonnegative on $\sigma$, its vanishing at an interior point $\xi_\sigma$ of $\sigma$ implies that  it vanishes on the whole of $\sigma$. Thus $\langle \xi',\alpha\rangle = \lambda_{\xi'} mr$ which implies that $(m,\alpha)\in \Lambda\setminus\Lambda_{\xi'}$. Therefore, we have shown $\Lambda \setminus \Lambda_{\xi_\sigma}\subseteq \Lambda\setminus\Lambda_{\xi'}$, i.e. $\Lambda_{\xi_\sigma}\supseteq\Lambda_{\xi'}$. This finishes the proof of the claim.

As a consequence of the claim, when $\sigma$ varies in the fan, $Z_\sigma$ enumerates the center of $\wt_\xi$ for all $\xi\in N_\bR$. 

Since $\wt_{\xi_\sigma}$ induces a product test configuration of $(X,\Delta)$ which is special, by Theorem \ref{thm:zhuang-special} there exists some $\gamma_\sigma\in(0,1)$ and some $G_\sigma\sim_\bQ -\gamma_\sigma(K_X+\Delta)$ such that $(X,\Delta+G_\sigma)$ is lc and $\wt_{\xi_\sigma}$ is its unique lc place. In particular, $Z_\sigma$ is the minimal lc center of $(X,\Delta+G_\sigma)$. Choose some $0<\varepsilon_\sigma\ll 1$ and some general $G'_\sigma\in|-K_X-\Delta|_\bQ$. By Kawamata subadjunction \cite{Kawamata-subadj}, we may write
\[
(K_X+\Delta+G_\sigma+\varepsilon_\sigma G'_\sigma)|_{Z_\sigma}\sim_\bQ K_{Z_\sigma}+\Gamma_\sigma
\]
for some divisor $\Gamma_\sigma\ge 0$ on $Z_\sigma$ such that $(Z_\sigma,\Gamma_\sigma)$ is klt. Moreover, if $D\ge 0$ is a $\bQ$-Cartier divisor on $X$ whose support does not contain $Z_\sigma$, then $(X,\Delta+G_\sigma+D)$ is lc in a neighbourhood of $Z_\sigma$ if $(Z_\sigma,\Gamma_\sigma+D|_{Z_\sigma})$ is lc (for a more precise version, see \cite{HMX-ACC}*{Proof of Theorem 4.2}). For any $g\in \Aut(X,\Delta)$, let $(Z_{\sigma,g},\Gamma_{\sigma,g})=(g\cdot Z_\sigma,g\cdot \Gamma_\sigma)$. We will eventually apply Lemma \ref{lem:Bertini in family} to the family $(Z_{\sigma,g},\Gamma_{\sigma,g})$. 

For now we state a technical result that is needed in our proof of Theorem \ref{t-DonConj}.

\begin{lem} \label{lem:v(D) estimate}
Assume that $\dim Z_\sigma\ge 1$. Let $\tau$ be a cone of the fan on $N_\bR$ such that $\sigma\subseteq \tau$. Let $\xi_0\in \sigma$, $\xi_1\in \tau$ and let $\xi_t=(1-t)\xi_0+t\xi_1$ for $t\in[0,1]$. Then for any $0\neq s\in R_m$ such that $Z_\sigma$ is not contained in the support of $D=(s=0)$, we have
\[
\wt_{\xi_t}(s)\le \frac{t\cdot A_{X,\Delta}(\wt_{\xi_1})}{\lct(Z_\sigma,\Gamma_\sigma;D|_{Z_\sigma})}.
\]
\end{lem}

\begin{proof}
For ease of notation, let $v_t=\wt_{\xi_t}$, $\xi=\xi_\sigma\in \sigma$, $G=G_\sigma$, $Z=Z_\sigma$, and $\Gamma=\Gamma_\sigma$. Using the weight decomposition, we may write $s=\sum_{\alpha\in M} s_\alpha$ and by definition 
\begin{equation} \label{eq-v}
v_t(s)=\min\{v_t(s_\alpha)\,\vert\,s_\alpha\neq 0\}.
\end{equation}
Let $M_s:=\{\alpha\in M \mid s_\alpha\neq 0 \textrm{ and }\langle \xi,\alpha\rangle = \lambda_{\xi} mr\}$.
Let $s_1:=\sum_{\alpha\in M_{s}} s_\alpha$. Note that by the definition and the linearity of $\lambda_\xi$ on $\sigma$, for each $\alpha\in M_s$ we necessarily have
\begin{equation} \label{eq:v_0=0}
    \langle \xi_0,\alpha\rangle = \lambda_{\xi_0} mr
\end{equation}
as well. By assumption, we have $M_s\neq\emptyset$ and $s_1\neq 0$ as otherwise $Z\subseteq \Supp(D)$. Since $\lambda_\xi$ is linear in $\tau$ and $\xi_0, \xi_1 \in \tau$, we know that $t\mapsto \lambda_{\xi_t}$ is linear for $t\in [0,1]$. Thus for each $\alpha\in M_s$ and $t\in [0,1]$ we have
\begin{equation} \label{eq:v_t=tv_1}
v_t(s_\alpha)=\langle\xi_t, \alpha\rangle-\lambda_{\xi_t}mr = (1-t) (\langle\xi_0, \alpha\rangle-\lambda_{\xi_0}mr ) + t (\langle\xi_1, \alpha\rangle-\lambda_{\xi_1}mr ) = tv_1(s_\alpha),
\end{equation}
where the last equality follows from \eqref{eq:v_0=0}.
Let $D_1=(s_1=0)$.
Then we have $D_1|_Z=D|_Z$ and 
\begin{equation}\label{eq:v_t(s)less}
v_t(s)\le v_t(s_1)=t\cdot v_1(s_1),
\end{equation}
where the first inequality follows from \eqref{eq-v}, and the second equality uses \eqref{eq:v_t=tv_1} and the fact that 
\[
v_t(s_1)= \min_{\alpha\in M_s} \{ v_t(s_\alpha)\} \quad \textrm{and}\quad v_1(s_1)= \min_{\alpha\in M_s} \{ v_1(s_\alpha)\}.
\]

Thus to prove the lemma, by \eqref{eq:v_t(s)less} it suffices to show that
\[
v_1(s_1)\le \frac{A_{X,\Delta}(v_1)}{\lct(Z,\Gamma;D_1|_Z)}, \quad \text{or equivalently}\quad \lct(Z,\Gamma;D_1|_Z)\le \frac{A_{X,\Delta}(v_1)}{v_1(D_1)}.
\]
But as $C_X(v_1)\cap Z\supseteq Z_\tau$ is non-empty, by subadjunction we have
$$\lct(Z,\Gamma;D_1|_Z)\le \lct_Z(X,\Delta+G;D_1)\le \lct_Z(X,\Delta;D_1)\le \frac{A_{X,\Delta}(v_1)}{v_1(D_1)},$$ where the two log canonical thresholds in the middle are taken in a neighbourhood of $Z$, and we are done.
\end{proof}

We are now ready to present the proof of Theorem \ref{t-DonConj}.

\begin{proof}[Proof of Theorem \ref{t-DonConj}]
First note that if $(X,\Delta)$ has a special degeneration to some log Fano pairs $(X_0,\Delta_0)$ with $\delta(X_0,\Delta_0) = \delta$ and the theorem holds for $(X_0,\Delta_0)$, then it also holds for $(X,\Delta)$ by the openness of K-semistability \cite{BLX-openness, Xu-quasimonomial}. The following claim shows that the process of special degenerations preserving the stability threshold will stabilize after finitely many steps. 

\textbf{Claim}. Any sequence of special degenerations
\[
(X,\Delta)=:(X^{(0)},\Delta^{(0)})\rightsquigarrow (X^{(1)},\Delta^{(1)})\rightsquigarrow (X^{(2)},\Delta^{(2)})\rightsquigarrow  \cdots \rightsquigarrow (X^{(i)},\Delta^{(i)}) \rightsquigarrow\cdots
\]
satisfying that $\delta(X^{(i)},\Delta^{(i)})=\delta$ and $(X^{(i)},\Delta^{(i)})\not\cong (X^{(i+1)},\Delta^{(i+1)})$ for every $i\geq 0$ must terminate after finitely many steps.

Assuming the claim, there exists a finite sequence of special degenerations $(X,\Delta) \rightsquigarrow \cdots \rightsquigarrow (X^{(k)},\Delta^{(k)})$ preserving stability thresholds such that any special degeneration $(X^{(k)},\Delta^{(k)})\rightsquigarrow (X^{(k+1)},\Delta^{(k+1)})$ preserving the stability threshold is of product type, i.e. $(X^{(k)},\Delta^{(k)})\cong (X^{(k+1)},\Delta^{(k+1)})$. Thus from the above argument using openness of K-semistability, we may replace $(X,\Delta)$ by $(X^{(k)},\Delta^{(k)})$. 

Next, we prove the claim. 
Since the set of log Fano pairs with fixed volume, finite rational coefficient set, and  $\delta$-invariant is bounded by Theorem \ref{t-alphabounded}, we can embed every $(X^{(i)},\Delta^{(i)})$ into a common projective space $\bP^N$ using a common multiple of its anti-log-canoncial divisor, such that each special degeneration $(X^{(i)},\Delta^{(i)})\rightsquigarrow (X^{(i+1)},\Delta^{(i+1)})$ is induced by a one parameter subgroup of $G:=\mathrm{PGL}_{N+1}$. Let $\bfP$ denote a locally closed subscheme of the relative Hilbert-Chow scheme of $\bP^N$ such that $\bfP$ is of finite type, and each Hilbert-Chow point $z_i:=[(X^{(i)}, \Delta^{(i)})]$ belongs to $\bfP$, see e.g. \cite[Section 6]{BLZ-opt-destabilization} or \cite[Section 4.1]{BHLLX-theta} for details. Then by construction, we know that $z_{i+1}\in \overline{G\cdot z_i}\setminus G\cdot z_i$. This implies that $\overline{G\cdot z_{i+1}}\subsetneq \overline{G\cdot z_i}$ as closed subsets of $\bfP$. Since $\bfP$ is of finite type by boundedness, it is a Noetherian topological space. As a result, the sequence 
$\overline{G\cdot z_0}\supsetneq \overline{G\cdot z_1}\supsetneq \cdots$ must terminate after finitely many steps. Thus the proof of the claim is finished. 

Thanks to the claim, in the rest of the proof we may assume that any special degeneration $(X_0,\Delta_0)$ of $(X,\Delta)$ with $\delta(X_0,\Delta_0) = \delta$ satisfies $(X_0,\Delta_0)\cong (X,\Delta)$, i.e. the degeneration is induced by a one parameter subgroup of $\Aut(X,\Delta)$.
Let $m\in\bN$ be sufficiently large and divisible and let $D_m\in \frac{1}{m}|-m(K_X+\Delta)|$ be general. Since $(X,\Delta+mD_m)$ is lc by Bertini's theorem, we have
\[
A_{X,\Delta}(v)\ge A_{X,\Delta+(1-\delta)D_m}(v)=A_{X,\Delta}(v)-(1-\delta)v(D_m)\ge A_{X,\Delta}(v)-\frac{1-\delta}{m}A_{X,\Delta}(v)
\]
for any valuation $v$ on $X$. Since $S_{X,\Delta+(1-\delta)D_m}(v)=\delta\cdot S_{X,\Delta}(v)$, this implies that 
\[
1-\frac{1-\delta}{m}=\left(1-\frac{1-\delta}{m}\right)\delta^{-1}\cdot \delta(X,\Delta)\le \delta(X,\Delta+(1-\delta)D_m)\le \delta^{-1}\cdot \delta(X,\Delta)= 1.
\]
By Theorem \ref{thm:ODC n+1/n}, there is a special degeneration, 
\[
(X,\Delta+(1-\delta)D_m)\rightsquigarrow (X_m,\Delta_m+(1-\delta)G_m)
\]
induced by some special divisorial valuation $v_m$, such that 
\begin{equation} \label{eq:delta equal}
    \delta(X_m,\Delta_m+(1-\delta)G_m)=\delta(X,\Delta+(1-\delta)D_m)\ge 1-\frac{1-\delta}{m}.
\end{equation}
This implies
\begin{align*}
    A_{X_m,\Delta_m}(v)\ge A_{X_m,\Delta_m+(1-\delta)G_m} (v) & \ge \left(1-\frac{1-\delta}{m}\right)S_{X_m,\Delta_m+(1-\delta)G_m}(v) \\
    & =\left(1-\frac{1-\delta}{m}\right)\delta\cdot S_{X_m,\Delta_m}(v)
\end{align*}
for all valuation $v$ over $X_m$ and hence
\begin{equation} \label{eq:delta lower bound}
    \delta(X_m,\Delta_m)\ge \left(1-\frac{1-\delta}{m}\right)\delta
\end{equation}
is bounded from below by some constants that only depend on $\delta$. By Theorem \ref{t-alphabounded}, we see that $(X_m,\Delta_m)$ belongs to a bounded family.  

Since the $\delta$-invariant function (more precisely $\min\{\delta(X,\Delta),1\}$) is lower semi-continuous and constructible in a bounded family \cite{BLX-openness}*{Theorem 1.1}, we have $\delta(X_m,\Delta_m) = \delta$ when $m$ is sufficiently large by \eqref{eq:delta lower bound}. Thus by our assumption on $(X,\Delta)$ at the beginning of the proof, we see that $v_m$ is induced by a one-parameter subgroup of $\Aut(X,\Delta)$. It follows that in order to prove the K-semistability of $(X,\Delta+(1-\delta)D_m)$ for $m\gg 0$, it is enough to show $A_{X,\Delta+(1-\delta)D_m}(v)\ge S_{X,\Delta+(1-\delta)D_m}(v)$, or equivalently,
\[
A_{X,\Delta}(v) \ge (1-\delta)v(D_m) + \delta\cdot S_{X,\Delta}(v)
\]
for all $v\in \Val_X$ that are induced by one-parameter subgroups of $\Aut(X,\Delta)$. 

Fix a maximal torus $\bT$ of $\Aut(X,\Delta)$. Since all maximal tori are conjugate and the functions $A_{X,\Delta}, S_{X,\Delta}$ are $\Aut(X,\Delta)$-invariant, it suffices to show that
\begin{equation} \label{eq:twisted K-ss}
    A_{X,\Delta}(v) \ge (1-\delta)v(g\cdot D_m) + \delta\cdot S_{X,\Delta}(v)
\end{equation}
for all $v\in\Val_X$ of the form $\wt_\xi$ ($\xi\in N_\bR$) and all $g\in \Aut(X,\Delta)$.

\medskip

By \cite{Z-equivariant}*{Theorem 1.5}, there exists an $\Aut(X,\Delta)$-invariant closed subvariety $W$ of $X$ such that $W$ is contained in $C_X(v)$ for any valuation $v$ computing $\delta(X,\Delta)$. Consider the simplicial fan structure on $N_\bR$ induced by the piecewise linear function $\xi\mapsto \lambda_\xi$ as before. By Bertini's theorem we may assume that $(X,\Delta+mD_m)$ is lc and $W\not\subseteq \Supp(D_m)$. This also implies that $(X,\Delta+m(g\cdot D_m))$ is lc and $W\not\subseteq \Supp(g\cdot D_m)$. By applying Lemma \ref{lem:Bertini in family} to the effective Cartier divisors $m D_m$ and the families of pairs $(Z_{\sigma,g},\Gamma_{\sigma,g})$ parametrized by $U=\Aut(X,\Delta)$ for all $\sigma$ such that $\dim Z_{\sigma}\ge 1$, we also know that there exists some constant $a>0$ independent of $m$ such that $\lct(Z_{\sigma,g},\Gamma_{\sigma,g};D_m|_{Z_{\sigma,g}})\ge ma$, or equivalently,
\begin{equation} \label{eq:lct bound uniform in m}
    \lct(Z_\sigma,\Gamma_\sigma;g\cdot D_m|_{Z_\sigma})\ge ma
\end{equation}
for all $\sigma$ satisfying $\dim Z_\sigma\ge 1$ and all $g\in \Aut(X,\Delta)$. 

Let $\tau_1,\cdots,\tau_k$ be the cones of maximal dimension in the fan of $N_{\bR}$. We shall analyze the behavior of $A_{X,\Delta}$, $S_{X,\Delta}$, and $v(g\cdot D_m)$ on each cone $\tau_i$. For each $i=1,\cdots,k$, let 
\[
\sigma_i = \{\xi\in \tau_i\,\vert\,A_{X,\Delta}(\wt_\xi)=\delta\cdot S_{X,\Delta}(\wt_\xi)\}.
\]
By Lemma \ref{lem:A and S linear on N} and the fact that $A_{X,\Delta}(v)\ge \delta\cdot S_{X,\Delta}(v)$ for all $v\in\Val_X^\circ$, we see that $\sigma_i$ is a face of $\tau_i$. Let $\sigma'_i\subseteq \tau_i$ be the smallest face such that $\tau_i$ is the convex hull of $\sigma_i$ and $\sigma'_i$ (such $\sigma'_i$ exists since $\tau_i$ is simplicial). In particular, we have $\sigma_i\cap \sigma'_i=\{0\}$ and therefore there exists some constant $\varepsilon_0\in(0,1)$ such that
\begin{equation} \label{eq:A>delta*S}
    A_{X,\Delta}(\wt_\xi)\ge \frac{\delta}{1-\varepsilon_0} \cdot S_{X,\Delta}(\wt_\xi)
\end{equation}
for all $i=1,\cdots,k$ and all $\xi\in\sigma'_i$. 

We now proceed to prove that \eqref{eq:twisted K-ss} holds for all $m\ge \max\{\frac{1-\delta}{\varepsilon_0},\frac{1-\delta}{a\varepsilon_0}\}$. Let $\xi\in N_\bR$ and let $v=\wt_\xi$. Let $\tau_i$ be a cone that contains $\xi$. There are three cases to consider.

\smallskip
\emph{Case 1}: $\sigma_i=\{0\}$. Then $\sigma'_i=\tau_i$. Since $(X,\Delta+m(g\cdot D_m))$ is lc, combined with \eqref{eq:A>delta*S} we have $$A_{X,\Delta}(v)-\delta\cdot S_{X,\Delta}(v)\ge \varepsilon_0\cdot A_{X,\Delta}(v)\ge m\varepsilon_0\cdot v(g\cdot D_m)\ge (1-\delta) v(g\cdot D_m)$$ for any $g\in \Aut(X,\Delta)$. Thus \eqref{eq:twisted K-ss} holds in this case.

\smallskip
\emph{Case 2}: $\sigma_i\neq\{0\}$ and $Z_{\sigma_i}$ is a point. Then we necessarily have $Z_{\sigma_i}=W$ and since $Z_{\tau_i}\subseteq Z_{\sigma_i}$ we have $Z_{\tau_i}=W$ as well. Since $W\not\subseteq \Supp(g\cdot D_m)$, we deduce that $\wt_\xi(g\cdot D_m)=0$ for any $\xi\in \tau_i$ and any $g\in \Aut(X,\Delta)$. Thus \eqref{eq:twisted K-ss} clearly holds in this case.

\smallskip
\emph{Case 3}: $\sigma_i\neq\{0\}$ and $\dim Z_{\sigma_i}\ge 1$. We can write $\xi=(1-t)\xi_0+t\cdot \xi_1$ for some $\xi_0\in\sigma_i$ and $\xi_1\in \sigma'_i$. Then by Lemma \ref{lem:v(D) estimate} and \eqref{eq:lct bound uniform in m}, we know that 
\[
v(g\cdot D_m)\le \frac{t}{ma}\cdot A_{X,\Delta}(v_1)
\]
for any $g\in\Aut(X,\Delta)$, where $v_1=\wt_{\xi_1}$. On the other hand, since $A_{X,\Delta}$ and $S_{X,\Delta}$ are linear on $\tau_i$ by Lemma \ref{lem:A and S linear on N}, we have
\[
A_{X,\Delta}(v)-\delta\cdot S_{X,\Delta}(v)\ge t\varepsilon_0\cdot A_{X,\Delta}(v_1)
\]
by \eqref{eq:A>delta*S}. Combining the two inequalities above we get
\[
A_{X,\Delta}(v)-\delta\cdot S_{X,\Delta}(v) \ge ma\varepsilon_0\cdot v(g\cdot D_m)\ge (1-\delta)v(g\cdot D_m)
\]
for all $g\in\Aut(X,\Delta)$. Thus \eqref{eq:twisted K-ss} holds in this case as well.

Thus we have proved that \eqref{eq:twisted K-ss} holds for all $v=\wt_\xi$. As explained earlier, this implies the K-semistability of $(X,\Delta+(1-\delta)D_m)$ when $m\gg 0$. The remaining part of the theorem follows from Lemma \ref{lem:interpolation}.
\end{proof}



\section{Examples}\label{s-examples}

It is still natural to ask which lc places of a given $\mathbb Q$-complement induce finitely generated associated graded rings. Unlike the divisorial case, where the finite generation is guaranteed essentially by \cite{BCHM}, in the higher rank case, the associated graded ring could generally be non-finitely generated. In fact, it was first discovered in \cite[Theorem 4.16]{AZ-index2} that on every smooth cubic surface there exist lc places of complements whose associated graded rings are not finitely generated.

In this section, we give a complete picture of the finite generation problem in an explicit example, that is, lc places of $(\bP^2, C)$ where $C$ is an irreducible nodal cubic curve. As one will see, even in this simple set-up, the locus of finitely generated lc places is fairly complicated, and there are infinitely many special degenerations of $\bP^2$. In addition, both non-finitely generated and finitely generated non-divisorial rational rank two lc places appear in the same simplex (given by a dlt modification of $(\bP^2, C)$).
This also illustrates the importance of considering special complements in previous discussions. It will be an interesting question to understand better how to locate the valuations on a dual complex with a finitely generated associated graded ring. 


\medskip

We fix the following notation. Let $o$ be the node of $C$. Choose an analytic coordinate $(z,w)$ of $\bP^2$ at $o$ such that the analytic local equation of $C$ is given by $(zw=0)$. For $t\in \bR_{>0}$, let $v_t$ be the monomial valuation of weights $(1, t)$ in the coordinate $(z,w)$. Since $C$ is normal crossing, we know that any lc place of $(\bP^2, C)$ is a multiple of $v_t$ or $\ord_C$.  Let $R:= R(\bP^2, \cO_{\bP^2}(3))$ be the anti-canonical ring of $\bP^2$. We also denote $X:=\bP^2$.

\begin{thm}\label{thm:nodalcubic}
With the above notation, the associated graded ring $\gr_{v_t} R$ is finitely generated if and only if $t\in \bQ_{>0}\cup (\frac{7-3\sqrt{5}}{2},\frac{7+3\sqrt{5}}{2})$. Moreover, $\Proj\,\gr_{v_t} R$ is a $\bQ$-Fano variety if and only if $t\in (\frac{7-3\sqrt{5}}{2},\frac{7+3\sqrt{5}}{2})$.
For a detailed description of these special degenerations of $\bP^2$, see Remark \ref{rem:P^2-deg}.
\end{thm}

The proof of Theorem \ref{thm:nodalcubic} is divided into several parts. We first recall a useful lemma.

\begin{lem}[\cite{Fujita-plt}*{Claim 4.3}]\label{lem:fujita}
Suppose $a, b$ are coprime positive integers.
Let $\mu: \tX\to X=\bP^2$ be an $(a,b)$-weighted blow up at a smooth point with exceptional divisor $E$. Let $\epsilon_X(E):= \max\{\lambda\in \bR_{\geq 0} \mid \mu^*(-K_X)-\lambda E\textrm{ is nef}\}$. Then we have 
\begin{equation}\label{eq:fujita-plt}
T_X(E)\cdot \epsilon_X(E)=9ab \quad \textrm{and}\quad  S_X(E)=\frac{T_X(E)+\epsilon_X(E)}{3}.
\end{equation}
\end{lem}

\begin{prop}\label{prop:S-non-fg}
If $t\geq \frac{7+3\sqrt{5}}{2}$ and $t\not\in \bQ$, then $\gr_{v_t} R$ is not finitely generated.
\end{prop}

\begin{proof}
By Theorem \ref{thm:f.g. criterion}, it suffices to show that the function $t\mapsto S_X(v_t)$ is not linear in any sub-interval of $[\frac{7+3\sqrt{5}}{2}, +\infty)$. 
We first compute the $S$-invariant for $v_t$ when $t=\frac{b}{a}>\frac{7+3\sqrt{5}}{2}$ where $a,b$ are coprime positive integers. Let $\mu: \tX\to X=\bP^2$ be the $(a,b)$-weighted blow up at $o$ in the analytic coordinates $(z,w)$. Let $E$ be the $\mu$-exceptional divisor. Then easy computation shows that $\tC:=\mu_*^{-1}C \sim \pi^*(-K_X)-(a+b)E$  and $(\tC^2)=9-\frac{(a+b)^2}{ab}<0$. Hence $[\tC]$ is extremal in $\overline{NE}(\tX)$ by \cite[Lemma 1.22]{KM98}. Thus we have the pseudo-effective threshold $T_X(E)=\ord_E(C)=a+b$ which implies $S_X(E)= \frac{(a+b)^2+9ab}{3(a+b)}$ by \eqref{eq:fujita-plt}. Since $v_t=\frac{1}{a}\ord_E$, for any rational $t>\frac{7+3\sqrt{5}}{2}$ we have 
\begin{equation}\label{eq:S-non-fg}
S_X(v_t)=\frac{1}{a}S_X(E)=\frac{t^2+11t+1}{3(t+1)}.
\end{equation}
Since the $S$-invariant is continuous in the dual complex \cite[Proposition 2.4]{BLX-openness}, \eqref{eq:S-non-fg} holds for any $t\in [\frac{7+3\sqrt{5}}{2}, +\infty)$. Thus $t\mapsto S_X(v_t)$ is not linear in any sub-interval.
\end{proof}


Next we turn to proving finite generation of $\gr_{v_t} R$ for $t\in [1, \frac{7+3\sqrt{5}}{2})$. The idea is to find a sequence of increasing rational numbers $1=t_0<t_1<\cdots< t_n < \cdots $ with $\lim_{n\to\infty} t_n= \frac{7+3\sqrt{5}}{2}$ such that the function $t\mapsto S_X(v_t)$ is linear in each interval $[t_n, t_{n+1}]$.

Let $(d_n)_{n\geq 0}$ be the following sequence of integers, where $d_0=1$, $d_1=1$, $d_2=2$, and $d_{n+1}=3d_n-d_{n-1}$. The sequence of $(d_n)_{n\geq 0}$
goes as $1,1,2,5,13,34, 89, \cdots$.
It is easy to see that $(d_n)$ satisfies the following properties.
\begin{itemize}
    \item Each $d_n$ is not divisible by $3$;
    \item $d_n=F_{2n-1}$ where $(F_n)$ is the Fibonacci sequence;
    \item $(1, d_n, d_{n+1})$ is a Markov triple, i.e. $1+d_n^2+d_{n+1}^2=3d_n d_{n+1}$;
    \item $d_n^2+1=d_{n-1} d_{n+1}$. 
\end{itemize}


Let $t_n:=\frac{d_{n+1}}{d_{n-1}}$ for $n\geq 1$ and $t_0:=1$. Then it is easy to see that $(t_n)$ is a strictly increasing sequence whose limit is $\frac{7+3\sqrt{5}}{2}$. 

\begin{prop}\label{prop:S-fg}
We have $S_X(v_t)= \frac{d_{n+1}}{d_n} + \frac{d_n}{d_{n+1}} t$ for $t\in [t_n, t_{n+1}]$ and $n\geq 0$.
\end{prop}

In order to prove Proposition \ref{prop:S-fg}, in the following lemma we find a sequence of very singular plane curves $D_n$ of degree $d_n$ such that they compute the $T$-invariant for $v_{t_n}$. Note that $D_3$ is precisely the singular plane quintic with an $A_{12}$-singularity (see e.g.  \cite{ADL19}*{Section 7.1}). 

\begin{lem}\label{lem:singular-plane-curve}
For each $n>0$, there exists an integral plane curve $D_n$ such that $\deg(D_n)=d_n$ and the Newton polygon of the defining function of $D_n$ in $(z,w)$ at $o$  is the line segment joining $(d_{n+1},0)$ and $(0, d_{n-1})$. 
\end{lem}

\begin{proof}
Let $\mu_n: \tX_n\to X$ be the $(d_{n-1}, d_{n+1})$-weighted blow-up in $(z,w)$ at $o$. Let $E_n$ be the $\mu_n$-exceptional divisor, so $\ord_{E_n}=d_{n-1}v_{t_n}$.
It is clear that 
$h^0(\bP^2, \cO(d_n))=\frac{d_n^2+3d_n+2}{2}$.
Using Pick's theorem, it is easy to compute that
$\colength(\fa_{d_{n-1} d_{n+1}}(\ord_{E_n}))=\frac{d_n^2+3d_n}{2}$.
Thus we have $h^0(\bP^2, \cO(d_n))>\colength(\fa_{d_{n-1} d_{n+1}}(\ord_{E_n}))$, which implies the existence of a plane curve $D_n$ of degree $d_n$ with $\ord_{E_n}(D_n)\geq d_{n-1} d_{n+1}$. 

Next, we show that the curve $D_n$ is integral. Assume to the contrary, then there exists an integral plane curve $D$ of degree $d<d_n$, such that 
$
\frac{\ord_{E_n}(D)}{d}\geq  \frac{\ord_{E_n}(D_n)}{d_n}\geq \frac{d_{n-1}d_{n+1}}{d_n}$.
In fact, we always have $\frac{\ord_{E_n}(D)}{d}>\frac{d_{n-1}d_{n+1}}{d_n}$ since $\frac{d_{n-1}d_{n+1}}{d_n}=d_n+\frac{1}{d_n}$ and $d<d_n$. 
Clearly $D\neq C$ since $
\frac{\ord_{E_n}(C)}{3}=\frac{d_{n-1}+d_{n+1}}{3}=d_n<\frac{d_{n-1}d_{n+1}}{d_n}$.
Computing local intersection numbers, 
yields $
(D\cdot C)_o\geq \ord_{E_n}(D) \left(\frac{1}{d_{n-1}}+ \frac{1}{d_{n+1}}\right)> \frac{d(d_{n-1}+d_{n+1})}{d_n}=3d$.
On the other hand, Bezout's theorem implies $(D\cdot C)_o\leq (D\cdot C)=3d$, a contradiction. Thus $D_n$ is integral. 

Finally, we show the Newton polygon statement. Suppose the Newton polygon of $D_n$ passes through $(p,0)$ and $(0,q)$. Then by computing local intersection numbers, we know that $3d_n=(D_n\cdot C)\geq (D_n\cdot C)_o = p+q$. 
On the other hand, $\ord_{E_n}(D_n)\geq d_{n-1}d_{n+1}$ implies that $p\geq d_{n+1}$ and $q\geq d_{n-1}$. Hence we must have $p= d_{n+1}$ and $q= d_{n-1}$. 
\end{proof}



\begin{proof}[Proof of Proposition \ref{prop:S-fg}]
We first treat the case when $n=0$, i.e. $t\in [1,2]$. 
Choose a suitable projective coordinates $[x_0,x_1,x_2]$ of $\bP^2$ such that $C=(x_0x_2^2=x_1^3+x_0x_1^2)$ and $o=[1,0,0]$. In the affine chart $[1,x_1,x_2]$ of $\bP^2$, let $z':=x_1-x_2$ and $w':= x_1+x_2$. Then after possibly rescaling and switching $(z,w)$, we may assume that $\ord_o(z- z')\geq 2$ and $\ord_o(w-w')\geq 2$. Let $u_t$ be the monomial valuation of weights $(1,t)$ in the coordinate $(z',w')$. We will show that $v_t=u_t$ for any $t\in [1,2]$. Since $u_t(z')=1\leq 2t=u_t(w'^2)$ and $u_t(w')=t\leq 2=u_t(z'^2)$, we know that $u_t(z)=u_t(z')=1$ and $u_t(w)=u_t(w')=t$. Hence for any $p,q\in \bZ_{\geq 0}$ we have $v_t(z^p w^q)=p+tq=u_t(z^p w^q)$. Then for any non-zero function $f\in \cO_{\bP^2,o}\setminus\{0\}$ with a Taylor expansion $f=\sum_{(p,q)\in \bZ_{\geq 0}^2} c_{p,q} z^p w^q$, we have
\[
v_t(f)=\min \{p+tq\mid c_{p,q} \neq 0\}= \min \{u_t(z^p w^q)\mid  c_{p,q} \neq 0\}\leq u_t(f).
\]
By switching $(z,w,v_t)$ and $(z',w', u_t)$, similar arguments show that $u_t(f)\leq v_t(f)$. Thus we have $v_t=u_t$ for any $t\in [1,2]$. Since $u_t$ is toric in the affine coordinate $(z',w')$, a standard toric computation shows that $S_X(v_t)=S_X(u_t)=1+t$.

From now on we may assume $n\geq 1$. 
By Lemma \ref{lem:singular-plane-curve}, we know that $(\tD_n^2)=(D_n^2)+(d_{n+1}d_{n-1})^2(E_n^2)=-1$ where $\tD_n:=(\mu_n)_*^{-1}D_n$. Thus $[\tD_n]$ is extremal in $\overline{NE}(\tX_n)$ by \cite[Lemma 1.22]{KM98}. Hence $T_X(\ord_{E_n})=\frac{3}{d_n} \ord_{E_n} (D_n)=\frac{3d_{n-1} d_{n+1}}{d_n}$. Thus Lemma \ref{lem:fujita} implies that $S_X(v_{t_n})=\frac{1}{d_{n-1}}S_X(E_n)=\frac{d_{n+1}}{d_n}+\frac{d_n}{d_{n-1}}$. So $t=t_n$ satisfies the statement of Proposition \ref{prop:S-fg}. 

Let $t_n':=\frac{d_{n+1}^2}{d_n^2}$ for $n\geq 1$, then we have $t_n<t_n'<t_{n+1}$. Let $\mu_n':\tX_n'\to X$ be the $(d_n^2, d_{n+1}^2)$-weighted blow-up in $(z,w)$ at $o$ with exceptional divisor $E_n'$. Then using Lemma \ref{lem:singular-plane-curve}, similar computation shows $(\tD_n'^2)=0$ where $\tD_n':=(\mu_n')_*^{-1}D_n$. Hence $[\tD_n']$ lies in the boundary of $\overline{NE}(\tX_n')$ by \cite[Lemma 1.22]{KM98}. Thus $T_X(\ord_{E_n'})=\frac{3}{d_n}\ord_{E_n'}(D_n)=3d_nd_{n+1}$, and Lemma \ref{lem:fujita} implies that $t=t_n'$ also satisfies the statement of Proposition \ref{prop:S-fg}. Since $S$-invariant is concave by Lemma \ref{lem:S linear imply f.g. interval}, the proof is finished.
\end{proof}

\begin{proof}[Proof of Theorem \ref{thm:nodalcubic}]
Choose the projective coordinates $[x_0,x_1,x_2]$ of $\bP^2$ as in the proof of Proposition \ref{prop:S-fg}. Then the  automorphism $\sigma$ of $(\bP^2, C)$ given by $\sigma([x_0, x_1,x_2]):=[x_0,x_1,-x_2]$ interchanges the two analytic branches of $C$ at $o=[1,0,0]$. Thus $\sigma_* v_t= t\cdot v_{t^{-1}}$ which implies that $v_{t^{-1}}$ and $v_t$ have isomorphic associated graded rings after a grading shift. So we may assume $t\in [1, +\infty)$ from now on. The non-finite generation of $\gr_{v_t} R$ when $t\in [\frac{7+3\sqrt{5}}{2},+\infty)\setminus \bQ$ is proven in Proposition \ref{prop:S-non-fg}. For $t\in (\frac{7+3\sqrt{5}}{2},+\infty)\cap \bQ$, Theorem \ref{t-complements} and Corollary \ref{cor:alpha<A/T} imply that $\gr_{v_t} R$ is finitely generated whose $\Proj$ is not klt as $A_X(v_t)=T_X(v_t)$. The finite generation of $\gr_{v_t}R$ for $t\in [1, \frac{7+3\sqrt{5}}{2})=\cup_{n\geq 0} [t_n, t_{n+1}]$ follows from Theorem \ref{thm:f.g. criterion} and Proposition \ref{prop:S-fg}. 

Finally, we show that $\Proj \,\gr_{v_t} R$ is a $\bQ$-Fano variety for $t\in [t_n, t_{n+1}]$. When $n=0$, i.e. $t\in [1,2]$, we know from the proof of Proposition \ref{prop:S-fg} that $v_t$ is toric in the projective coordinate $[x_0,x_1-x_2, x_1+x_2]$. Thus $\Proj_{v_t} R  \cong \bP^2$. Thus we may assume $n\geq 1$ in the rest of the proof. From computations in the proof of Proposition \ref{prop:S-fg}, we know that $A_X(E_n)=\epsilon_X(E_n)<T_X(E_n)$ and $A_X(E_n')<\epsilon_X(E_n')=T_X(E_n')$. Since both $\tX_n$ and $\tX_n'$ are of Fano type, we know that $-K_{\tX_n}-E_n$ and $-K_{\tX_n'}-E_n'$ are nef and hence semiample. Thus by Bertini's theorem we can find $\bQ$-complements $G_n$ and $G_n'$ of $(\tX_n, E_n)$ and $(\tX_n',E_n')$ respectively, such that $(\tX_n, E_n+G_n)$ and $(\tX_n', E_n'+G_n')$ are plt. Hence Theorem \ref{thm:zhuang-special} implies that both $E_n$ and $E_n'$ are special divisors as they satisfy $A<T$, and the desired $\bQ$-complements of plt type are given by $(\mu_n)_* G_n$ and $(\mu_n')_* G_n'$. Since $v_{t_n}$ is a rescaling of $\ord_{E_n}$, it induces a special degeneration of $\bP^2$. By the last paragraph of the proof of Lemma \ref{lem:S linear imply f.g. general}, we know that $\gr_{v_t} R\cong \gr_{v_{t_n'}}R\cong \gr_{E_n'} R$ for any $t\in (t_n, t_{n+1})$. Thus $\Proj\,\gr_{v_t} R$ is a $\bQ$-Fano variety as $E_n'$ is special.
\end{proof}

\begin{rem}\label{rem:P^2-deg}
Using similar arguments to \cite[Proof of Proposition 7.4]{ADL19}, one can show that $\Proj\, \gr_{v_{t_n}} R$ for $n\geq 1$ is isomorphic to the weighted hypersurface 
\[
(x_0 x_3 = x_1^{d_{n+1}}+ x_2^{d_{n-1}})\subset \bP(1,d_{n-1}, d_{n+1}, d_n^2).
\]
Such a Manetti surface is a common partial smoothing of $\bP(1, d_{n-1}^2, d_n^2)$ and $\bP(1, d_{n}^2, d_{n+1}^2)$. Similarly, for any $t\in (t_n, t_{n+1})$ one can show that $\Proj\,\gr_{v_t} R\cong\Proj\,\gr_{v_{t_n'}} R\cong \bP(1, d_{n}^2, d_{n+1}^2)$. This provides infinitely many special degenerations of $\bP^2$ which are unbounded.
\end{rem}

\bibliography{ref}

\end{document}